\input amstex
\mag=\magstep1
\documentstyle{amsppt}
\nologo
\NoBlackBoxes
\NoRunningHeads
\vcorrection{-22pt}

\def\fp{\flushpar}
\font\tenptmbit=cmmib10 
\font\eightptmbit=cmmib8
\define\bk#1{\text{\tenptmbit #1}}
\define\eightptbk#1{\text{\eightptmbit{#1}}}
\define\inbox#1{$\boxed{\text{#1}}$}
\define\underbarl#1{\lower 1.4pt \hbox{\underbar{\raise 1.4pt \hbox{#1}}}}
\define\ce#1{\lceil#1\rceil} 
\define\tp#1{\negthinspace\ ^t#1}
\define\dg#1{(d^{\circ}\geqq#1)}
\define\Dg#1#2{(d^{\circ}(#1)\geqq#2)}
\define\llg{\langle\hskip -1.5pt\langle}
\define\rrg{\rangle\hskip -1.5pt\rangle}

\define\jrac#1#2{\dfrac{\lower 2pt\hbox{$#1$}}{\raise 2pt\hbox{$#2$}}}
\define\sz#1{^{\raise 1pt\hbox{$\sssize #1$}}\hskip -1pt}
\define\lr#1{^{\sssize\left(#1\right)}}
\define\br#1{^{\sssize\left[#1\right]}}
\define\do#1{^{\sssize\left<#1\right>}}
\define\ord{\text{\rm ord}}
\define\twolower#1{\hbox{\raise -2pt \hbox{$#1$}}}

\define\hookdownarrow{\kern 0pt
 \hbox{\vbox{\offinterlineskip \kern 0pt
  \hbox{$\cap$}\kern 0pt
  \hbox{\hskip 3.5pt$\downarrow$}\kern 0pt}}
}  
\font\large=cmr10 scaled \magstep2
\def\sumaste{\sideset \and{}^{{\text{\large *}}} \to \sum}

\font\twelveptrm=cmr9 scaled \magstep2

\font\sixptrm=cmr6

\font\thirteenptbf=cmbx13 

\def\boxit#1{\vbox{\hrule\hbox{\vrule\kern3pt
    \vbox to 43pt{\hsize 182pt\kern3pt#1\eject\kern3pt\vfill}
    \kern1pt\vrule}\hrule}} 
\def\oVrule{\vrule width .5pt}
\def\oHrule{\hrule height .5pt}
\def\lbox#1#2#3{\kern 0pt 
\dimen1=#1pt \dimen2=#2pt%
\advance\dimen1 by -1.0pt
\dimen3=\dimen1%
\advance\dimen3 by -16pt
\advance\dimen2 by -1.0pt
\dimen4=\dimen2%
\advance\dimen4 by -20pt
\vbox to #1pt
{
\hsize #2pt
\oHrule
\hbox to #2pt
{
\vsize \dimen1
\oVrule
\vbox to \dimen1
{
\hsize \dimen2
\vskip 8pt
\hbox to \dimen2
{
\vsize \dimen3
\hskip 10pt
%
%
\vbox to \dimen3{\hsize \dimen4\fp #3\vfil}%
\hfil\hskip 10pt
}%
\vskip 8pt
}%
\oVrule%
}%
\oHrule%
}\kern 0pt
}%
\def\cbox#1#2#3{\kern 0pt 
\dimen1=#1pt \dimen2=#2pt%
\advance\dimen1 by -1.0pt
\dimen3=\dimen1%
\advance\dimen3 by -3pt
\advance\dimen2 by -1.0pt
\dimen4=\dimen2%
\advance\dimen4 by -3pt
\vbox to #1pt%
{%
\hsize #2pt
\oHrule%
\hbox to #2pt%
{%
\vsize \dimen1%
\oVrule%
\vbox to \dimen1%
{%
\hsize \dimen2%
\vskip 1.5pt
\hbox to \dimen2%
{%
\vsize \dimen3%
\hskip 1.5pt\hfil
\vbox to \dimen3{\hsize \dimen4\vfil\hbox{#3}\vfil}%
\hfil\hskip 1.5pt
}%
\vskip 1.5pt
}%
\oVrule%
}%
\oHrule%
}\kern 0pt
}%
\def\cboxit#1#2#3{$\hbox{\lower 2.5pt \hbox{\cbox{#1}{#2}{#3}}}$}

\document
\baselineskip=13.3pt 

{\ }
\vskip 100pt
\centerline{\thirteenptbf Theory of The Generalized Bernoulli-Hurwitz Numbers}
\vskip 7pt
\centerline{\thirteenptbf for The Algebraic Functions of Cyclotomic Type}
\vskip 7pt
\centerline{\thirteenptbf and The Universal Bernoulli Numbers}
\vskip 263pt
\centerline{\twelveptrm Yoshihiro \^ONISHI}
\vskip 17pt
\pageno=0
\newpage
\pageno=1
\fp
{\bf Introduction.}\ \
Let  $C$  be a smooth algebraic curve of genus  $g$, 
and assume it is defined by
  $$
  y^2=x^{2g+1}-1, \ \ \ \ \text{or}\ \ \ \   y^2=x^{2g+1}-x 
  \tag{1}
  $$
with the unique point  $\infty$  at infinity. 
We consider the integral
  $$
  u=\int_\infty^x\frac{x^{g-1}dx}{2y}.
  \tag{2}
  $$
As this is an integral of the first kind, namely of a holomorphic $1$-form, 
this is finite valued everywhere.  
Since the integral (2) has a zero at  $\infty$  of order $1$, 
there is an inverse function  $u\mapsto x$  near  $u=0$. 
If  $g = 1$, namely  $C$  is an {\it elliptic curve}, 
then it is meromorphically extended to the whole complex plane of  $u$  and 
such the function  $u\mapsto x$  is no other than 
the elliptic function  $\wp(u)$  attached to  $C$.  

On the other hand, if  $g=2$, the function  $u\mapsto x$  is never extended globally, 
we must think about the set of  $g$  integrals
  $$
  u_1=\int_\infty^x\frac{dx}{2y},\ \ \ 
  u_2=\int_\infty^x\frac{xdx}{2y},\ \ \ 
  \cdots\ , \ 
  u_g=\int_\infty^x\frac{x^{g-1}dx}{2y},
  \tag{3}
  $$
moreover about the sum of the right hand side 
for  $g$  points  $(x_1,y_1)$, $(x_2,y_2)$, $\cdots$, $(x_g,y_g)$, 
and study them as functions of several variables.  
Such the idea is an authority of treatise of Abelian functions back from Jacobi.

Now we should be notice that, in the case of  $g=0$  for the first type in (1), 
the function  $u\mapsto x$  is just the {\it trigonometric function}  $1/\sin^2(u)$.  
The Bernoulli numbers  $B_{2n}$, 
one of the most important set of numbers in mathematics, 
(especially, in the theory of cyclotomic fields) 
is just the Laurent coefficient of this function:
  $$
  \frac{1}{\sin^2(u)}=\frac{1}{u^2}
  +\sum_{n=1}^{\infty}(-1)^{n-1}\frac{2^{2n}B_{2n}}{2n}\frac{u^{2n-2}}{(2n-2)!}.
  \tag{4}
  $$
Furthermore, in the case of  $g=1$  and  $C$  is  $y^2=x^3-x$, 
very important numbers  $E_{4n}$, the Hurwitz numbers, are defined similarly by
  $$
  \wp(u)=\frac1{u^2}
  +\sum_{n=1}^{\infty}\frac{2^{4n}E_{4n}}{4n}\frac{u^{4n-2}}{(4n-2)!}.  
  \tag{5}
  $$
The numbers  $E_{4n}$  are used to study the Gauss number field.  

But Jacobi's theory does not work well for an extension of Bernoulli and Hurwitz numbers, 
and the author were not able to go further for  $g=2$.  
Thus the author simply studied the (2) and he found surprising facts: 
The Laurent coefficients of the function  $u\mapsto x$  have very natural extensions 
of the three main properties for Bernoulli numbers, 
von Staudt-Clausen theorem, von Staudt second theorem, and Kummer congruence relation 
(and their analogies \cite{{\bf Hu1}}, \cite{{\bf Hu2}}, \cite{{\bf L}} for Hurwitz numbers) 
(Theorems 6.1.1, 6.2.1, 7.1.1, and 7.2.1). 

We have here rather modern description of our situation. 
We assume  $g\geqq 1$.  
We construct a standard universal Abelian covering  $\kappa^{-1}\iota(C)$  of  $C$  
in the  $g$  dimensional complex space  $\bold C^g$  as follows. 
The concrete definition is given at the Section 4. 
Here the pull-back of  $\infty$  to  $\bold C^g$  coincides 
with the lattice  $\Lambda\subset\bold C^g$  of all the periods of (3), 
and  $\bold C^g/\Lambda$  is (the  $\bold C$-rational points of) 
the Jacobian variety  $J$  of  $\bold C$.  
\def\kappainviotaC
{\lower 6pt \hbox{\cbox{17}{50}{\hbox{$\kappa^{-1}\iota(C)$}}}}
\def\surjection{\hskip -5pt\longrightarrow\hskip -5pt\rightarrow}
\def\surj#1{\offinterlineskip \vbox{\hbox{\raise 5pt \hbox{$#1$}}\hbox{$\surjection$}}}
  $$
  \matrix
  \kappainviotaC & \surjection   & \ C \hfill \\
   \hbox{\vrule width 0pt height 2pt}    &               &          \\
  \hookdownarrow &               & \ \ \hookdownarrow\raise 5pt\hbox{$\iota$} \hfill\\
  \bold C^g      & \surj{\kappa} & \bold C^g/\Lambda=J
  \endmatrix
  $$
For  $g=1$, the two vertical embeddings are bijections.  
We take the  $g$-th coordinate  $u_g$  as a natural local parameter around 
the origin  $\bold 0=(0,\cdots, 0)\in\bold  C^g$ $=$ 
${(u_1, u_2, \cdots, u_g)}$  of  $\kappa^{-1}\iota(C)$.  
If  $u=(u_1, u_2, \cdots, u_g)$  varies on  $\kappa^{-1}\iota(C)$,  
we denote  $x(\bk{u})$   the inverse function  $u\mapsto x$  determined by (3).  
We define  $y(\bk{u})$  similarly. 
Under this situation, the Laurent coefficients at  $u=0$  with respect to  $u_g$  
of the functions  $x(\bk{u})$  and  $y(\bk{u})$  on  $\kappa^{-1}\iota(C)$,  
which we call generalized Bernoulli-Hurwitz numbers, 
have the quite similar properties as Bernoulli and Hurwitz numbers have. 
Especially we show von Staudt-Clausen type theorem, 
an extension of von Staudt second theorem, 
and Kummer type congruence relation for the generalized Bernoulli-Hurwitz numbers. 
This paper contains various numerical examples 
including the Bernoulli and Hurwitz numbers. 
They would be convenient to the reader.

In order to explain our proofs of the theorems above, 
we introduce the universal Bernoulli numbers
\footnote{July 2003, the author was translating the first Japanese version 
of this paper into English,
and he found a gap in his proof of the Kummer type congruence relation. 
To repair the gap, he met the universal Bernoulli numbers, 
and especially a certain type of Kummer congruence relation
obtained by Professor Arnold Adelberg. }. 
As the Bernoulli numbers are associated to the formal group 
with the formal logarithm  $t\mapsto \log(1+t)$, 
the universal Bernoulli numbers are associated to the universal formal group, 
which is investigated mainly in stable homotopy theory. 
The universal Bernoulli numbers are universal as the name shows, 
and they give our generalized Bernoulli-Hurwitz numbers by a specialization.  
Therefore, roughly speaking, 
all the properties which are satisfied by the universal Bernoulli numbers 
would be also satisfied by the generalized Bernoulli-Hurwitz numbers. 
For instance, Clarke's theorem (Proposition 2.3.1) for the universal Bernoulli numbers,
which is a beautiful unification of analogies of von Staudt-Clausen theorem 
and von Staudt second theorem, 
is essential in our proof of the associated Clarke type theorem 
(Theorems 6.1.1 and 6.2.1) for the generalized Bernoulli-Hurwitz numbers. 

Additionally, our proof of the Clarke type theorem for generalized Bernoulli-Hurwitz numbers 
needs a method that was introduced by Carlitz in \cite{{\bf Ca1}}.  
The explicit addition formula and multiplication formula for an elliptic function are
traditional method to investigate Hurwitz numbers, 
and their natural extensions for  $g>1$  in Jacobi's theory could not work well, 
because the inverse function  $x(\bk{u})$  is out of his theory.  
Surprisingly, Carlitz's method in \cite{{\bf Ca1}} keeps out of such formulae
\footnote{It is plausible by \cite{{\bf Ca1}} and \cite{{\bf Ca4}} 
that Carlitz made much effort to find analogous numbers 
as the Hurwitz numbers for higher genus cases.  
He, however, did not have our approach and could not succeeded. 
The author think that the age when he was active is much poor 
in detailed knowledge about Jacobian varieties. 
Though to find our results needed such knowledge, our paper never uses 
Jacobian varieties at least logically.
}.
We use his method extensively in the Sections 10 and 11. 
His method is a way to compare different kinds of powers of a power series 
which seems to go back to Lagrange inversion theorem
\footnote{This was pointed out by Takayuki Oda. 
We need another type of Lagrange inversion theorem (Proposition 1.2.1) 
in order to have results on the universal Bernoulli numbers.
}.  
Since the divisor of $n$-th power of an algebraic function  $\varphi$  
is just the $n$-plication of the divisor of the function  $\varphi$, 
we can regard his method treats essentially the addition and multiplication formulae
\footnote{This is S.Matustani's opinion.
}.

Furthermore, we show a Kummer type congruence relation (Theorem 3.1.1) for
the universal Bernoulli numbers, which is one of the main results 
of this paper. 
This result gives a best possible estimate about the modulus for such the congruences, 
and implies the Kummer congruence relation (Corollary 3.1.3) discovered by Adelberg. 
Though this result had contained an open lemma (namely, Proposition 3.2.1 (2)), 
Seidai Yasuda proved it February, 2004.  
He also proved the Kummer type congruence relation that had been a conjecture 
in the previous version of this paper.  
He proved the congruence by studying the {\it formal group} obtained by 
the universal formal group by specializing the universal Bernoulli numbers to 
the generalized Bernoulli-Hurwitz numbers. 
At the result, we could not use the congruence for the universal Bernoulli numbers 
in order to prove the one for the generalized Bernoulli-Hurwitz numbers, 
because the former is much weaker than the later. 
So the reader who is interested only in the generalized Brnoulli-Hurwitz numbers 
is suggested to skip the subsections 3.1 3.4 3.5.  

Since to write up the general proof for the widest class 
of the curves in our theory is so complicated, 
we describe in a general form only the statements of main results,
and restrict the context of the paper mainly for the genus two curve
  $$
  y^2=x^5-1.
  $$
Our results also hold for the curves of more general type as in 17.1.  
We call such a curve to be of cyclotomic type.  
The author believes that the reader who understands the arguments 
of this paper can easily apply them generally for the curves of cyclotomic type. 
Moreover, our results would be extended to any curve over the $p$-adic integer ring 
with certain prime  $p$ such that it is completely ramified at infinity 
and the formal completion of its Jacobian variety at the origin contains,
as a factor, the formal completion of the curve at  $\infty$.  
Of course, our proof is effective for the Bernoulli numbers and Hurwitz numbers, 
and is regarded as a new proof for such the classical numbers.

The plan of this research is initialized by a suggestion of M. Koike 
in my visit to Kyushu University several years ago, 
and by an inspiration from the nice book \cite{{\bf AIK}}. 
Professor T. Asai and S. Matsutani gave helpful comments 
and continuously encouraged the author.  
Lemma 18.3.2 was proved by Hiroshi Suzuki.  
Professor A. Adelberg read the draft of the Sections from 1 to 3, 
and he pointed out several errors. 
These three Sections was established by his helpful comments. 
I would like to express my thanks to the related persons.

This paper is a complete revised edition by including Yasuda's proofs, 
of an old manuscript which contained certain open lemma about Kummer type congruence
for the universal Bernoulli numbers, and in which the Kummer type congruence
for the generalized Bernoulli-Hurwitz numbers was a conjecture yet. 
Though this edition contains some personal opinions of the author and it was named solely, 
the reader should understand this is a joint work of Yasuda and \^O.  

The author understand by this work that Hurwitz is deeply related to the Leonid meteor storms
\footnote{Around the year Hurwitz's papers \cite{{\bf Hu1}} and \cite{{\bf Hu2}} were published, 
the mother comet Tempel-Tuttle of the Leonid meteors got across the Earth's orbit 
(the date would be in June, 1898). 
The day of Hurwitz's death, 18th December 1919, is a day 
peoples were able to watch the the Leonid storms. 
The author found later that the initial phenomenon leading the results of this paper 
was discovered just at the peak of the meteors at Japan, 
very early before day-light of 19th December 2002. 
The author strongly hopes to watch spectacular display from the Leonid meteors in future.
}.

\newpage \fp
{\bf Contents}
\fp
\vskip 7pt
\fp
\inbox{1} Preliminaries from combinatorics \fp\ \ \ 
  1.1. Fundamentals on factorials\fp\ \ \ 
  1.2. Lagrange inversion formula\fp\ \ \ 
  1.3. Fundamentals on binomial coefficients
\fp
\vskip 1.5pt
\fp
\inbox{2} The universal Bernoulli numbers and their properties \fp\ \ \ 
  2.1. Definition of the universal Bernoulli numbers\fp\ \ \ 
  2.2. Schur function type expression of the universal Bernoulli numbers\fp\ \ \ 
  2.3. Clarke's theorem
\fp
\vskip 1.5pt
\fp
\inbox{3} The Kummer type congruence for the universal Bernoulli numbers \fp\ \ \ 
  3.1. Main theorem\fp\ \ \ 
  3.2. The first tool for the proof (Yasuda)\fp\ \ \ 
  3.3. The second tool for the proof\fp\ \ \ 
  3.4. Proof of the Kummer type congruence relation\fp\ \ \ 
  3.5. The Kummer-Adelberg type congruence relation
\fp
\vskip 1.5pt
\fp
\inbox{4} Hyperelliptic functions \fp\ \ \ 
  4.1. Fundamentals\fp\ \ \ 
  4.2. Hyperelliptic functions and their variable
\fp
\vskip 1.5pt
\fp
\inbox{5} Differential equations \fp\ \ \ 
  5.1. Generalities\fp\ \ \ 
  5.2. The case of  $y(\bk{u})^2=x(\bk{u})^{2g+1}-1$\fp\ \ \ 
  5.3. The case of  $y(\bk{u})^2=x(\bk{u})^{2g+1}-x(\bk{u})$
\fp
\vskip 1.5pt
\fp
\inbox{6} The Clarke type theorem \fp\ \ \ 
  6.1. The case of  $y(\bk{u})^2=x(\bk{u})^{2g+1}-1$ \fp\ \ \ 
  6.2. The case of  $y(\bk{u})^2=x(\bk{u})^{2g+1}-x(\bk{u})$ \fp\ \ \ 
  6.3. A remark on  $A_p$
\fp
\vskip 1.5pt
\fp
\inbox{7} The Kummer type congruence relation \fp\ \ \ 
  7.1. The case of  $y(\bk{u})^2=x(\bk{u})^{2g+1}-1$ \fp\ \ \ 
  7.2. The case of  $y(\bk{u})^2=x(\bk{u})^{2g+1}-x(\bk{u})$ 
\fp
\vskip 1.5pt
\fp
\inbox{8} Hurwitz-integral series \fp\ \ \ 
  8.1. Definition and basic properties \fp\ \ \ 
  8.2. Hurwitz-integrality of  $1/x^{1/2}(\bk{u})$ \fp\ \ \ 
  8.3. Hurwitz-integrality of  $1/y^{1/5}(\bk{u})$ 
\fp
\vskip 1.5pt
\fp
\inbox{9} Outline of the proofs \fp\ \ \ 
  9.1. The Clarke type theorem\fp\ \ \ 
  9.2. The Kummer type congruence relation
\fp
\vskip 1.5pt
\fp
\inbox{10} The Clarke type theorem on  $x(\bk{u})$  over  $\bold Z[\frac12]$ \fp\ \ \ 
  10.1. The Clarke type theorem on  $x(\bk{u})^{1/2}$ \fp\ \ \ 
  10.2. Relation between coefficients of  $x(\bk{u})^{1/2}$  and of  $x(\bk{u})$ \fp\ \ \ 
  10.3. Relation between coefficients of  $x(\bk{u})$  and of  $x(\bk{u})^{3/2}$ \fp\ \ \ 
  10.4. Relation between coefficients of  $x(\bk{u})^{3/2}$   and of $x(\bk{u})^2$ \fp\ \ \ 
  10.5. The Clarke type theorem on  $x(\bk{u})^2$  over  $\bold Z[\frac12]$
\fp
\vskip 1.5pt
\fp
\inbox{11} The Clarke type theorem on  $y(\bk{u})$  over  $\bold Z[\frac15]$ \fp\ \ \ 
  11.1. The Clarke type theorem on  $y(\bk{u})^{1/5}$ \fp\ \ \ 
  11.2. Relation between coefficients of  $y(\bk{u})^{1/5}$  and of  $y(\bk{u})^{2/5}$ \fp\ \ \ 
  11.3. Relation between coefficients of  $y(\bk{u})^{2/5}$  and of  $y(\bk{u})^{3/5}$ \fp\ \ \ 
  11.4. Relation between coefficients of  $y(\bk{u})^{3/5}$  and of  $y(\bk{u})^{4/5}$ \fp\ \ \ 
  11.5. Relation between coefficients of  $y(\bk{u})^{4/5}$  and of  $y(\bk{u})$ \fp\ \ \ 
  11.6. The Clarke type theorem on  $y(\bk{u})$  over  $\bold Z[\frac15]$
\fp
\vskip 1.5pt
\fp
\inbox{12} The Clarke type theorem on  $x(\bk{u})$ and  $y(\bk{u})$  over  $\bold Z$ \fp\ \ \ 
  12.1. Relation between coefficients of  $y(\bk{u})$  and of  $x^2(\bk{u})$ \fp\ \ \ 
  12.2. From  $x^2(\bk{u})$  to  $x(\bk{u})$
\fp
\vskip 1.5pt
\fp
\inbox{13} Proof of the Kummer type congruence relation (Yasuda) \fp\ \ \ 
  13.1. Honda's theorem and formal groups  \fp\ \ \ 
  13.2. Hochschilt's formula and Honda's theorem \fp\ \ \ 
  13.3. Proof of the Kummer type congruence relation
\fp
\vskip 1.5pt
\fp
\inbox{14} Other Kummer type congruence relations \fp\ \ \ 
  14.1. On generalized Bernoulli-Hurwitz numbers of higher order \fp\ \ \ 
  14.2. On the Hurwitz coefficients of  $t(u_g)$  and of  $s(u_g)$ \fp\ \ \ 
  14.3. The Vandiver-Carlitz type congruence relation
\fp
\vskip 1.5pt
\fp
\inbox{15} Numerical examples for classical numbers \fp\ \ \ 
  15.1. Bernoulli numbers \fp\ \ \
  15.2. Hurwitz numbers for the curve  $y^2=x^3-1$ \fp\ \ \ 
  15.3. Hurwitz numbers for the curve  $y^2=x^3-x$
\fp
\vskip 1.5pt
\fp
\inbox{16} Numerical examples for new numbers \fp\ \ \ 
  16.1. $x(\bk{u})$  of the curve  $y^2=x^5-1$ \fp\ \ \ 
  16.2. $y(\bk{u})$  of the curve  $y^2=x^5-1$ \fp\ \ \ 
  16.3. $x(\bk{u})$  of the curve  $y^2=x^5-x$ \fp\ \ \ 
  16.4. $x(\bk{u})$  of the curve  $y^2=x^7-1$ \fp\ \ \ 
  16.5. $x(\bk{u})$  of the curve  $y^2=x^7-x$ \fp\ \ \ 
  16.6. $x(\bk{u})$  of the curve  $y^2=x^9-1$
\fp
\vskip 1.5pt
\fp
\inbox{17} Non-hyperelliptic curves \fp\ \ \ 
  17.1. Algebraic curves with completely ramified at infinity\fp\ \ \ 
  17.2. $x(\bk{u})$  of the curve  $y^3=x^5-1$ \fp\ \ \ 
  17.3. $x(\bk{u})$  of the curve  $y^4=x^3-x$
\fp
\vskip 1.5pt
\fp
\inbox{18} Appendices \fp\ \ \ 
  18.1. Relations on binomial coefficients\fp\ \ \ 
  18.2. Links with certain Eisenstein type series without justification\fp\ \ \ 
  18.3. Related formula with Hochschilt's formula\fp\ \ \ 
  18.4. Problems (or exercises?)
\fp
\vskip 1.5pt
\fp
References

\newpage \fp
{\bf Convention on notations. }\ \ \ \fp
\vskip 5pt
\fp
(1) For a rational number  $\alpha$, 
we denote by  $\lfloor \alpha \rfloor$  the largest integer not exceed  $\alpha$,
and by  $\lceil \alpha \rceil$  the smallest integer not smaller than  $\alpha$.  
\vskip 5pt
\fp
(2) We use the notation
  $$
  (z)_n=z(z-1)\cdots(z-n+1)
  $$
for an integer n = 0.  
Here the range of z is determined by each situation. 
The binomial coefficiets is written by this as
  $$
  \binom{z}{n}=\frac{(z)_n}{n!}.
  $$
\vskip 5pt
\fp
(3) Generalizing the factorial symbol, 
we denote by  $n!!$  the product  $n(n-2)(n-4)\cdots$  
of the sequence with step  $-2$  from  $n$  to  $1$ or $2$. 
Similarily  $n!!!$  means the product of the sequence of positive integers 
from  $n$  with step  $-3$.  
Moreover, for instance, we denote  $n!!!!!=n!\lr{5}$.  
Hence  $12!\lr{5}=12\cdot 7\cdot 2$.  
\vskip 5pt
\fp
(4) If  $p$  is a prime and the $p$-part of given rational number  $r$  is  $p^e$, 
then we write  $e=\ord_pr$.  
If  $\tau$  is a polynomial (possibly in several variables) with rational coefficients,
then we denote by  $\ord_p\tau$  the least number of ordpr for all the coefficients  $r$  of  $\tau$.  
\vskip 5pt
\fp
(5) For a prime number  $p$  and an integer  $a$, 
we denote by  $a|_p$  the neglect of  $p$-part of  $a$, 
namely  $a|_p=a/p^{\ord_pa}$.  
\vskip 5pt
\fp
(6) If  $F(z)$  is a formal power seriesnwith respect to  $z$,  
we denote by  $[z^n]F(z)$  its coefficient of  $z^n$. 
We use also  $[\frac{z^n}{n!}]F(z):=n![z^n]F(z)$  
or  $[\frac{z}{n}]F(z):=n[z^n]F(z)$.  
\vskip 5pt
\fp
(7) For a formal power series  $\varphi(z)$  with respect to  $z$ 
(permitting negative-power terms), 
we call  $[\frac{z^n}{n!}]\varphi(z)$ ($n\geqq 0$) the Hurwitz coefficient 
for  $\varphi(z)$  of  $z^n$.  
We say  $[\frac{z}{n}]\varphi(z)$  its Carlitz coefficient.
\vskip 5pt
\fp
(8) Let  $R$  be a commutative ring.  
We denoted by  $R\llg v \rrg$  the ring of the formal power series 
with respect to  $z$  consisting of only terms of non-negative order 
such that all of whose Hurwitz coefficients belong to  $R$.
\vskip 5pt
\fp
(9) In an expression of a power series with respect to  $z$, 
the symbol  $(d^{\circ}(z)\geqq m)$  stands for 
its part of the terms of degree at least  $m$.  
When  $m$  is obvious, we simply denote them as usual by  $+\cdots$.

\newpage \input 01-1.o
         \fp
{\bf 1.2.   Lagrange inversion formula.}\ \
For a power series  $F(z)$  with respect to  $z$, 
we denote by  $[z^n]F(z)$  its coefficient of  $z^n$.  
The following formula is called Lagrange inversion formula.  
\vskip 5pt
\fp
\lbox{95}{360}
{{\bf Proposition 1.2.1.}\ \ \it
Let  $\varphi(u)=u+\cdots$  be a power series with respect to  $u$ 
consists of only the terms of non-negative degree such that its coefficient of  $u$  
is  $1$.  
Let  $\psi(t)=\varphi^{-1}(t)$   be its formal inverse series, 
namely the power series with respect to  $t$  such that  $\varphi(\psi(t))=t$.  
Then 
\fp
  $$
  [u^n]\Big(\frac{u}{\varphi(u)}\Big)^{\ell}
  =\frac{\ell}{\ell-n}[t^n]\Big(\frac{\psi(t)}{t}\Big)^{\ell-n}.
  $$
}
\vskip 2pt
\fp
Proofs are found in  \cite{{\bf Co}, pp.148-153} 
and the reference for \cite{{\bf Ad1}, Proposition 2.1}.  
\vskip 12pt

         \fp
{\bf 1.3. Fundamentals on binomial coefficients.}\ \ 
We prove the following Lemma.  
We do not use this until 10.1.  
\vskip 5pt
\fp
\lbox{67}{360}
{{\bf Lemma 1.3.1.}\ \ \it
Let  $n$  and  $q$  are positive integers, 
and  $r$  be an integer such that  $r<q$.  
Then
  $$
  \frac{(qn-r)!\lr{q}}{(qn)!\lr{q}}
  \in\bold Z\Big[\jrac 1q\Big].
  $$
}
\vskip 5pt
\fp
{\it Proof}. \ 
Let  $\xi=q\{(1+z)^{1/q}-1\}$. 
So that  $z=(1+\frac 1q \xi)^q-1$.  
We apply the Lagrange inversion formula (1.2.3) for  $\ell=-1$, 
and have
  $$
  [z^n]\bigg(\frac{z}{q\{(1+z)^{1/q}-1\}}\bigg)^{\hskip -2pt -1}
  \hskip -4pt 
  =[\xi^n]\bigg(\frac{(1+\frac 1q \xi)^q-1}{\xi}\bigg)^{\hskip -3pt -1-n}
  \hskip -3pt \tfrac{d}{d\xi}\{(1+\tfrac 1q \xi)^q-1\}. 
  $$
The right hand side obviously belongs to  $\bold Z[\frac{1}{q}]$.  
Hence the left hand side too.  
Looking at the coefficients of
  $$
  (1+z)^{r/q}=\{(1+z)^{1/q}\}^r,
  $$
we see the statement.  
\qed
\vskip 5pt
Here, we show this by another way.  
Although the following Lemma is rather weaker than 1.3.1, 
we have 1.3.1 by varying  $p$  because of the equality
  $$
  \binom{\frac{r}{q}-1}{n}
  =\Big(\frac{r- q}{q}\Big)
   \Big(\frac{r-2q}{q}\Big)
   \cdots
   \Big(\frac{r-nq}{q}\Big)
   \Big/n!
  =(-1)^n\frac{(qn-r)!\lr{q}}{(qn)!\lr{q}}.
  \tag{1.3.2}
  $$
\vskip 3pt
\fp
\lbox{66}{360}{{\bf Lemma 1.3.3.}\ \ \it
Let  $p$  be a prime number,  $n\geqq 0$  be an integer, and  $z\in \bold Z_{(p)}$.  
Then 
  $$
  \binom{z}{n}\in \bold Z_{(p)}. 
  $$
}
\vskip 5pt
\fp
{\it Proof }. \ 
For a given integer  $n\geqq 0$, 
the function  $z\mapsto n$  is a continuous map 
from  $\bold Z_p$  to  $\bold Q_p$.  
Since  $\bold Z$  is dense in  $\bold Z_p$  
and  $\binom{z}{n}\in\bold Z$  for  $z\in\bold Z$, 
the statement follows.
\qed
\vskip 35pt

\newpage \fp
\inbox{\bf 2}\ \ 
{\bf The universal Bernoulli numbers and their properties.}\ \ 
\vskip 4pt
\fp
{\bf 2.1. Definition of the universal Bernoulli numbers.}\ \  
Let  $f_1$, $f_2$, $\cdots$  be infinitely many indeterminates.  
We consider the power series 
  $$
  u=u(t)=t+\sum_{n=1}^{\infty}f_{\sssize n}\frac{t^{n+1}}{n+1}
  \tag{2.1.1}
  $$
and its formal inverse series
  $$
  t=t(u)=u-f_1\frac{u^2}{2!}+(3{f_1}^2-2f_2)\frac{u^3}{3!}+\cdots,
  \tag{2.1.2}
  $$
namely such the series that  $u(t(u))=u$.  
Then we define  $\hat{B}_{n2}\in\bold  Q[f_1, f_2, \cdots]$.  
by
  $$
  \frac{u}{t(u)}
  =\sum_{n=0}^{\infty}\hat{B}_n\frac{u^n}{n!}
  \tag{2.1.2}
  $$
and call them the universal Bernoulli numbers (of order 1). 
If we set  $f_n=(-1)^n$, then  $u(t)=\log(1 + t)$  
and  $t(u)=e^t-1$.   
Then  $\hat{B}_n$  is just the classical Bernoulli number  $B_n$.  
\vskip 15pt

         \fp
{\bf 2.2. Schur function type expression of the universal Bernoulli numbers. }\ \     
We introduce more symbols.  
For a finite sequence  $U=(U_1, U_2, \cdots)$  of non-negative integers, 
we denote  $w(U)=\sum_jjU_j$, and call it the weight of  $U$.  
The number  $d(U)=\sum_jU_j$  is called the degree of  $U$.   
We can regard  $U$  a partition of  $w(U)$.   
For simplicity, we write as 
  $$
  U!=U_1!U_2!\cdots, \ \ \ \binom{d}{U}=\frac{d!}{U!}.
  \tag{2.2.1}
  $$
Using notations  $\Lambda^{\sssize U}=2^{\sssize U_1}3^{\sssize U_2}4^{\sssize U_3}\cdots$,
and  $f^{\sssize U}={f_1}^{\sssize U_1}{f_2}^{\sssize U_2}{f_3}^{\sssize U_3}\cdots$, 
we define
  $$
  \gamma_{\sssize U}
  =\Lambda^{\sssize U}U!
  \tag{2.2.2}
  $$
If we set  $h(t)=(\psi(t)/t)-1$, then
  $$
  \Big(\frac{\psi(t)}{t}\Big)^s
  =\big(1+h(t)\big)^s
  =\sum_{d=0}^{\infty}\binom{s}{d}h^d(t), \hskip 10pt
  h^d(t)=\sum_{d(U)=d}\binom{d}{U}
  \frac{f^{\sssize U}}{\Lambda^{\sssize U}}\hskip 2pt t^{\sssize w(U)}. 
  \tag{2.2.3}
  $$
Hence, by writing as
  $$
  \tau_{\sssize U}=(-1)^{d(U)-1}\frac{(w(U)+d(U)-2)!}{\gamma_{\sssize U}}
  \tag{2.2.4}
  $$
and using 1.2.1 for  $\ell=1$, we have the following expression of  $\hat{B}_n$.  
\vskip 5pt
\fp
\lbox{65}{360}
{{\bf Proposition 2.2.5.}\ \ \it We have
  $$
  \frac{\hat{B}_n}{n}
   =\hskip -2pt\sum_{w(U)=n}
   \tau_{\sssize U}f^{\sssize U}.
  $$
}
\vskip 5pt
\fp
Following Haigh (\cite{{\bf Cl}, p.594, $\ell$.-12}), 
we call this the Schur function type expression of  $\hat{B}_n$.

\newpage \fp
{\bf 2.3. Clarke's theorem.}\ \ \ 
Here we describe Clarke's theorem on the universal Bernoulli numbers.  
For a prime number  $p$  and an integer  $a$, 
we denote by  $a|_p$  the neglect of  $p$-part of  $a$, 
namely  $a|_p=a/p^{\ord_pa}$.  
\vskip 5pt
\lbox{224}{360}
{{\bf Proposition 2.3.1.}\ \ \it We have
  $$
  \align
  \hat{B}_1&=\frac12f_1, \\
  \frac{\hat{B}_2}{2}&=-\frac14{f_1}^2+\frac13f_2, \\
  \frac{\hat{B}_n}{n}
  &\equiv 
  \cases
  \dsize\sum_{\Sb n=a(p-1) \\ p:\text{\sixptrm\hskip 2pt prime}\endSb}\hskip -10pt 
  \dfrac{{a|_{\sssize p}}^{-1}\ \text{\rm mod}\ p^{1+\ord_pa}}{p^{1+\ord_pa}}\ 
  {f_{p-1}}^a
      \hskip 30pt \text{\rm ( if  $n\equiv 0 \mod 4$ )} & \\
  \dfrac{{f_1}^{n-6}{f_3}^2}{2}-\dfrac{n{f_1}^n}{8}
  + \hskip -10pt
  \dsize\sum_{\Sb n=a(p-1) \\ p:\text{\sixptrm\hskip 2pt odd prime}\endSb}\hskip -10pt 
  \dfrac{{a|_{\sssize p}}^{-1}\ \text{\rm mod}\ p^{1+\ord_pa}}{p^{1+\ord_pa}}\ 
  {f_{p-1}}^a & \\
      \hskip 124pt 
      \text{\rm ( if  $n\neq 2$  and  $n\equiv 2 \mod 4$ )} & \\
  \dfrac{{f_1}^n+{f_1}^{n-3}f_3}{2} \hskip 50pt
      \text{\rm ( if $n\neq 1$  and  $n\equiv 1$, $3 \mod 4$ )} & 
  \endcases  \\
  &\mod \text{\lq\lq\hskip 2pt the set of weight  $n$  polynomials in  $\bold Z[f_1, f_2, \cdots]$"}. 
  \endalign
  $$
}
\fp
\vskip 5pt
\fp
The proof is elementary, but is given by quite involved use of 2.2.5. 
For the reader who is interested in should be referred to \cite{{\bf Cl}}.  
\vskip 5pt
\fp
{\bf Remark 2.3.2.}\ 
(1) For example, if  $n\equiv 0$ mod $4$, 
this proposition shows that
  $$
  \hat{B}_n\equiv -\sum_{\Sb p:\text{\sixptrm prime}\\ p-1|n\endSb}
   \frac{{f_{p-1}}^{n/(p-1)}}{p}
  \mod{\bold Z[f_1, f_2, \cdots]}, 
  $$
and if additionally  $p-1\hskip -4pt\not\hskip -1pt|\hskip 2pt n$  
then  $\hat{B}_n/n\in\bold Z_{(p)}[f_1,f_2,\cdots]$. 
Indeed for a prime  $p$, if  $n=a(p-1)$  and   $\ord_pa=\nu$, then
  $$
  \frac{n\cdot\big({a|_p}^{-1}\hskip -5pt\mod{p^{1+\nu}}\big)}{p^{1+\nu}}
  \equiv 1-\frac{1}{p} \mod{p^{\nu}\bold Z}
  $$
and
  $$
  \frac{n\cdot\big({a|_p}^{-1}\hskip -5pt\mod{p^{1+\nu}}\big)}{p^{1+\nu}}
  \equiv -\frac{1}{p} \mod{\bold Z}.
  $$
Hence the former statement.  
The later is obvious.  
The former is an analogy of the von Staudt-Clausen theorem for the classical Bernoulli numbers, 
and the later is an analogy of the von Staudt second theorem.  
Clarke's theorem is a beautiful unification of these analogy of these two theorems. 
\fp
(2) Although the modulus is whole the ring  $\bold Z[f_1, f_2, \cdots]$  in  \cite{{\bf Cl}},  
we obviously see that it may be replaced as above because of 2.2.5.

\newpage \fp
\inbox{\bf 3}\ \ 
{\bf The Kummer type congruences for the universal Bernoulli numbers.}\ \ 
\vskip 5pt
\fp
{\bf 3.1. Main theorem.}\ \ \ 
The universal Bernoulli numbers satisfy the congruence
relation of Kummer's original type modulo  $p^{\lfloor a/2 \rfloor}$.  
\vskip 5pt
\fp
\lbox{77}{360}
{{\bf Theorem 3.1.1.}\ \ \it
Let  $p$  be a prime, $a$  and  $n$  be positive integers. 
Assume that  $n>a$  and  $n\not\equiv 0$  {\rm mod}  $(p-1)$. 
Then
  $$
  \sum_{r=0}^a\binom{a}{r}(-f_{p-1})^{a-r}
           \frac{\hat{B}_{n+r(p-1)}}{n+r(p-1)} 
  \equiv 0 \mod \ p^{\lfloor a/2 \rfloor}\bold Z_{(p)}[f_1,f_2,\cdots].
  $$
}
\vskip 5pt
\fp
{\bf Remark 3.1.2.}\ \ 
(1) If  $a = 1$  and  $n>a=1$  with  $n\not\equiv 0$, $1$ mod $(p-1)$, 
then the congruence above holds modulo  $p$.  
This fact appeared by \cite{{\bf Ad1}, Theorem 3.2} in the first time. 
A shorter proof and an extension to the case where  $n\equiv 1$ mod $(p-1)$  
are given in  \cite{{\bf Ad2} ,Theorem 1}.
\fp
(2) For an odd prime  $p\geqq 7$, let consider  $U$  such that  $U_1=p$, $U_{2p-1}=(p - 3)/2$,  
and the other entries are  $U_j=0$.  
Then  $w(U)=p+(p-5)(2p-1)/2\equiv -1$ mod $(p-1)$.  
For this  $U$, we can easily show that
  $$
  \ord_p(\tau_{\sssize U})=(p-5)/2 (=\lfloor (p-4)/2 \rfloor).
  $$
Note that we have  $n>a$  for  $a=p-4$  and  $n=w(U)$.   
Looking at (3.4.1) below, we understand that the estimate 3.1.1 is best possible.  
\fp
(3) In the example above for  $p=5$, we have  $\ord_5(\tau_{\sssize U})=0$. 
Then  $n=w(U)=5\equiv 1$ mod $(5 - 1)$ 
and this is one of the excluded case in \cite{{\bf Ad1}, Theorem 3.2}.   
Keeping this case in mind and slightly modifying the proof of 3.3.1 below, 
we can show that the congruence in 3.1.1 holds modulo  $p$ for  $a=1$  
provided the additional condition  $n\not\equiv 1$ mod $(p - 1)$.  

We prove in 3.5 the following congruence of Adelberg 
(\cite{{\bf Ad3}, (i) of the Theorem}) directly from 3.1.1 and the Remark 3.1.2(1).   
\vskip 5pt
\fp
\lbox{69}{360}
{{\bf Corollary 3.1.3.}\ {\rm (Adelberg's congruence)} \it   
If  $n\not\equiv 0$, $1$  {\rm mod}  $(p-1)$  and  $n>a$, then
  $$
  {f_{p-1}}^{p^{a-1}}\hskip -3pt\cdot
  \frac{\hat{B}_{n\vphantom{n+p^a(p-1)}}}{n}
  \equiv  
  \frac{\hat{B}_{n+p^{a-1}(p-1)}}{n+p^{a-1}(p-1)}
  \mod p^a\bold Z_{(p)}[f_1,f_2,\cdots].
  $$
}
\vskip 3pt
\fp
The Theorem 3.1.1 gives a similar but more complicated
congruence if  $n\equiv 1$ mod $(p-1)$, that is  \cite{{\bf Ad3}, (ii) of the Theorem}.

\newpage \fp
{\bf 3.2. The first tool for the proof (Yasuda).}\ \ \ 
\fp
We need two Propositions for the proof of 3.1.1. 
The following is the first one. 
\fp
\vskip 3pt
\fp
\lbox{168}{360}
{{\bf Proposition 3.2.1.}\ \ \it
Let  $p$  be a prime, and  $a$, $n$  be non-negative integers.  
Let 
  $$
  M=\cases
    \ord_p(n!) \ \ & \text{{\rm (} if  $n\geqq ap$ {\rm )}}, \\
           a-\lfloor n/p \rfloor + \ord_p(n!)
          -  \lfloor(a-\lfloor n/p \rfloor)/p\rfloor 
                          & \text{{\rm (} if  $n<ap$ {\rm )}}.
    \endcases
  $$
Then we have the following congruences.  
\fp
{\rm (1)} If  $q$ is a non-negative integer, then
  $$
  \sum_{r=0}^a\frac{\left((r+q)p+n\right)!}{(r+q)!p^{r+q}}
  \binom{a}{r}
  \equiv 0 \mod p^M;
  $$
\fp
{\rm (2)} If  $r_0$  is an integer such that  $0<r_0\leqq a$  and  $n\geqq r_0p$  is an integer, 
then
  $$
  \sum_{r=r_0}^a\frac{\left((r-r_0)p+n\right)!}{(r-r_0)!p^{r-r_0}}
  \binom{a}{r}
  \equiv 0 \mod p^M. 
  $$
}
\vskip 2pt
\fp
At the beginning of this work, while (1) had been proved by the author, 
(2) was still a conjecture.  
Recently Yasuda gave a unified proof of (1) and (2) in \cite{{\bf Ya1}}.  
This subsection is entirely due to his proof.  For the proof we define, 
for two real numbers  $\alpha>0$  and  $\beta$,
  $$
  \bold Z_p\llg v \rrg^{(\alpha),\hskip 1pt\beta}
  =\Big\{
    \varphi \in \bold Q_p[[v]]\hskip 2pt \Big|
    \hskip 2pt  \big(\tfrac{d}{dv}\big)^{pa}\varphi
    \in p^{\lceil \alpha(a+\beta)\rceil}\bold Z_p\llg v\rrg
    \ \text{for all}\  a\geqq 0
   \Big\}. 
  \tag{3.2.2}
  $$
\vskip 3pt
\fp
{\bf Remark 3.2.3.}\ 
(1) The reader may be easier to understand the proof 
if he regard  $\bold Z_p\llg v \rrg^{(\alpha),\hskip 1pt \beta}$  
as if  \lq\lq$p^{\alpha\beta}\bold Z_p\llg p^{\alpha/p}v\rrg$".  
\fp
(2) We only use the case of  $\alpha=1-\frac{1}{p}$  for the proof of 3.2.1.  
\vskip 5pt
\fp
\lbox{164}{360}{{\bf Lemma 3.2.4.}\ \it 
{\rm (1)} $\bold Z_p\llg v\rrg^{(\alpha),\hskip 1pt\beta}$  is a  $\bold Z_p$-submodule in  $\bold Q_p[[v]]$.  
\fp
{\rm (2)} $\bold Z_p\llg v\rrg^{(\alpha),\hskip 1pt\beta}
       \subset 
       p^{\lceil \alpha\beta\rceil}\bold Z_p\llg v\rrg$, 
Especially, if  $\beta\geqq 0$  
then  $\bold Z_p\llg v\rrg^{(\alpha),\hskip 1pt\beta}\subset \bold Z_p\llg v\rrg$. 
\fp
{\rm (3)} Let  $f(v)=\sum_{n=0}^{\infty}a_nv^n\in \bold Q_p[[v]]$.  
$f(v)\in \bold Z_p\llg v\rrg^{(\alpha),\hskip 1pt \beta}$  
if and only if  $\sum_{n=0}^ma_nv^n\in \bold Z_p\llg v\rrg^{(\alpha),\hskip 1pt \beta}$  
for all  $m\geqq 0$.  
\fp
{\rm (4)} If  $\alpha_1\geqq \alpha_2>0$  
then  $\bold Z_p\llg v\rrg^{(\alpha_1),\hskip 1pt\beta}\subset \bold Z_p\llg v\rrg^{(\alpha_2),\hskip 1pt\beta}$. 
\fp
{\rm (5)} If  $\beta_1\geqq \beta_2$  
then  $\bold Z_p\llg v\rrg^{(\alpha),\hskip 1pt\beta_1}\subset \bold Z_p\llg v\rrg^{(\alpha),\hskip 1pt\beta_2}$. 
\fp
{\rm (6)} If  $n\in \bold Z$ (possibly negative) 
then  $p^n\bold Z_p\llg v\rrg^{(\alpha),\hskip 1pt\beta}
\subset \bold Z_p\llg v\rrg^{(\alpha),\hskip 1pt\beta+\frac{n}{\alpha}}$. 
\fp
{\rm (7)} If  $\varphi\in \bold Z_p\llg v\rrg^{(\alpha),\hskip 1pt\beta}$  
then  $\frac{d}{dv}\varphi\in \bold Z_p\llg v\rrg^{(\alpha),\hskip 1pt\beta}$. 
\fp
{\rm (8)} If  $\varphi\in \bold Z_p\llg v\rrg^{(\alpha),\hskip 1pt\beta}$  
then  $\big(\frac{d}{dv}\big)^p\varphi\in \bold Z_p\llg v\rrg^{(\alpha),\hskip 1pt\beta+1}$. 
\fp
{\rm (9)} Let  $m\geqq 0$  and  $b\in \bold Q_p$. 
$bv^m\hskip -2pt/m! \in \bold Z_p\llg v\rrg^{(\alpha),\hskip 1pt\beta}$  
if and only if  
$\ord_p(b)\geqq \alpha\big(\lfloor \frac{m}{p}\rfloor+\beta\big)$.   
}
\fp
{\it Proof.} \  
All the statements are easily shown by the definition (3.2.2).  
\qed

\newpage

\vskip 3pt
\fp
\lbox{56}{360}{{\bf Corollary 3.2.5.}\ \it  
If  $0<\alpha\leqq 1$, $\varphi\in\bold Z_p\llg v\rrg^{(\alpha),\hskip 1pt\beta}$  then 
  $$
  \Big(\Big(\jrac{d}{dv}-1\Big)^p+1\Big)\varphi
  \in \bold Z_p\llg v\rrg^{(\alpha),\hskip 1pt\beta+1}.
  $$
}
\vskip 3pt
\fp
{\it Proof.} \ 
For our purpose it suffice only to prove the case  $p\neq 2$.  
The Leibniz rule shows
  $$ 
  \Big(\Big(\jrac{d}{dv}-1\Big)^p+1\Big)\varphi
  =\Big(\jrac{d}{dv}\Big)^p\varphi
   +\sum_{j=1}^{p-1}(-1)^j\binom{p}{j}\Big(\jrac{d}{dv}\Big)^j\varphi. 
  $$
This belongs to  $\bold Z_p\llg v\rrg^{(\alpha),\hskip 1pt\beta+1}+p\bold Z_p\llg v\rrg^{(\alpha),\hskip 1pt\beta}$  
Therefore the statement follows from 3.2.4 (6), (5) and the assumption for  $\alpha$.  
\qed
\vskip 7pt
\fp
\lbox{81}{360}{{\bf Lemma 3.2.6.}\ \ \it
Let  $0<\alpha\leqq 1$.  \fp
{\rm (1)} $\bold Z_p[[v]] \subset \bold Z_p\llg v\rrg^{(\alpha),\hskip 1pt 0}$.  \fp
{\rm (2)} If  $\varphi_1\in\bold Z_p\llg v\rrg^{(\alpha),\hskip 1pt \beta_1}$  
and     $\varphi_2\in\bold Z_p\llg v\rrg^{(\alpha),\hskip 1pt \beta_2}$, 
then  $\varphi_1\varphi_2\in\bold Z_p\llg v\rrg^{(\alpha),\hskip 1pt \beta_1+\beta_2}$. \fp
Particularly, if  $0<\alpha\leqq 1$, 
then  $\bold Z_p\llg v\rrg^{(\alpha),\hskip 1pt 0}$  is 
a sub $\bold Z_p$-algebra in  $\bold Z_p\llg v\rrg$, 
and  $\bold Z_p\llg v\rrg^{(\alpha),\hskip 1pt \beta}$  is 
a  $\bold Z_p\llg v\rrg^{(\alpha),\hskip 1pt 0}$-module.  
}
\vskip 3pt
\fp
{\it Proof.} \ 
(1) Suppose  $\varphi\in \bold Z_p[[v]]$. 
Since  $\frac{1}{(pa)!}\big(\frac{d}{dv}\big)^{pa}\varphi\in \bold Z_p[[v]]\subset \bold Z_p\llg v\rrg$    
for arbitrary  $a\geqq 0$, we have 
  $$
  \Big(\jrac{d}{dv}\Big)^{pa}\varphi\in(pa)!\bold Z_p\llg v\rrg 
                              =p^aa!\bold Z_p\llg v\rrg
                             \subset p^a\bold Z_p\llg v\rrg. 
  $$
Because  $\alpha\leqq 1$, this is contained in  $p^{\lceil \alpha a\rceil}\bold Z_p\llg v\rrg$.  
Hence  $\varphi\in\bold Z_p\llg v\rrg^{(\alpha),\hskip 1pt 0}$.  
\fp
(2) By the Leibniz rule, we see
  $$
  \align
  \Big(&\jrac{d}{dv}\Big)^{pa}(\varphi\psi) 
  =\sum_{j=0}^{pa}
     \binom{pa}{j}
      \Big(\hskip -1pt\Big(\jrac{d}{dv}\Big)^j     \varphi\Big)
      \Big(\hskip -1pt\Big(\jrac{d}{dv}\Big)^{pa-j}\psi\Big)    \\
  &=\sum_{j=0}^a
     \binom{pa}{pj}
      \Big(\hskip -1pt\Big(\jrac{d}{dx}\Big)^{pj}   \varphi\Big)
      \Big(\hskip -1pt\Big(\jrac{d}{dv}\Big)^{pa-pj}\psi\Big)    
   +\hskip -6pt\sum_{0\leqq j \leqq pa, p\hskip -2pt\not\hskip 2pt|j}
    \hskip -4pt
     \binom{pa}{j}
      \Big(\Big(\jrac{d}{dv}\Big)^j     \varphi\Big)
      \Big(\Big(\jrac{d}{dv}\Big)^{pa-j}\psi\Big).
  \endalign
  $$
If  $p\hskip -2pt\not|\hskip 2ptj$  then  $p\big|\binom{pa}{j}$.  
The definition  (3.2.2) and 3.2.4 (7) imply
  $$
  \align
  \Big(\jrac{d}{dv}\Big)^{pa}(\varphi\psi)
 &\in\sum_{j=0}^a p^{\lceil \alpha(j+\beta_1)\rceil}
                  p^{\lceil \alpha(a-j+\beta_2)\rceil} 
                  \bold Z_p\llg v\rrg \\
 & \hskip -20pt
   + \sum_{0\leqq j\leqq pa,\  p\hskip -2pt\not\hskip 2pt| j} 
   \hskip -5pt
                 pp^{\lceil \alpha(\lfloor \frac{j}{p}\rfloor+\beta_1)\rceil}
                  p^{\lceil \alpha(\lfloor \frac{pa-j}{p}\rfloor+\beta_2)\rceil}
                  \bold Z_p\llg v\rrg \\
 & \subset p^{\lceil \alpha(a+\beta_1+\beta_2)\rceil}\bold Z_p\llg v\rrg 
       + pp^{\lceil \alpha(a-1+\beta_1+\beta_2)\rceil}\bold Z_p\llg v\rrg. 
  \endalign
  $$
By  $\alpha\leqq 1$, this is contained in   $p^{\lceil \alpha(a+\beta_1+\beta_2)\rceil}\bold Z_p\llg v\rrg$, 
and we have shown that 
  $$
  \varphi\psi\in\bold Z_p\llg v\rrg^{(\alpha),\hskip 1pt \beta_1+\beta_2}
  $$
as desired. 
\qed

\newpage

\fp
\lbox{53}{360}{{\bf Lemma 3.2.7.} \ \it
We have
  $$
  \exp\Big(v+\jrac{v^p}{p}\Big)\in \bold Z_p\llg v\rrg^{(1-\frac{1}{p}),\hskip 1pt 0}.
  $$
}
\vskip 3pt
\fp
{\it Proof.} \ 
Thanks to 3.2.6 (2) and (3), it suffice to show 
\fp\hskip 10pt
(1)  $\exp\big(\sum_{n=0}^{\infty}\frac{v^{p^n}}{p^n}\big)\in \bold Z_p[[v]]$G
\fp\hskip 10pt
(2)  $\exp\big(\hskip -3pt -\frac{v^{p^n}}{p^n}\big)\in \bold Z_p\llg v\rrg^{(1-\frac{1}{p}), \hskip 1pt 0}$  
for all  $n\geqq 2$.  \fp
Since the property (1) is well-known and proved in \cite{{\bf Ho}, p.238, 5.4} 
and \cite{{\bf R}, p.388, Theorem}, 
we show only (2).  
To do so, by 3.2.4 (3), it suffice to show
\fp\hskip 10pt
$(2)'$ $\frac{1}{m!}\big(\frac{v^{p^n}}{p^n}\big)^m\in\bold Z_p\llg v\rrg^{(1-\frac{1}{p}),\hskip 1pt 0}$  
for all  $n\geqq 2$  and  $m\geqq 0$.  \fp
This  $i2j'$  is equivalent to 
\fp\hskip 10pt
$(2)''$ $\ord_p\big(\frac{(p^nm)!}{m!p^{mn}}\big)\geqq (1-\tfrac{1}{p})p^{n-1}m$ 
for all  $n\geqq 2$  and  $m\geqq 0$, \fp
because of 3.2.4 (9).  
As 
  $$
  \ord_p\Big(\frac{(p^nm)!}{m!}\Big)
  =\sum_{j=0}^{n-1}p^jm
  \geqq p^{n-1}m+m,
  $$
to show $(2)''$ it suffice to check 
  $$
  p^{n-1}m+m-mn\geqq \big(1-\tfrac{1}{p}\big)p^{n-1}m, 
  $$
namely  $m(p^{n-2}-n+1)\geqq 0$.  
This is easily checked.  
\qed
\vskip 15pt
\fp
{{\it Proof of} 3.2.1.} \ 
Let  $n\geqq 0$  and  $m=-\lceil \frac{n}{p} \rceil$.  
By 3.2.4 (9), we see
  $$
  \jrac{v^n}{n!}\in\bold Z_p\llg v\rrg^{(1-\frac{1}{p}),\hskip 1pt m}\cap\bold Z_p\llg v\rrg.
  $$
By 3.2.7 and 3.2.6 (2) we have  $\frac{v^n}{n!}\exp\big(v+\frac{v^p}{p}\big)
\in\bold Z_p\llg v\rrg^{(1-\frac{1}{p}),\hskip 1pt m}\cap\bold Z_p\llg v\rrg$; 
and by 3.2.5 we have 
  $$
  \Big(\Big(\jrac{d}{dv}-1\Big)^p+1\Big)^a\Big(\jrac{v^n}{n!}\exp\Big(v+\jrac{v^p}{p}\Big)\Big)
  \in\bold Z_p\llg v\rrg^{(1-\frac{1}{p}),\hskip 1pt m+a}\cap\bold Z_p\llg v\rrg
  $$
for all  $a\geqq 0$.  
This is contained in  $p^{\lceil(1-\frac{1}{p})(m+a)\rceil}\bold Z_p\llg v\rrg\cap\bold Z_p\llg v\rrg$  
by the property 3.2.4 (2).  
Since 
  $$
  \Big(\Big(\jrac{d}{dv}-1\Big)^p+1\Big)^a\Big(\jrac{v^n}{n!}\exp\Big(v+\jrac{v^p}{p}\Big)\Big)
  =\exp(v)\Big(\Big(\jrac{d}{dv}\Big)^p+1\Big)^a\Big(\jrac{v^n}{n!}\exp\Big(\jrac{v^p}{p}\Big)\Big),
  $$
and  $\exp(v)$  is invertible in  $\bold Z_p\llg v\rrg$, we see 
  $$
  \Big(\Big(\jrac{d}{dv}\Big)^p+1\Big)^a\Big(\jrac{v^n}{n!}\exp\Big(\jrac{v^p}{p}\Big)\Big)
  \in p^{\lceil(1-\frac{1}{p})(m+a)\rceil}\bold Z_p\llg v\rrg\cap\bold Z_p\llg v\rrg. 
  $$
Substituting  $m=-\lfloor \frac{n}{p}\rfloor$, we have finally
  $$
  \Big(\Big(\jrac{d}{dv}\Big)^p+1\Big)^a\Big(\jrac{v^n}{n!}\exp\Big(\jrac{v^p}{p}\Big)\Big)
  \in 
  \cases
    \bold Z_p\llg v\rrg & \text{( if  $pa<n$ )}, \\
    p^{a-\lfloor \frac{n}{p}\rfloor-\lfloor(a-\lfloor \frac{n}{p}\rfloor)/p\rfloor}
    \bold Z_p\llg v\rrg & \text{( if  $pa\geqq n$ )}.\\
  \endcases
  $$
Looking at the coefficients of  $v^{n+qp}/(n+qp)!$  or  $v^{n-r_0p}/(n-r_0p)!$   
of the above after multiplying it by  $n!$, 
we get the desired congruences.  
\qed

\newpage \fp
{\bf 3.3. The second tool for the proof.}\ \ \ 
\fp
We show the second Proposition used in the proof of 3.1.1.  
\vskip 3pt
\fp
\lbox{84}{360}
{{\bf Proposition 3.3.1.}\ \ \it 
Let  $p$  be an odd prime and  $U$  is a partition 
with  $U_{p-1}=0$, $d(U)\neq 0$.  
Then  $\tau_{\sssize U}$  defined by {\rm (2.2.4)} satisfies   
  $$
  \ord_p (\tau_{\sssize U})
  \geqq \bigg\lfloor \frac{w(U)+d(U)-2}{2p} \biggr\rfloor.  
  $$
}
\vskip 5pt
\fp
{\bf Remark 3.3.2.} \ \ 
As we mentioned in 3.1.2 (2), 
we have  $\ord_p(\tau_{\sssize U})=(p-5)/2$  
for any prime  $p\geqq 5$  and any partition  $U$  such 
that  $U_1=p$, $U_{2p-1}=(p-5)/2$, with the others  $U_j=0$.    
For this  $U$, since  $w(U)=p+\frac{(2p-1)(p-5)}{2}$, 
$d(U)=p+\frac{p-5}{2}$, 
and  $\lfloor (w(U)+d(U)-2)/(2p) \rfloor
  =\big\lfloor (p^2-3p-2)/(2p) \big\rfloor
  =\lfloor \frac{p-5}{2}+\frac{p-1}{p} \rfloor
  =(p-5)/2$,  
the estimate above is best possible.  
\vskip 5pt
\fp
{\it Proof.} \  
For simplicity, we denote  $w(U)=n$  and  $d(U)=d$.  
As  $d\neq 0$, $n+d-2>0$.  
We show the estimate under the condition  $U_{2p-1}\neq 0$  as follows: 
  $$
  \allowdisplaybreaks
  \align
 & \ord_p (\tau_{\sssize U})  
  =\ord_p ( (n+d-2)! )-\ord_p ( \gamma_{\sssize U} ) \\
 &=\ord_p \bigg( \Big(-2+\sum_{j\neq p-1}(j+1)U_j \Big)!\bigg)
    - \hskip -7pt
      \sum_{(\epsilon, k)\neq (1,1)}
      \hskip -10pt kU_{\epsilon p^k-1}
    - \sum_{j\neq p-1}\ord_p (U_j!) \\
 & \hskip 40pt 
     (\text{ In the bellow  $\epsilon$  runs through the positive integers coprime to  $p$. }) \\
 & \geqq   
   \ord_p 
     \bigg(\Big(
      - 2
      + \hskip -7pt
        \sum_{j\neq p-1,~2p-1}\hskip -7pt jU_j 
      + 2p U_{2p-1}
     \Big)!\bigg)
      - \hskip -7pt
        \sum_{(\epsilon, k)\neq (1,1)}
        \hskip -7pt kU_{\epsilon p^k-1} 
      - \ord_p(U_{2p-1}!) \\
 &= \ord_p 
     \bigg(\Big(
      - 2
      + \sum_{p\not\hskip 2pt | j+1} jU_j 
      + \hskip -15pt
        \sum_{(\epsilon, k)\neq (1,1),~(2,1)}
        \hskip -15pt (\epsilon p^k-1)U_{\epsilon p^k-1}
      + 2p U_{2p-1}
     \Big)!\bigg) \\
 &  \hskip 50pt
   -\hskip -8pt
    \sum_{(\epsilon, k)\neq (1,1)}
    \hskip -8pt kU_{\epsilon p^k-1}
   -\ord_p(U_{2p-1}!) \\
 &\geqq  
    \sum_{\nu=1}^{\infty}
     \bigg\lfloor
      \frac{1}{p^{\nu}}
      \bigg(
      - 2 
      + \sum_{p\not\hskip 2pt | j+1} jU_j
      + \hskip -7pt
        \sum_{(\epsilon, k)\neq (1,1),~(2,1)} 
        \hskip -7pt (\epsilon p^k-1)U_{\epsilon p^k-1}
      + 2pU_{2p-1} 
      \bigg)
     \bigg\rfloor \\
 & \hskip 20pt
   -\hskip -7pt\sum_{(\epsilon, k)\neq (1,1)}
    \hskip -7pt kU_{\epsilon p^k-1} 
   -\ord_p(U_{2p-1}!) 
    \hskip 20pt (\ \because 
    \ord_p(N!)
    =\sum_{\nu=1}^{\infty}
     \big\lfloor \tfrac{N}{p^{\nu}} \big\rfloor\ ) \\
 &=~
     \bigg\lfloor
      \frac{1}{p}
      \bigg(
      - 2 
      + \sum_{p\not\hskip 2pt | j+1} jU_j
      + \hskip -7pt
        \sum_{(\epsilon, k)\neq (1,1),~(2,1)} 
        \hskip -7pt (\epsilon p^k-1)U_{\epsilon p^k-1}
      + 2pU_{2p-1} 
      \bigg)
     \bigg\rfloor \\
 & \hskip 25pt
    +\sum_{\nu=2}^{\infty}
     \bigg\lfloor
      \frac{1}{p^{\nu}}
      \bigg(
      - 2 
      + \sum_{p\not\hskip 2pt | j+1} jU_j
      + \hskip -7pt
        \sum_{(\epsilon, k)\neq (1,1),~(2,1)} 
        \hskip -7pt (\epsilon p^k-1)U_{\epsilon p^k-1}
      + 2pU_{2p-1} 
      \bigg)
     \bigg\rfloor \\
 & \hskip 60pt 
    -\hskip -7pt\sum_{(\epsilon, k)\neq (1,1)}
     \hskip -7pt kU_{\epsilon p^k-1} 
    -\ord_p(U_{2p-1}!) \\
 &\geqq 
     \bigg\lfloor
      \frac{1}{p}
      \bigg(
      - 2 
      + \sum_{p\not\hskip 2pt | j+1} jU_j
      + \hskip -7pt
        \sum_{(\epsilon, k)\neq (1,1),~(2,1)} 
        \hskip -7pt (\epsilon p^k-1)U_{\epsilon p^k-1}
      + 2pU_{2p-1} 
      \bigg)
     \bigg\rfloor \\
 &\hskip 25pt
    +\sum_{\nu=2}^{\infty}
     \bigg\lfloor
      \frac{- 2 + 2pU_{2p-1}}{p^{\nu}}
     \bigg\rfloor
    -\hskip -7pt\sum_{(\epsilon, k)\neq (1,1)}
     \hskip -7pt kU_{\epsilon p^k-1} 
    -\ord_p(U_{2p-1}!) \\
 &= \bigg\lfloor
      \frac{1}{p}
      \bigg(
      - 2 
      + \sum_{p\not\hskip 2pt | j+1} jU_j
      + \hskip -7pt
        \sum_{(\epsilon, k)\neq (1,1),~(2,1)} 
        \hskip -7pt (\epsilon p^k-1)U_{\epsilon p^k-1}
      + 2pU_{2p-1} 
      \bigg)
     \bigg\rfloor \\
   &\hskip 25pt
    -\hskip -7pt\sum_{(\epsilon, k)\neq (1,1)} 
     \hskip -7pt kU_{\epsilon p^k-1}
    +\sum_{\nu=2}^{\infty}
     \bigg\lfloor
      \frac{- 2 + 2pU_{2p-1}}{p^{\nu}}  
     \bigg\rfloor 
    -\ord_p(U_{2p-1}!) \\
 &= \bigg\lfloor
      \frac{1}{p}
      \bigg(
      - 2 
      + \sum_{p\not\hskip 2pt | j+1} jU_j
      + \hskip -7pt
        \sum_{(\epsilon, k)\neq (1,1),~(2,1)} 
        \hskip -7pt (\epsilon p^k-kp-1)U_{\epsilon p^k-1}
      + 2pU_{2p-1} 
      \bigg)
     \bigg\rfloor  \\
 &  \hskip 25pt
    -U_{2p-1}
    +\sum_{\nu=2}^{\infty}
     \bigg\lfloor
      \frac{- 2 + 2pU_{2p-1}}{p^{\nu}}  
     \bigg\rfloor 
    -\ord_p(U_{2p-1}!) \\
 &= \bigg\lfloor
      \frac{1}{p}
      \bigg(
      - 2 
      + \sum_{p\not\hskip 2pt | j+1} jU_j
      + \hskip -7pt
        \sum_{(\epsilon, k)\neq (1,1),~(2,1)} 
        \hskip -7pt (\epsilon p^k-kp-1)U_{\epsilon p^k-1}
      + 2pU_{2p-1} 
      \bigg)
     \bigg\rfloor \\ 
 &  \hskip 25pt
    - U_{2p-1}
    + \ord_p((-2+2pU_{2p-1})!)
    -\bigg\lfloor
      \frac{- 2 + 2pU_{2p-1}}{p}  
     \bigg\rfloor 
    - \ord_p(U_{2p-1}!).  
  \endalign
  $$
Here we can replace the first term  $-\frac{2}{p}$  by  $-\frac{1}{p}$,   
because if in the first bracket  $\lfloor \ \ \rfloor$  
the sum of the other terms is an integer then both of  $-\frac{2}{p}$  
and  $-\frac{1}{p}$  contribute as  $-1$, 
and if such the sum is not an integer its fractional part is larger than 
or equals to  $\frac{1}{p}$.  
Moreover, since  $j\leqq (j + 1)/2$  about the term  $jU_j$  in the second term sum, 
and  $\epsilon p^k -kp-1>\epsilon p^k/2$   if  $(\epsilon, k)\neq(1,1)$, $(2,1)$  
for the third term sum, we see
  $$
  \allowdisplaybreaks
  \align
 &\geqq \bigg\lfloor
      \frac{1}{2p}
      \bigg(\hskip -3pt
      - 2 
      + \hskip -3pt\sum_{p\not\hskip 2pt | j+1} (j+1)U_j
      + \hskip -7pt
        \sum_{(\epsilon, k)\neq (1,1),~(2,1)} 
        \hskip -14pt \epsilon p^kU_{\epsilon p^k-1}
      + 2pU_{2p-1} 
      \bigg) 
      + U_{2p-1}
     \bigg\rfloor 
    -U_{2p-1} \\
 &  \hskip 25pt
    + \ord_p((-2+2pU_{2p-1})!)
    -\bigg\lfloor
      \frac{- 2 + 2pU_{2p-1}}{p}  
     \bigg\rfloor 
    - \ord_p(U_{2p-1}!) \\
 &= \bigg\lfloor
      \frac{-2+w(U)+d(U)}{2p}
     \bigg\rfloor \\
 &  \hskip 6pt
    +\ord_p(2U_{2p-1})!)+2U_{2p-1}-\ord_p(2pU_{2p-1})
    -\bigg\lfloor
      \frac{- 2 + 2pU_{2p-1}}{p}  
     \bigg\rfloor 
    - \ord_p(U_{2p-1}!) \\
&= \bigg\lfloor
      \frac{-2+w(U)+d(U)}{2p}
     \bigg\rfloor \\
 &  \hskip 6pt
    +\ord_p(2U_{2p-1})!)+2U_{2p-1}-\ord_p(2U_{2p-1})-1
    -(-1+2U_{2p-1}) 
    - \ord_p(U_{2p-1}!) \\
 &\geqq 
    \bigg\lfloor
      \frac{-2+w(U)+d(U)}{2p}
    \bigg\rfloor
   +\ord_p\bigg(\frac{(2U_{2p-1}-1)!}{U_{2p-1}!}\bigg).
  \endalign
  $$
Hence, our proof is completed for the case  $U_{2p-1}\neq 0$. 
If  $U_{2p-1}=0$, by substituting this at the beginning of this calculation, 
we can more easily prove the desired estimate.
\qed

\newpage \fp
{\bf 3.4. Proof of the Kummer type congruence relation.}\ \ \ 
Let us establish the proof of 3.1.1. 
By 2.2.5 the left hand side of the congruence in 3.1.1 is written as follows:
  $$
  \sum_{r=0}^a\binom{a}{r}(-1)^r {f_{p-1}}^{\hskip -8pt a-r}
           \frac{\hat{B}_{n+r(p-1)}}{n+r(p-1)} 
 =\sum_{r=0}^a\binom{a}{r}(-1)^r {f_{p-1}}^{\hskip -8pt a-r}
  \hskip -15pt\sum_{w(U)=n+r(p-1)}\hskip -10pt \tau_{\sssize U}f^U,    
  \tag{3.4.1}
  $$
where
  $$
  \tau_{\sssize U}=(-1)^{d(U)-1}\frac{(w(U)+d(U)-2)!}{\gamma_{\sssize U}}. 
  $$
By picking up  ${f_{p-1}}^r$  or  ${f_{p-1}}^{r-r_0}$  from  $f^U$, 
we see 
  $$
  \aligned
  =\sum_{r=0}^a\binom{a}{r}&(-1)^r{f_{p-1}}^{a-r} \\
  &\Bigg\{\hskip -2pt
    \sum_{w(U)=n}
         \hskip -7pt\tau_{\sssize U[r]}f^U{f_{p-1}}^r
   +\sum_{r_0=1}^r\hskip -15pt
         \hskip -10pt
           \sum_{\hskip 6pt\Sb\hskip 20pt w(U)=n+r_0(p-1) \\ U_{p-1}=0\endSb}
              \hskip -20pt
                 \tau_{\sssize U[r-r_0]}f^U{f_{p-1}}^{r-r_0}
  \Bigg\},
  \endaligned
  \tag{3.4.2}
  $$
where  $U[r]$  means to increase  $r$  the  $(p-1)$-st entry of  $U$.  
After exchanging the sums about  $r$  and about  $U$, 
by writing down   $\tau_{\sssize U[r]}$  and   $\tau_{\sssize U[r-r_0]}$, we have
  $$
  \aligned
 =&\sum_{w(U)=n}\hskip -5pt
         \frac{f^U{f_{p-1}}^{a}}{\gamma_{\sssize U|_{p-1}}} 
       \sum_{r=0}^a
        \binom{a}{r}(-1)^{d(U[r])+r-1} 
         \frac{\{w(U[r])+d(U[r])-2\}!}{p^{r+U_{p-1}}(r+U_{p-1})!} \\
   &\hskip 10pt
     + \ \ \sum_{r_0=1}^a\hskip -22pt
        \sum_{\Sb \hskip 24pt w(U)=n+r_0(p-1) \\ U_{p-1}=0 \endSb}
    \hskip -22pt
         \frac{f^U{f_{p-1}}^{a-r_0}}{\gamma_{\sssize U}} 
     \hskip 5pt\Bigg\{
         \sum_{r=r_0}^a\binom{a}{r}(-1)^{d(U[r-r_0])+r-1} \\
  &\hskip 130pt
        \cdot
        \frac{\{w(U[r-r_0])+d(U[r-r_0])-2\}!}{p^{r-r_0}(r-r_0)!} \Bigg\}.
  \endaligned
  \tag{3.4.3}
  $$
Here  $U|_{p-1}$  means the  $(p-1)$-st entry of  $U$  by  $0$, 
so that  $\gamma_{\sssize U|_{p-1}}$  in the former sum means to delete 
the  $(p-1)$-st part  $p^{U_{p-1}}U_{p-1}!$  of  $\gamma_{\sssize U}$.  
For each  $r$, note that
  $$
  \aligned
  \gamma_{\sssize U[r]|_{p-1}}
 &=\gamma_{\sssize U|_{p-1}},  \\
  \gamma_{\sssize U[r]|_{p-1}}p^{r+U_{p-1}}(r+U_{p-1})!
 &=\gamma_{\sssize U[r]}, 
  \endaligned
  \tag{3.4.4}
  $$
and that, in the later sum, 
  $$
  \gamma_{\sssize U}
  =\left(2^{U_1}\cdots(p-1)^{U_{p-2}}(p+1)^{U_p}\cdots\right)\cdot
               \left(U_1!\cdots U_{p-2}!U_p!\cdots\right)
  \tag{3.4.5}
  $$
does not contain  $(p-1)$-st part. 
We denote the first and second sum of (3.4.3) by  $\sum_1$  and  $\sum_2$, respectively.  
and denote
  $$
  =\ {\sum}_1+{\sum}_2
  =\ \sum_{w(U)=n}
   \hskip -7pt S_1(U)
  \ \ + \ \ \sum_{r_0=1}^a
   \hskip -20pt
   \sum_{\Sb \hskip 20pt w(U)=n+r_0(p-1) \\ U_{p-1}=0 \endSb}
   \hskip -20pt S_2(U).
  \tag{3.4.6}
  $$
Then about  $S_1(U)$  for such a  $U$  that  $n=w(U)$, 
we have 
  $$
  \aligned
  w(U[r])+d(U[r])-2
  &=n+(p-1)r+d(U)+r-2 \\
  &=n+pr+d-2 \\
  &=(r+U_{p-1})p+n-pU_{p-1}+d(U)-2.
  \endaligned
  \tag{3.4.7}
  $$
Here we note that 
  $$
  \aligned
  n-pU_{p-1}+d(U)-2
 &=\big(n-(p-1)U_{p-1}\big)+(d(U)-U_{p-1})-2 \\
 & =w(U|_{p-1})+d(U|_{p-1})-2.
  \endaligned
  \tag{3.4.8}
  $$
About  $S_2(U)$  for such a  $U$  that  $w(U)=n+r_0(p-1)$, we see 
  $$
  \aligned
  w(U[r-r_0])+&d(U[r-r_0])-2 \\
  &=n+r_0(p-1)+(p-1)(r-r_0)+d(U)+(r-r_0)-2 \\
  &=(r-r_0)p+n+r_0p+d(U)-r_0-2.
  \endaligned
  \tag{3.4.9}
  $$
We divide  $\sum_1$ (resp.  $\sum_2$) into two parts 
according to  $n-pU_{p-1}+d(U)-2$  in (3.4.8) 
(resp. $n+r_0p+d(U)-r_0-2$ in (3.4.9)) is  $\geqq ap$  or  $< ap$, 
and denote the sum as  $\sum'_1 + \sum''_1$  
(resp. $\sum'_2 + \sum''_2$).  
Here we pay attention in (3.4.8) that  $n-pU_{p-1}+d(U)-2>0$ 
because of  $n\not\equiv 0$ mod $(p-1)$.  
\fp
(a) For the sum  $\sum'_1$, since  $n-pU_{p-1}+d(U)-2\geqq ap$, we have 
  $$
  \aligned
  \ord_p(S_1(U))
  &\geqq -\ord_p(\gamma_{\sssize U|_{p-1}})+\ord_p((n-pU_{p-1}+d-2)!) \\
  &\geqq \bigg\lfloor \frac{n-pU_{p-1}+d(U)-2}{2p}\bigg\rfloor \\
  &\geqq \bigg\lfloor \frac{ap}{2p}\bigg\rfloor \\
  &=\Big\lfloor \frac{a}{2}\Big\rfloor
  \endaligned
  \tag{3.4.10}
  $$
by 3.2.1 (1) and 3.3.1. 
\fp
(b) For the sum  $\sum''_1$, by denoting  $N=n-pU_{p-1}+d(U)-2$  and  $N=pb+e$ ($0\leqq e<p$), 
we see  $b<a$  because  $0<n-pU_{p-1}+d(U)-2< ap$.  
Therefore, by 3.3.1, we have
  $$
  \aligned
   \ord_p&(N!)
  -\ord_p(\gamma_{\sssize U|_{p-1}}) \\
 &=\ord_p((w(U|_{p-1})+d(U|_{p-1})-2)!)
  -\ord_p(\gamma_{\sssize U|_{p-1}}) \\
 &=\ord_p(\tau_{\sssize U|_{p-1}}) \\
 &\geqq \lfloor N/(2p) \rfloor.
  \endaligned
  \tag{3.4.11}
  $$
This and 3.2.1 (1) give
  $$
  \aligned
  \ord_p(S_1(U))
  &\geqq 
     a-\Big\lfloor \frac{N}{p} \Big\rfloor
    +\ord_p(N!) 
    -\Big\lfloor\frac{a-\lfloor N/p\rfloor}{p}\Big\rfloor 
    -\ord_p(\gamma_{\sssize U|_{p-1}})\\
  &=\bigg\lfloor \frac{N}{2p} \bigg\rfloor 
  +a-\Big\lfloor \frac{N}{p} \Big\rfloor
    -\Big\lfloor\frac{a-\lfloor N/p\rfloor}{p}\Big\rfloor \\
  &\geqq\Big\lfloor \frac{b}{2} \Big\rfloor+a-b
   -\Big\lfloor \frac{a-b}{p}\Big\rfloor \\
  &> \Big(\frac{b}{2}-1\Big)-b+a-\frac{a-b}{p} \\
  &=-1+\frac{a}{2}-\frac{(a-b)(p-2)}{2p} \\
  &=-1+\frac{a}{2}.  
  \endaligned
  \tag{3.4.12}
  $$
As the initial side is an integer, it must be  $\geqq \lfloor a/2 \rfloor$.  
Hence  $\ord_p(\sum_1)\geqq \lfloor a/2 \rfloor$.  
We can prove  $\ord_p(\sum_2)\geqq \lfloor a/2 \rfloor$    
by using 3.2.1 (2) instead of 3.2.1 (1).  
About such the proof, we should keep in mind that, 
for the case  $n+r_0p+d(U)-r_0-2\leqq ap$  in (3.4.9), 
in order to use 3.2.1 (2), 
it must be  $n+r_0p+d(U)-r_0-2\geqq r_0p$.  
However, since  $d(U)\geqq 1$, this condition is satisfied 
if  $n+r_0p+1-r_0-2\geqq r_0p$, 
namely  $n-r_0-1\geqq 0$.  
The latest condition follows from our assumption  $n>a$  
for all  $r_0=1$, $\cdots$, $a$.  
\qed
\vskip 10pt
\fp
{\bf Remark 3.4.13.} \ \ 
If we replace the condition  $a<n$  in 3.1.1 by  $a=n$, 
the Kummer type congruence relation for the generalized Bernoulli-Hurwitz numbers 
described below stands up modulo  $p^{n-1}$  
and this modulus is best possible. 
This fact means that the condition  $a<n$  is essential, 
because the generalized Bernoulli-Hurwitz
numbers are specialization of the universal Bernoulli numbers.

\newpage \fp
{\bf 3.5. The Kummer-Adelberg type congruence relation.}\ \ 
We prove the congruence relation 3.1.3 of Adelberg by using 3.1.1 and 3.1.2(1).   
\vskip 3pt
\fp
{{\it Proof of} \hskip 2pt 3.1.3.} \ 
We note that if  $p = 3$, the condition  $n\neq 0$, $1$ mod  $(p-1)$  
is always false and the statement of 3.1.3 is vacuous.  
So that, we may suppose  $p\geqq 5$.  
We prove the congruence relation by induction on  $a$.  
As we already mentioned the case  $a=1$  in 3.1.2 (1), 
and assume  $a>1$.  
By taking  $p^{a-1}$  as  $a$ in 3.1.1, we see
  $$
  \sum_{r=0}^{p^{a-1}}
  (-1)^r\binom{p^{a-1}}{r\ \ \ }{f_{p-1}}^{p^{a-1}-r}
  \frac{\hat{B}_{n+r(p-1)}}{n+r(p-1)}
  \equiv 0 \mod{p^{\lfloor p^{a-1}/2 \rfloor}}.
  \tag{3.5.1}
  $$
Here, we have  $\lfloor p^{a-1}/2\rfloor > a$  since  $a\geqq 2$  and  $p\geqq 5$.    
If  $r\neq 0$, $p^{a-1}$, then
  $$
  \aligned
  \ord_p\binom{p^{a-1}}{r\ \ \ }
  &=\frac{S_p(p^{a-1}-r)+S_p(r)-S_p(p^{a-1})}{p-1} \\
  &=\frac{S_p(p^{a-1}-r)+S_p(r)-1}{p-1}
  \endaligned
  \tag{3.5.2}
  $$
by (1.1.2). 
We denote  $\nu=\ord_pr$. 
Let
  $$
  \align
  p^{a-1}-r&=d_{a-2}p^{a-2}+d_{a-3}p^{a-3}+\cdots+d_1p+d_0 \hskip 20pt
  (0\leqq d_j\leqq p-1), \\
          r&=h_{a-2}p^{a-2}+h_{a-3}p^{a-3}+\cdots+h_1p+h_0 \hskip 20pt
  (0\leqq h_j\leqq p-1)
  \endalign
  $$
be their  $p$  base digit expressions. 
Then, obviously,
  $$
  d_j+h_j=
  \cases
  p-1,& (a-2\geqq j \geqq \nu+1),\\
  p,  & (j=\nu), \\
  0,  & (\nu-1\geqq j \geqq 0).
  \endcases
  \tag{3.5.3}
  $$
Namely,  $S_p(p^{a-1}-r)+S_p(r)=(p-1)(a -2-\nu)+p$.  
Thus  $\ord_p\binom{p^{a-1}}{r\ \ \ }=a-1-\nu$.    
Keeping in mind that  $p$  is an odd number, and consider the sums
  $$
  \aligned
   (-1)^r     &   \binom{p^{a-1}}{ r\ \ \ }{f_{p-1}}^{p^{a-1}-r}
   \frac{\hat{B}_{n+r(p-1)}}{n+r(p-1)} \\
 &+(-1)^{p^{a-1}-r} \binom{p^{a-1}}{p^{a-1}-r}{f_{p-1}}^{r}
   \frac{\hat{B}_{n+(p^{a-1}-r)(p-1)}}{n+(p^{a-1}-r)(p-1)}
  \endaligned
  \tag{3.5.4}
  $$
for  $0\leqq r\leqq (p^{a-1}-1)/2$.  
Summing up the consideration above, 
and the assumption of induction applied to  $\nu=\ord_pr=\ord_p(p^{a-1} -r)$, 
we see that (3.5.4) is divisible by  $p^a$.  
Hence 3.1.3.
\qed

\newpage \fp
\inbox{\bf 4}\ \ {\bf Hyperelliptic functions}
\vskip 3pt
\fp
{\bf 4.1. Fundamentals.} \ \ \ 
We consider a hyperellptic curve  $C$ of genus  $g$.  
We assume  $C$  is defined by  $y^2=f(x)$, where
  $$
  f(x)=\lambda_0x^{2g+1}+\lambda_1x^{2g}+\cdots+\lambda_{2g+1}
  $$
is a polynomial of  $x$  over  $C$, and  $f(x)=0$  has no multiple roots.  
We set  $\lambda_0=1$  for our convenience.  
The curve  $C$  is regarded as a non-singular algebraic curve with  
unique point  $\infty$  at infinity.  
As is well-known, the set of 
     $$
     \dfrac{x^{j-1}dx}{2y} \quad ij=1, \ldots, gj 
     \tag{4.1.1}
     $$
forms a basis of the differential forms of the first kind on  $C$.  
Let  $\left[\omega' \ \omega''\right]$  be the period matrix 
defined by taking a suitable set of generators of the fundamental group of  $C$, 
and let
     $$
     \Lambda
     := \omega' \tp{\left[\matrix \bold Z & \bold Z & \cdots & \bold Z
                  \endmatrix \right]}
     +  \omega''\tp{\left[\matrix \bold Z & \bold Z & \cdots & \bold Z
                  \endmatrix \right]}\ (\subset \bold C^g)
     $$
be the lattice in  $\bold C^g$  with respect to the differentials (4.1.1).  
We denote by  $J$  the jacobian variety of  $C$, 
and by  $\text{Sym}^g(C)$   the symmetric product of  $g$  copies of  $C$.  
Then we have a birational map
  $$
  \align
  \text{Sym}^g(C) &\rightarrow \text{Pic}^{\circ}(C)=J \\
  (P_1, \ldots, P_g)&\mapsto \text{the class of}\ P_1+\cdots +P_g-g\cdot\infty.
  \endalign 
  $$ 
We identify  $J$  as analytic manifolds with  $\bold C^g$.  
We denote the natural map  $\bold C^g\mapsto \bold C^g/\Lambda = J$  by  $\kappa$.  
The map
  $$
  \iota : Q \mapsto Q-\infty
  $$
is an embedding of  $C$  into  $J$.  
The pull-back  $\kappa^{-1}\iota(C)$   of  $\iota$  with respect to  $\kappa$  
is a {\it universal Abelian covering} of the curve  $C$.  
In analytic level, the birational  map above sends 
the each element  $(P_1, \cdots, P_g)\in \text{Sym}^g(C)$  to 
the point  $\bk{u}\mod{\varLambda}\in\bold C^g/\Lambda$, 
where
  $$
  \bk{u}=(u_1, \ldots, u_g)=
  \bigg(\int_{\infty}^{P_1}+\cdots+\int_{\infty}^{P_g}\bigg)
  (\omega_1, \ \ldots,  \ \omega_g).
  \tag{4.1.2}
  $$
The following notational convention is important.  
\vskip 3pt 
\fp
\lbox{28}{360}{{\bf Convention} : 
In this paper, we denote  $u_g$  simply by  $u$.   
}

\newpage

\fp
{\bf 4.2. The hyperelliptic functions and their variable.} \ 
In this paper, for each point  $\bk{u}\in\kappa^{-1}\iota(C)$   we denote by 
  $$
  (x(\bk{u}), \hskip 4pt y(\bk{u}))
  $$
the coordinate  $(x, y)$  of  $C$  such that  
$\kappa(\bk{u})=\iota\left(x(\bk{u}), \ y(\bk{u})\right)$.  
We call a rational expression of  $x(\bk{u})$  and  $y(\bk{u})$  
hyperelliptic function regarding as a function on  $\kappa^{-1}\iota(C)$.  
We confirm the following fundamental fact. 
\vskip 5pt
\fp
\lbox{90}{360}{\kern 0pt
{\bf Lemma 4.2.1.} \ \ \it 
Low degree terms of the Laurent developments of  $x(\bk{u})$  and  $y(\bk{u})$  
at  $\bk{u}=(0, \cdots, 0)$  with respect to  $u=u_g$  is given by
  $$
  x(\bk{u})=\frac 1{u^2}+\Dg{u}{0}, \quad
  y(\bk{u})=-\frac 1{u^{2g+1}}+\Dg{u}{-2g+1}.  
  $$
}
\vskip 3pt
\fp{\it Proof} \  
We take  $t=\dfrac 1{\sqrt{x}}$  as a local parameter at  $\infty$. 
Here we consider the brunch such that  $t>0$  for  $x>0$.    
We suppose that  $u\in\kappa^{-1}\iota(C)$  is sufficiently 
near to  $(0, 0, \ldots, 0)$, 
and that the three kinds of coordinates  $t$, $\bk{u}=(u_1, \ldots,  u_g)$,  $(x, y)$  
correspond to the same point of  $C$. 
Then
  $$
  \aligned
  u=u_g
  =&\int_{\infty}^{(x,y)}\frac {x^{g-1}dx}{2y} 
  =\int_{\infty}^{(x,y)}
   \frac{   x^{-3/2}dx   }
        {   
           2\sqrt{  
                     1+\lambda_1\frac 1x+\cdots
                     +\lambda_{2g+1}\frac 1{x^{2g+1}}
                  }
        } \\
  =&\int_0^t \frac{t^3\cdot\left(-\frac{2}{t^3}\right)dt}{2+\dg1}
  =-t+\Dg{t}2.
  \endaligned
  $$
Thus  $x(\bk{u})=\frac {1}{u^2}+\Dg{u}{-1}$.  
The similar fact is shown for  $y(\bk{u})$. 
By the definition, 
we see  $x(-\bk{u})=x(\bk{u})$, $y(-\bk{u})=-y(\bk{u})$  
The proof have been completed.
\qed
\vskip 5pt

The following properties are shown by the similar argument and omit their proofs.  
\vskip 3pt
\fp
\lbox{120}{360}
{{\bf Lemma 4.2.2.} \ \ \it
Let  $\bk{u}=(u_1, u_2, \ldots, u_g)$  be a variable on  $\kappa^{-1}\iota(C)$.  
Then
    $$
    \aligned
    u_1&=\tfrac1{2g-1} {u_g}^{2g-1} + (d^{\circ}(u_g)\geqq 2g), \\ 
    u_2&=\tfrac1{2g-3} {u_g}^{2g-3} + (d^{\circ}(u_g)\geqq 2g-2),\\
    &\cdots\cdots\cdots \\
    u_{g-1}&=\tfrac13 {u_g}^3 + (d^{\circ}(u_g)\geqq 4). 
    \endaligned
    $$
}

\fp
These facts suggest that
\vskip 3pt
\fp
\lbox{42}{360}
{It is natural to take  $u=u_g$  as the variable for 
the hyperelliptic functions locally near the  $\bk{u}=(0, \cdots, 0)$. 
}
\vskip 3pt
\fp

\newpage \fp
\inbox{\bf 5}\ \ {\bf Differential equations}
\vskip 3pt
\fp{\bf 5.1. Generalities. } \ 
For the hyperelliptic curve  $y(\bk{u})^2=f\big(x(\bk{u})\big)$ 
($f$  is a separable polynomial of degree  $2g+1$),  
the definition of  $u=u_g$  gives
  $$
  \frac{du}{dx}(\bk{u})=\frac{x^{g-1}}{2y}(\bk{u}).
  \tag{5.1.1}
  $$   
After squaring this, substituting the defining equation above of  $C$  into it, 
we have  $\big(\frac{du}{dx}\big)^{\hskip -1pt 2}=\frac{x^{2g-2}}{4f}$.  
Namely, 
\vskip 3pt
\fp
\lbox{48}{360}{
  $$
  x(\bk{u})^{2g-2}x'(\bk{u})^2=4f\big(x(\bk{u})\big) 
   \ \ \ \ \text{( $'$\  means  $\frac{d}{du}$ )}.
  \tag{5.1.2}
  $$
}
\fp
This (5.1.2) is just the analogy of  $\wp'(u)^2=4\wp(u)^3-g_2\wp(u)-g_3$  
in our theory
\footnote{Carlitz studied in \cite{{\bf Ca4}} about the solution of 
the differential equation  
$\big(\frac{dx}{du}(u)\big)^{\sssize 2}=$
\lq\lq a degree six polynomial of $x(u)$"  ($u\in\bold C$)  
instead of  $\wp'(u)^2=4\wp(u)^3-g_2\wp(u)-g_3$. 
The paper looks like not succeeded.
}.
Now we define the numbers  $C_n$  by the Laurent development of  $x(\bk{u})$  
with respect to  $u$: 
  $$
  x(\bk{u})
  =\frac1{u^2}+\frac{c_{-1}}{u}
  +\sum_{n=2}^{\infty}\frac{C_n}{n}\frac{u^{n-2}}{(n-2)!}
  \tag{5.1.3}
  $$
These  $C_n$  are analogy of Bernoulli and Hurwitz numbers.  
Though the definition of Bernoulli and Hurwitz numbers 
are based on the determination of their $2$-parts, 
such the property is a problem in the future for our theory. 
Of course, the recursion relation for  $C_n$  is obtained from (5.1.2). 
For the function  $y(\bk{u})$,  
we define also the numbers  $D_n$  by its Laurent expansion at  $u=0$ 
of  $y(\bk{u})$  with respect to  $u$:
  $$
  y(\bk{u})
  =\frac{-1}{u^{2g+1}}
  +\frac{d_{-2g}}{u^{2g}}
  +\cdots
  +\frac{d_{-1}}{u}
  +\sum_{n=2g+1}^{\infty}\frac{D_n}{n}\frac{u^{n-2g-1}}{(n-2g-1)!}.
  \tag{5.1.4}
  $$
Since the differential equation for  $y(\bk{u})$  is 
given by  $du=x^{g-1}dx/2y$, 
we also have the recursion relation for  $D_n$.

\newpage \input 05-2.o
         \input 05-3.o
\newpage \fp
\inbox{\bf 6}\ \ {\bf The Clarke type theorem}
\vskip 5pt
\fp
{\bf 6.1. The case of  $\bk{y}(\bk{u})^2=\bk{x}(\bk{u})^{2g+1}-{\bold 1}$.} \ 
In this case the Clarke type theorems for  $C_{(4g+2)n}$  and  $D_{(4g+2)n}$  
are stated as follows:
\vskip 5pt
\fp
\lbox{170}{360}
{{\bf Theorem 6.1.1.}\ \ \it
For each numbers  $C_{(4g+2)n}$  and  $D_{(4g+2)n}$  defined 
in {\rm (5.1.3)} and {\rm (5.1.4)}, respectively, for the curve  $y^2=x^{2g+1}-1$,  
we have  
  $$
  \align
  \frac{C_{(4g+2)n}}{(4g+2)n} 
  &\equiv
  - \hskip -18pt
   \sum_{\hskip 15pt\Sb p\equiv 1\hskip -3pt\mod{(2g+1)} \\ (4g+2)n=a(p-1) \endSb}
   \hskip -18pt
   \dfrac{{a|_p}^{-1}\mod{p^{1+\ord_pa}}}{p^{1+\ord_pa}}
   {A_p}^a
   \ \ \mod{\bold Z}, \\
  \frac{D_{(4g+2)n}}{(4g+2)n} 
  &\equiv
  -\hskip -18pt
   \sum_{\hskip 15pt \Sb p\equiv 1\hskip -3pt\mod{(2g+1)} \\ (4g+2)n=a(p-1) \endSb}
   \hskip -18pt
   \dfrac{{((2g)!a)|_p}^{-1}\mod{p^{1+\ord_pa}}}{p^{1+\ord_pa}}
   {A_p}^a
   \ \ \mod{\bold Z}, 
  \endalign
  $$
where
  $A_p=(-1)^{(p-1)/(4g+2)}\cdot \dbinom{(p-1)/2}{(p-1)/(4g+2)}$. 
}
\vskip 2pt
\fp
Obviously, these results show the following von Staudt-Clause type theorem and
an extension of von Staudt second theorem for  $C_{(4g+2)n}$  and  $D_{(4g+2)n}$.  
\vskip 5pt
\fp
\lbox{187}{360}
{{\bf Corollary 6.1.2.}\ \ \it
{\rm (1)} For each numbers  $C_{(4g+2)n}$  and  $D_{(4g+2)n}$,  
there are integers  $G_{(4g+2)n}$  and  $H_{(4g+2)n}$  such that 
  $$
  \align
  C_{(4g+2)n}
  &=\hskip -5pt\sum_{\Sb p\equiv 1\hskip -3pt\mod{(2g+1)} \\ p-1|(4g+2)n \endSb}
   \hskip -10pt
   \dfrac{{\ \ A_p}^{(4g+2)n/(p-1)}}{p}
   \ \ +\ G_{(4g+2)n}, \\
  D_{(4g+2)n} 
  &=\hskip -5pt\sum_{\Sb p\equiv 1\hskip -3pt\mod{(2g+1)} \\ p-1|(4g+2)n \endSb}
   \hskip -10pt
   \dfrac{((2g)!^{-1}\ \text{\rm mod}\ p)\ {A_p}^{(4g+2)n/(p-1)}}{p}
   \ \ +\ H_{(4g+2)n},
  \endalign
  $$
where  $A_p$  is the same one defined in {\rm 6.1.1}.  \fp
{\rm (2)} If  $(p-1)\hskip -3pt\not|\hskip 2pt (4g+2)n$, 
then  $C_{(4g+2)n}/((4g+2)n)$  and  $D_{(4g+2)n}/((4g+2)n)$  belong to  $\bold Z_p$.  
}
\vskip 5pt
\fp
{\it Proof.}\ 
Let  $p\equiv 1$ mod $(2g+1)$, $(4g+2)m=a(p-1)$, and  $\ord_pa=\nu$.  
Then 
  $$
  \frac{(4g+2)m\cdot\big({a|_p}^{-1}\hskip -5pt\mod{p^{1+\nu}}\big)}{p^{1+\nu}}
  \equiv 1-\frac{1}{p} \mod{p^{\nu}\bold Z}
  $$
and 
  $$
  \frac{(4g+2)m\cdot\big({a|_p}^{-1}\hskip -5pt\mod{p^{1+\nu}}\big)}{p^{1+\nu}}
  \equiv -\frac{1}{p} \mod{\bold Z}.
  $$
This yields the results.
\qed

\newpage 
\vskip 10pt

\fp
{\bf 6.2. The case of  $\bk{y}(\bk{u})^2=\bk{x}(\bk{u})^{2g+1}-\bk{x}(\bk{u})$.} \ 
In this case, the Clarke type theorem for  $C_{4gn}$  and  $D_{4gn}$  are 
stated as follows:
\vskip 3pt
\fp
\lbox{165}{360}
{{\bf Theorem 6.2.1.}\ \ \it
For each numbers  $C_{4gn}$  and  $D_{4gn}$  defined in {\rm (5.1.3)} and  {\rm (5.1.4)}, 
respectively, for the curve  $y^2=x^{2g+1}-x$, we have 
  $$
  \align
  \frac{C_{4gn}}{4gn} 
  &\equiv
   -\hskip -18pt
   \sum_{\hskip 15pt\Sb p\equiv 1\hskip -3pt\mod{4g} \\ 4gn=a(p-1) \endSb}
    \hskip -10pt
   \dfrac{{a|_p}^{-1}\mod{p^{1+\ord_pa}}}{p^{1+\ord_pa}}
   {A_p}^a
   \ \ \mod{\bold Z}, \\
  \frac{D_{4gn}}{4gn} 
  &\equiv
   -\hskip -18pt
   \sum_{\hskip 15pt \Sb p\equiv 1\hskip -3pt\mod{4g} \\ 4gn=a(p-1) \endSb}
   \hskip -10pt
   \dfrac{{((2g)!a)|_p}^{-1}\mod{p^{1+\ord_pa}}}{p^{1+\ord_pa}}
   {A_p}^a
   \ \ \mod{\bold Z},
  \endalign
  $$
where   $A_p=(-1)^{(p-1)/(4g)}\cdot\dbinom{(p-1)/2}{(p-1)/(4g)}$. 
}
\fp
These results also shows the following von Staudt-Clausen type theorem and an
extension of von Staudt second theorem for  $C_{4gn}$  and  $D_{4gn}$.  
\vskip 3pt
\fp
\lbox{163}{360}{{\bf Corollary 6.2.2.}\ \ \it  
{\rm (1)} For each numbers  $C_{4gn}$  and  $D_{4gn}$, 
there are integers  $G_{4gn}$  and  $H_{4gn}$  such that
  $$
  \align
  C_{4gn}
  &=\hskip -5pt\sum_{\Sb p\equiv 1\hskip -3pt\mod{4g} \\ p-1|4gn \endSb}
   \hskip -5pt
   \dfrac{{\ \ A_p}^{4gn/(p-1)}}{p}
   \ \ +\ G_{4gn},\\
  D_{4gn} 
  &=\hskip -5pt\sum_{\Sb p\equiv 1\hskip -3pt\mod{4g} \\ p-1|4gn \endSb}
   \hskip -5pt
   \dfrac{((2g)!^{-1}\ \text{\rm mod}\ p)\ {A_p}^{4gn/(p-1)}}{p}
   \ \ +\ H_{4gn},
  \endalign
  $$
where  $A_p$  is the same one defined in {\rm 6.2.1}. \fp
{\rm (2)} If  $(p-1)\hskip -2pt\not|\hskip 2pt 4gn$, then  $C_{4gn}/(4gn)$  and  $D_{4gn}/(4gn)$  
belong to  $\bold Z_p$. 
}
\vskip 2pt
\fp
{\it Proof. } \ 
Similar to 6.1.2. 
\qed
\vskip 8pt
\fp
{\bf 6.3. A remark on  $\bk{A}_p$.} \ 
Let  $p$  be a prime number such that  $p\equiv 1$  mod  $(2g + 1)$
or  $p\equiv 1$  mod  $4g$  according to  $C$   is defined by  $y^2 = x^{2g+1}-1$  
or  $y^2=x^{2g+1}-x$.  

We take the set of (4.1.1), namely
  $$
  \Big(\frac{1}{2y},\ \frac{x}{2y}, \ \cdots, \frac{x^{g-1}}{2y}\Big)
  $$
as the basis of the differential forms of the first kind on  $C$  mod  $p$.  
Then the Hasse-Witt matrix ($g\times g$  matrix) with respect to this basis is 
a diagonal matrix (\cite{{\bf Yu}, p.381}), 
and its  $(g, g)$-entry is just the number  $A_p$.  
Katz pointed out in the case of Hurwitz numbers that  $A_p$  is equal to 
the Hasse invariant (the unique entry of the Hasse-Witt matrix!) \cite{{\bf Ka}, p.2}.  
Our results are quite natural extension of that fact.

\newpage \fp
\inbox{\bf 7}\ \ {\bf Kummer type congruence relation}
\vskip 5pt
\fp
{\bf 7.1. The case of  $\bk{y}(\bk{u})^2=\bk{x}(\bk{u})^{2g+1}-{\bold 1}$.} \ 
In this case, the Kummer's original type congruence relations for  $C_{(4g+2)n}$  
and  $D_{(4g+2)n}$  are stated as follows:
\vskip 5pt
\fp
\lbox{173}{360}
{{\bf Theorem 7.1.1.}\ \ \it 
Let  $C_{(4g+2)n}$  and  $D_{(4g+2)n}$  be the numbers defined 
in {\rm (5.1.3)} and {\rm (5.1.4)}, 
respectively, for the curve  $y^2=x^{2g+1}-1$.  
For a prime  $p\equiv 1\mod{(2g + 1)}$  and positive integers  $a$  and  $n$  
such that  $(4g + 2)n-2\geqq a$, 
if  $(p-1)\not\vert\ (4g+2)n$, then we have 
  $$
  \align
  &\sum_{r=0}^a\dbinom{a}{r}(-A_p)^{a-r}
  \cdot\frac{C_{(4g+2)n+r(p-1)}}{(4g+2)n+r(p-1)}\equiv 0 \mod {p^a\bold Z_{(p)}}, \\
  &\sum_{r=0}^a\dbinom{a}{r}(-A_p)^{a-r}
  \cdot\frac{D_{(4g+2)n+r(p-1)}}{(4g+2)n+r(p-1)}\equiv 0 \mod {p^a\bold Z_{(p)}},
  \endalign
  $$
where
  $$
  A_p=(-1)^{(p-1)/(4g+2)}
  \cdot\dbinom{(p-1)/2}{(p-1)/(4g+2)}.
  $$
}\fp
\vskip 5pt
\fp
These congruence relations are just the same form as in 
Kummer's original paper \cite{{\bf Ku}} and in the case of Hurwitz numbers \cite{{\bf L}, p.193, (26)}. 
\vskip 5pt
\fp
{\bf Remark 7.1.2}. \ 
Under the facts (10.1.1), (10.1.2), and (10.1.3) proved later, 
we understand that 7.1.1 is satisfied modulo  $p^{\lfloor a/2\rfloor}$.  
\vskip 10pt

         \input 07-2.o
\newpage \fp
\inbox{\bf 8}\ \ {\bf Hurwitz-integral series}
\vskip 5pt
\fp
{\bf 8.1. Definition and basic properties.}\ \ \ 
We describe the notion of Hurwitz integrality and its properties.  
In this subsection the letter  $R$  always means a subring of 
the field  $\bold Q_p$  of the  $p$-adic numbers for a fixed prime  $p$.  
In our practice, $R$  will be  $\bold Z[\frac{1}{q}]$  with non-zero integer  $q$, 
the localization  $\bold Z_{(p)}$  of the integer ring  $\bold Z$  at  $p$, 
or its $p$-adic completion  $\bold Z_p$.   
\vskip 3pt
\fp
\lbox{100}{360}
{{\bf Definition 8.1.1.}\ \ \ 
Let  $z$  be an indeterminate or a variable. Let
  $$
  h(z)=\sum_{n=0}^{\infty}h_n\frac{z^n}{n!}\ \ \ \ (\ h_n\in \bold Q_p\ )
  $$
be a power series with respect to  $z$.   
If all the coefficients  $h_n$  belong to  $R$, 
then we say  $h(z)$  is Hurwitz integral over  $R$.   
The ring consists of the Hurwitz-integral series over  $R$  
with respect to  $z$  is denoted by  $R\llg z\rrg$.  
}
\vskip 3pt
\fp
It is easily seen that  $R\llg z\rrg$  is an integral domain, 
and is closed under the operations  $d/dz$  and  $\int_0^z\cdot\hskip 2pt dz$.  
A series  $h(z)\in R\llg z\rrg$  is a unit in  $R\llg z\rrg$  
if and only if its constant term is a unit in  $R$.  
Moreover the following properties are easily shown (see \cite{{\bf Hu2}, Section 1}).  
\vskip 5pt
\fp
\lbox{260}{360}{\kern 0pt
{\bf Proposition 8.1.2.}\ \ \it 
Let 
   $$
   h(z)=\sum_{n=0}^{\infty}h_n\frac{z^n}{n!}\ \ \ \ (\ h_n\in \bold Q_p\ ).
   $$
\fp
{\rm (1)} If the first  $n$  coefficients  $h_0$, $\cdots$, $h_{n-1}$  belong to  $R$, 
and there is a polynomial  $F(T_0, T_1, \cdots, T_{n-1})$  
of  $n$  variables over  $R$  such that there exists a relation
   $$ 
   h\lr{n}(z)=F(h(z), h'(z), \cdots, h\lr{n-1}(z))
   $$
on derivatives of  $h(z)$,  then  $h(z)\in R \llg z\rrg$.  
\fp
{\rm (2)} If  $h(z)\in R\llg z\rrg$, $h_0=0$, and  $h_1=1$,  
then the formal inverse series
   $$
   z=h^{-1}(w)=w+\cdots
   $$
of  $w=h(z)$  also belongs to  $R\llg z\rrg$.   
\fp
{\rm (3)} If  $h(z)\in R\llg z\rrg$, $h_0=0$, and  $h_1=1$, 
then for any positive integer  $m$,
   $$
   \frac{h(z)^m}{m!}
   $$
also belongs to  $R\llg z\rrg$.  
}
\vskip 3pt
\fp
Although there are so many criteria for Hurwitz integrality, 
they will be complicated to state.  
See 8.3 for example.

\newpage \input 08-2.o
         \fp
{\bf 8.3. Hurwitz-integrality of  ${\bold 1}/\bk{y}^{1/5}(\bk{u})$.} \ \ \
We also show the following fact. 
\vskip 5pt
\fp
\lbox{41}{360}
{{\bf Proposition 8.3.1.}\ \ \it
On the curve  $y^2=x^{2g+1}-1$  or  $y^2=x^{2g+1}-x$, 
we have 
\fp\hskip 80pt  
$s:=-1/{y(\bk{u})}^{1/(2g+1)} (=u+\cdots)\ \in\  \bold Z\llg u\rrg$. 
}
\fp
{\it Proof.} \ 
We restrict again the argument within only the curve  $y^2=x^5-1$.  
By (5.1.1), we see that
  $$
  du
  =\frac{xdx}{2y}
  =\frac{x\tfrac{dx}{dy}dy}{2y}
  =\frac{xdy}{2y\tfrac{dy}{dx}}
  =\frac{xdy}{5x^4}
  =\frac{dy}{5x^3}
  =\frac{dy}{5(y^2+1)^{3/5}},
  \tag{8.3.2}
  $$
and
  $$
  \frac{dy}{du}=5(y^2+1)^{3/5}D
  \tag{8.3.3}
  $$
Writing simply  $s'=ds/du$, $s''=d^2s/du^2$, $\cdots$, 
we have  $dy/du=-5s^{-6}ds/du$.  
Hence
  $$
  s'=-(1+s^{10})^{3/5}.
  \tag{8.3.4}
  $$
This formula gives the following equation by induction:  
For each integer  $n\geqq 1$, the  $n$-th derivative of  $s$  is written 
as a finite sum of the form
  $$
  s\lr{n}=\sum_{j}(1+s^{10})^{L_{nj}/5}P_{nj}(s,s',s'',\cdots, s\lr{n-1}),
  \tag{8.3.5}
  $$
where each  $P_{nj}$  means a polynomial of  $n$  variables over  $\bold Z$, 
and  $L_{nj}$  is an integer.  
Summing up this and  $s(0)=0$, 
we see that  $s'(0)$, $s''(0)$, $s\lr{3}(0)$, $\cdots$  are all integers.
Therefore we conclude that
  $$
  1/{y(\bk{u})}^{1/5}=-u-48\cdot 9!\ \frac{\dsize u^{\sssize 11}}{11!\ \ }+\cdots\in \bold Z\llg u\rrg
  \tag{8.3.6}
  $$
for the curve  $y^2=x^5-1$.  
\qed
\vskip 3pt
\fp

\newpage \fp
\inbox{\bf 9}\ \ \ {\bf Outline of the proof} \ \ \ 
\vskip 5pt
\fp
We sketch the proofs of 6.1.1 and 7.1.1 (also of 6.2.1 and 7.2.1) 
by taking the case of the curve  $y^2=x^5-1$  as an example. 
Recall our convention that, for a given power series  $\varphi(z)$  
with respect to  $z$,  
each coefficient  $\big[\frac{z^n}{n!}\big]\varphi(z)$ ( $n\geqq 0$ ) 
is called a {\it Hurwitz coefficient} of it.  
In the following, any Hurwitz coefficients are those of developments
with respect to  $u$.
\vskip 5pt
\fp
{\bf 9.1. The Clarke type theorem. } \ 
The proof of Clarke type theorem 6.1.1 is
divided into the following three steps.  
\vskip 3pt
\fp
\lbox{131}{360}
{\underbarl{{\it Step} 1}. \ \ 
We prove the Clarke type theorem 
for the Hurwitz coefficients  $C_n/n$  of  $x(\bk{u})$  
and  $C\lr{4}_n/(n)_4$  of  $x^2(\bk{u})=(1/t)^4$  over  $\bold Z[\frac12]$.  
\fp
\vskip 5pt
\fp
\underbarl{{\it Step} 2}. \ \  
We prove the Clarke type theorem 
for the Hurwitz coefficients  $D_n/n$  of  $y(\bk{u})$  
over  $\bold Z[\frac15]$. \fp
\vskip 5pt
\fp
\underbarl{{\it Step} 3}. \ \  
We make a connection of the results above.  
By using simple formulae concerning  $D=d/du$  
among the functions  $x(\bk{u})$, $x^2(\bk{u})$, $y(\bk{u})$, 
we can connect the Hurwitz coefficients of these functions, 
and show both of the results in Steps 1 and 2 are valid 
over  $\bold Z=\bold Z[\frac12]\cap\bold Z[\frac15]$.  
}
\vskip 3pt
\fp
More detailed outline is as follows.
\fp
\underbarl{{\it Step} 1}:  
We prove the Clarke type theorem (10.1.7) 
for the Hurwitz coefficients  $C\lr{1}_{10m}/(10m)$  of  $x^{1/2}(\bk{u})$  
over  $\bold Z[\frac12]$.  
Then we show
\vskip 3pt
\fp 
\hskip 10pt
a relation (10.2.4) between the Hurwitz coefficients 
of  $x(\bk{u})$  and of  $x^{1/2}(\bk{u})$, 
\fp 
\hskip 10pt
a relation (10.3.4) between the Hurwitz coefficients 
of  $x^{3/2}(\bk{u})$  and of  $x(\bk{u})$, 
\fp
\hskip 10pt
a relation (10.4.4) between the Hurwitz coefficients 
of  $x^2(\bk{u})$  and of  $x^{3/2}(\bk{u})$. 
\vskip 3pt
\fp
These relations yield the Clarke type theorems 
for these Hurwitz coefficients over  $\bold Z[\frac12]$ 
by using the integrality property (10.1.3) of the Carlitz coefficients 
(see the Convention) of the formal inverse series  $u$  
with respect to  $t=x^{-1/2}(\bk{u})$.  

The Clarke type theorem above 
for the Hurwitz coefficients of  $x^2(\bk{u})=\big(x^{1/2}(\bk{u})\hskip -1pt\big)^4$  
(namely, (10.5.3)) is important in Step 3.
\vskip 3pt
\fp
\underbarl{{\it Step} 2}: 
In 11.1, we prove the Clarke type theorem (11.1.7) for the Hurwitz numbers
$D\lr{1}_{10m}/(10m)_4$  of  $y^{1/5}(\bk{u})$  
by using the integrality (11.1.3) of the Carlitz coefficients
of the formal inverse series  $u$  with respect to  $s=y^{-1/5}(\bk{u})$.  
On the other hand, we show that
\vskip 3pt
\fp
\hskip 10pt
a relation (11.2.4) between the Hurwitz coefficients of  $y^{2/5}(\bk{u})$  and of  $y^{1/5}(\bk{u})$, 
\fp
\hskip 10pt
a relation (11.3.4) between the Hurwitz coefficients of  $y^{3/5}(\bk{u})$  and of  $y^{2/5}(\bk{u})$, 
\fp
\hskip 10pt
a relation (11.4.4) between the Hurwitz coefficients of  $y^{4/5}(\bk{u})$  and of  $y^{3/5}(\bk{u})$,
\fp
\hskip 10pt
a relation (11.5.4) between the Hurwitz coefficients of  $y(\bk{u})$  and of  $y^{4/5}(\bk{u})$. 
\vskip 3pt
\fp
Summing up these relations, 
we can connect the Hurwitz coefficients  $D_{10m}/(10m)$  
of  $y(\bk{u})=\left(y^{1/5}(\bk{u})\right)^5$  
to  $D\lr{1}_{10m}/(10m)$ (see (11.6.1) and (11.6.2)).  
Therefore we have the Clarke type theorem (11.6.3) 
for  $D_{10m}/(10m)$  over  $\bold Z[\frac15]$,  
and finish the Step 2.
\vskip 3pt
\fp
\underbarl{{\it Step} 3}: Let  $D=d/du$.  
We can show by using
  $$
  (Dx^2)(\bk{u})=4y(\bk{u})
  $$
that  $D_{10m}$   and  $C\lr{4}_{10m}$  are essentially the same numbers.  
Since   $D_{10m}/(10m)$ $\in\bold Z_{(2)}$,
the denominator of  $C\lr{4}_{10m}/(10m)_4$  contains neither a power of  $2$  nor of  $5$.   
Hence, we have the Clarke type theorem (12.1.4) for  $D_{10m}/(10m)$  over  $\bold Z$,  
namely the second formula in 6.1.1.  
Finally, the formula  
  $$
  D^2x(\bk{u})=6x^2(\bk{u})+\frac{4}{x^3(\bk{u})}
  $$
with  $D=d/du=(2y/x)d/dx$, and the Hurwitz-integrality (8.2.4) of  $1/x^3(\bk{u})$ 
connect  $C\lr{4}_{10m}/(10m)_4$  with  $C_{10m}/(10m)$.   
So that we have the desired Clarke type theorem for  $C_{10m}/(10m)$, 
namely the first formula in 6.1.1.
\vskip 3pt

Incidentally, Prof. T.Oda pointed out the method used frequently 
in the proof above, 
that is a method to compare the Hurwitz coefficients of two different power
of a given power series with the Carlitz coefficients of 
the formal inverse series of it, 
is a variant of Lagrange inversion formula 1.2.1.  
The reader who is interested in such the method is referred to \cite{{\bf W}, pp.128-133} 
(Lagrange-B\"urmann theorem) or \cite{{\bf Co}, pp.148-153}.  

There is also a congruence in \cite{{\bf Ad1}, Theorem 3.4} 
linking through various order universal Bernoulli numbers, 
namely the universal case of Hurwitz coefficients of various powers of  $1/t(u)$  
with respect to  $u$.  
While the author believe it would closely relate to our method, 
he could not use it.  

\newpage

\fp
{\bf 9.2. The Kummer type congruence relation.} \  
Let  $p\equiv 1$ mod  $5$  be a fixed prime number.  
Although we describe only about  $C_{10m}$, 
the theorem for  $D_{10m}$  is proved similarly.  
The outline of the proof is as follows.
\vskip 3pt
\fp
\underbarl{{\it Step} 1}: \ \ 
We show that the formal inverse series  $u(t)\in \bold Z_p\llg t \rrg$  
of  $u\mapsto t=x(\bk{u})^{-1/2}$  satisfies the Honda's criterion (13.1.10) (see (13.1.10)).  
This fact implies the formal group law  $F$  whose formal logarithm is  $u(t)$, 
namely
  $$
  F(t_1, t_2):={u}^{-1}(u(t_1)+u(t_2))
  $$
is defined over  $\bold Z_p$ (see 13.1.1).  
\vskip 3pt
\fp
\underbarl{{\it Step} 2}: \ \  
By using Hochshilt's formula and Honda's property (13.1.10), 
we see that
  $$
  \Big(\big(\tfrac{d}{du}\big)^p-A_p\tfrac{d}{du}\Big)t(u)\in p\bold Z_p[[t(u)]] 
  \tag{9.2.1}
  $$
for  $t=t(u)=x(\bk{u})^{-1/2}\in\bold Z_p\llg u\rrg$, in 13.2.  
\vskip 3pt
\fp
\underbarl{{\it Step} 3}: \ \ 
Let  $\xi\in\bold Z_p$  be a primitive  $(p-1)$-st root of  $1$.  
The multiplication by  $\xi$  
  $$
  F_{\xi}(t):={u}^{-1}\big(\xi u(t)\big)=\xi+\cdots.
  $$
determined from  $F$  belongs 
to  $\xi t+t^2\bold Z_p[[t]]$  by the results in the Step 1.  
Now we use the notation  $x\langle u \rangle=x(\bk{u})\in\frac{1}{u^2}+\bold Q_p[[u]]$  
in order to avoid confusion.  
Thanks to the existence of such  $F_{\xi}(t)$, 
we can easily show that 
  $$
  x\langle u(t)\rangle-\xi^2x\langle \xi u(t)\rangle\in \bold Z_p[[t]].
  \tag{9.2.2}
  $$
Adding (9.2.1) and (9.2.2), we have, for any integer  $a>0$, that
  $$
  \Big(\big(\tfrac{d}{du}\big)^p-A_p\tfrac{d}{du}\Big)^a
  \big(x\langle u(t)\rangle-\xi^2x\langle\xi u(t)\rangle\big) 
  \in p^a\bold Z_p[[t]]
  \subset p^a\bold Z_p\llg u \rrg.
  \tag{9.2.3}
  $$
The coefficient of  $u^{10n-a-2}/(10n-a-2)!$  is 
just the left hand side of the first formula in 7.1.1, 
and the poof is completed.
\vskip 5pt

We mention here that the Kummer type congruence relations 
for  $C\lr{\nu}_{10m}$  and  $D\lr{\nu}_{10m}$  
without the division by  $10m$  are proved in 14.1.

\newpage \fp
\inbox{\bf 10}\ \ {\bf The Clarke type theorem on  $\bk{x}(\bk{u})$  over  $\bold Z[\frac12]$.} \ \ 
\vskip 5pt
\fp
{\bf 10.1. The Clarke type theorem on  $\bk{x}(\bk{u})^{1/2}$.} \ 
Let us consider the formal inverse series of  $u\mapsto \sqrt{x(\bk{u})}$,   
namely the power series  $u$  with respect to
  $$
  t=\frac{-1}{x^{1/2}(\bk{u})}.
  $$
That is 
  $$
  \align
  u
  &=\int_{\infty}^x \frac{xdx}{2y}  \\
  &=\int_0^{-t}
    \frac{
            \frac{1}{t^2}\cdot\left(-\frac{2dt}{t^3}\right)
         }
         {
            2\sqrt{  \frac{1}{t^{10}}-1  }
         }  \\
  &=-\int_0^{-t}
      \frac{1}{\sqrt{1-t^{10}}}dt \\
  &=-\int_0^{-t}
     \bigg(
      1+\sum_{m=1}^{\infty}
        (-1)^m\binom{-\frac{1}{2}}{m}t^{10m}
     \bigg)dt \\
  &=t+\sum_{m=1}^{\infty}
        (-1)^m\binom{-\frac{1}{2}}{m}\frac{t^{10m+1}}{10m+1}.
  \endalign
  $$
For convenience of quotation, we write again this;
  $$
  u=t+\sum_{m=1}^{\infty}
      (-1)^m\binom{-\frac{1}{2}}{m}\frac{t^{10m+1}}{10m+1}. 
  \tag{10.1.1}
  $$
We denote the Carlitz coefficients by
  $$
  f_{10m}=(-1)^m\binom{-\frac{1}{2}}{m}.
  \tag{10.1.2}
  $$
This notation corresponds to that in the definition of 
the universal Bernoulli numbers.  
By 1.3.1 and (1.3.2), we see that
  $$
  f_{10m}\in\bold Z[\tfrac12].
  \tag{10.1.3}
  $$
The other coefficients  $f_n$  are  $0$.  
Then the divided universal Bernoulli number  $\hat{B}_{10m}/(10m)$  
is specialized to  ${C\lr{1}_{10m}}/(10m)$,  
and the expression of 2.2.5 yields that
  $$
  \frac{C\lr{1}_{10m}}{10m}\ \in\ 3!\bold Z[\tfrac12]
  \tag{10.1.4}
  $$
as follows.  
To prove this, we apply 3.3.1 for  $p=2$  and  $p=3$.  
Since  $f_{2-1}=0$  and  $f_{3-1}=0$  in this case, 
we may consider only the partitions  $U$  such that  $U_{2-1}=U_{3-1}=0$.  
Here, as
  $$
  w(U)=10m\geqq 10, \ \ 
  d(U)\geqq 1
  $$
we see
  $$
  \ord_2(\tau_{\sssize U})
  \geqq\lfloor \tfrac{10+1-2}{4}\rfloor
  =2, \ \ \ 
  \ord_3(\tau_{\sssize U})
  \geqq\lfloor \tfrac{10+1-2}{6}\rfloor
  =1.
  $$
Hence, in the expression of 2.2.5 for  $C\lr{1}_{10m}/(10m)$, 
all its coefficients are  divisible by 3!.   
Adding to (10.1.3), we conclude (10.1.4).

If  $p=10m+1$  is a prime number, 
then the coefficient  $f_{p-1}$  of  $t^{10m+1}/(10m + 1)$  modulo  $p$  
coincide with  $(2,2)$-entry of the Hasse-Witt matrix with respect to 
our canonical base of the differential of the first kind on the curve  $C$  
of reduction modulo  $p$, namely
  $$
  f_{p-1}=(-1)^{(p-1)/10}\binom{-\frac{1}{2}}{\frac{p-1}{10}}
  \equiv A_p \mod{p}.
  \tag{10.1.5}
  $$
For this, see 18.1.  

Under the consideration above, we apply the Clarke's theorem 2.3.1 for the
function   $u\mapsto x(\bk{u})^{1/2}$.  
Let  
  $$
  \frac1t=\sum_{m=0}^{\infty}\frac{C\lr{1}_{10m}}{10m}
                             \frac{{u_2}^{10m-1}}{(10m-1)!}
  \tag{10.1.6}
  $$
Proposition 2.3.1 and (10.1.4) show that
\vskip 5pt
\fp
\lbox{65}{360}{
  $$
  \frac{1}{3!}\frac{C\lr{1}_{10m}}{10m}
    \in 
    \hskip -10pt
     \sum_{\Sb p\equiv 1\hskip -3pt\mod{5} \\ 10m=a(p-1) \endSb}
   \hskip -10pt\frac{{(3!a)|_p}^{-1}\hskip -5pt\mod p^{1+\ord_pa}}{p^{1+\ord_pa}}{A_p}^{a}
  +\bold Z[\tfrac12].
  \tag{10.1.7}
  $$
}
\vskip 1pt
\fp

\newpage \fp
{\bf 10.2.  Relation between coefficients of  $\bk{x}(\bk{u})^{1/2}$  and of  $\bk{x}(\bk{u})$.} \ \ 
For the curve  $y^2=x^5-1$  and 
  $$ 
  t=\frac{-1}{x(\bk{u})^{1/2}}\ \ (=u+\cdots),
  $$
we let 
  $$
  \frac{1}{t^2}
  =x(\bk{u})
  =\frac{1}{u^2}\sum_{m=1}^{\infty}\frac{C\lr{2}_{10m}}{(10m)_2}\frac{{u}^{10m-2}}{(10m-2)!}.
  \tag{10.2.1}
  $$
Then
  $$
  \aligned
  \int_0^u
  \left(\frac{1}{t^2}-\frac{1}{u^2}\right)du 
  &=\int_0^{u}
    \left(
    \sum_{m=1}^{\infty}\frac{C\lr{2}_{10m}}{(10m)_2}\frac{u^{10m-2}}{(10m-2)!}
    \right)du \\
  &=\sum_{m=1}^{\infty}\frac{C\lr{2}_{10m}}{(10m)_2}\frac{u^{10m-1}}{(10m-1)!}.
  \endaligned
  \tag{10.2.2}
  $$
On the other hand, after differentiating (10.1.1) with respect to  $u$, 
dividing by  $t^2$, we have
  $$
  \frac{1}{t^2}
  =\frac{1}{t^2}\frac{dt}{du}
    +\sum_{m=1}^{\infty}(-1)^m\binom{-\frac{1}{2}}{m}
                           t^{10m-2}\frac{dt}{du}.
  \tag{10.2.3}
  $$
Namely, 
  $$
  \int_0^u
   \left(
    \frac{1}{t^2}-\frac{1}{u^2}
   \right)du_2
  =\left(-\frac{1}{t}
   +  \sum_{m=1}^{\infty}(-1)^m\binom{-\frac{1}{2}}{m}
                           \frac{t^{10m-1}}{10m-1}\right)
   + \frac{1}{u}.
  $$
By equating this with (10.2.2) and by using (10.1.6), we see that
  $$
  \align
  &\sum_{m=1}^{\infty}\frac{C\lr{2}_{10m}}{(10m)_2}\frac{u^{10m-1}}{(10m-1)!} \\
  &\hskip 30pt =\left(
  -\sum_{m=0}^{\infty}\frac{C\lr{1}_{10m}}{10m}\frac{u^{10m-1}}{(10m-1)!}
  +\sum_{m=1}^{\infty}(-1)^m\binom{-\frac{1}{2}}{m}
                           \frac{t^{10m-1}}{10m-1}
  \right)
  +\frac{1}{u}.
  \endalign
  $$
Since  $C\lr{1}_0=1$, we can remove the terms of negative degree, and we have
  $$
   \sum_{m=1}^{\infty}\frac{C\lr{1}_{10m}}{10m}\frac{u^{10m-1}}{(10m-1)!}
  +\sum_{m=1}^{\infty}\frac{C\lr{2}_{10m}}{(10m)_2}\frac{u^{10m-1}}{(10m-1)!}
  =\sum_{m=1}^{\infty}(-1)^m\binom{-\frac{1}{2}}{m}\frac{t^{10m-1}}{10m-1}.
  $$
The right hand side of this belongs to  $\bold Z[\frac{1}{2}]\llg u\rrg$  
because of 1.3.1, 8.1.2 (3), and (8.2.4). 
Adding this with (10.1.3), we see that
  $$
  \frac{C\lr{1}_{10m}}{10m} 
  +\frac{C\lr{2}_{10m}}{(10m)_2} 
      \in (10\hskip -1pt\cdot\hskip -3pt 1\hskip -2pt -\hskip -1pt 2)!\bold Z[\tfrac12]
      \subset 3!\bold Z[\tfrac12].
  \tag{10.2.4}
  $$
As  $C_{10m}/(10m)=C\lr{2}_{10m}/(10m)_2$  by the definition, we have
  $$
  \frac{C_{10m}}{10m}\in -\frac{C\lr{1}_{10m}}{10m}+3!\bold Z[\tfrac12]
  \tag{10.2.5}
  $$
by (10.2.4).

We remark that we can conclude the Clarke type theorem for  $C_{10m}$  
over  $\bold Z[\frac{1}{2}]$  by (10.2.5) above and (10.1.7).  
We need, however, in order to prove the Clarke type theorem over  $\bold Z$  
to take a rather long way mentioned in 9.1.

\newpage \fp
{\bf 10.3. Relation between coefficients of  $\bk{x}(\bk{u})$  and of  $\bk{x}(\bk{u})^{3/2}$.} \ \ 
We prove here a relation of Hurwitz coefficients of  $x(\bk{u})$  and of  $x^{3/2}(\bk{u})$  
for the curve  $y^2=x^5-1$.  
For
  $$ 
  t=\frac{-1}{x(\bk{u})^{1/2}}\ \ (=u+\cdots),
  $$
we consider 
  $$
  \aligned
  \frac{1}{t^3}
  &=-x^{3/2}(\bk{u})  \\
  &=\sum_{m=0}^{\infty}\frac{C\lr{3}_{10m}}{(10m)_3}
                       \frac{u^{10m-3}}{(10m-3)!}
  \ \ \ \ (C\lr{3}_0=1).
  \endaligned
  \tag{10.3.1}
  $$
We have 
  $$
  \aligned
  \int_0^u
  \left(\frac{1}{t^3}-\frac{1}{u^3}\right)du 
  &=\int_0^u
    \left(
     \sum_{m=1}^{\infty}\frac{C\lr{3}_{10m}}{(10m)_3}
                        \frac{u^{10m-3}}{(10m-3)!}
    \right)du \\
  &=\sum_{m=1}^{\infty}\frac{C\lr{3}_{10m}}{(10m)_3}
                                    \frac{u^{10m-2}}{(10m-2)!}.
  \endaligned
  \tag{10.3.2}
  $$
On the other hand, after differentiating (10.1.1) with respect to  $u$, 
dividing by  $t^3$, we have
  $$
  \frac{1}{t^3}
  =\frac{1}{t^3}\frac{dt}{du}
    +\sum_{m=1}^{\infty}(-1)^m\binom{-\frac{1}{2}}{m}
                           t^{10m-3}\frac{dt}{du}. 
  \tag{10.3.3}
  $$
Namely, 
  $$
  \int_0^u
   \left(
    \frac{1}{t^3}-\frac{1}{u^3}
   \right)du
  =\left(-\frac12\frac{1}{t^2}
   +  \sum_{m=1}^{\infty}(-1)^m\binom{-\frac{1}{2}}{m}
                           \frac{t^{10m-3}}{10m-3}\right)
   + \frac12\frac{1}{{u_2}^2}. 
  $$
By equating this with (10.3.2) and by using (10.2.1), we see that
  $$
  \align
  &\sum_{m=1}^{\infty}\frac{C\lr{3}_{10m}}{(10m)_3}\frac{u^{10m-2}}{(10m-2)!} \\
  &\hskip 30pt =\left(
   -\frac12\sum_{m=0}^{\infty}
                      \frac{C\lr{2}_{10m}}{(10m)_2}\frac{u^{10m-2}}{(10m-2)!}
  +\sum_{m=1}^{\infty}(-1)^m\binom{-\frac{1}{2}}{m}
                           \frac{t^{10m-2}}{10m-2}
  \right)
  +\frac12\frac{1}{u^2}.
  \endalign
  $$
Since  $C\lr{2}_0=1$, we can remove the terms of negative degree, and we have
  $$
  \sum_{m=1}^{\infty}\frac{ C\lr{2}_{10m}}{(10m)_2}\frac{u^{10m-2}}{(10m-2)!}
 +\sum_{m=1}^{\infty}\frac{2C\lr{3}_{10m}}{(10m)_3}\frac{u^{10m-2}}{(10m-2)!}
  =2\sum_{m=1}^{\infty}(-1)^m\binom{-\frac{1}{2}}{m}
                           \frac{t^{10m-2}}{10m-2}.
  $$
The right hand side belongs to  $\bold Z[\frac{1}{2}]\llg u\rrg$  
because of 1.3.1, 8.1.2 (3), and (8.2.4).
Adding this with (10.1.3), we see that
  $$
  \frac{C\lr{2}_{10m}}{(10m)_2}+\frac{2C\lr{3}_{10m}}{(10m)_3} 
      \in (10\hskip -1pt\cdot\hskip -3pt 1\hskip -2pt -\hskip -1pt 3)!\bold Z[\tfrac12]
      \subset 3!\bold Z[\tfrac12].
  \tag{10.3.4}
  $$

\newpage \fp
{\bf 10.4. Relation between coefficients of  $\bk{x}(\bk{u})^{3/2}$  and of  $\bk{x}(\bk{u})^2$.} \ \ \ 
We prove here a relation of Hurwitz coefficients of  $x(\bk{u})^{3/2}$  
and of  $x^2(\bk{u})$  for  $y^2=x^5-1$.  
For
  $$ 
  t=\frac{-1}{x(\bk{u})^{1/2}}\ \ (=u+\cdots),
  $$
we consider
  $$
  \aligned
  \frac{1}{t^4}
  &=x^2(\bk{u})  \\
  &=\sum_{m=0}^{\infty}\frac{C\lr{4}_{10m}}{(10m)_4}\frac{{u}^{10m-4}}{(10m-4)!}
  \ \ \ \ (C\lr{4}_0=1).
  \endaligned
  \tag{10.4.1}
  $$
Then
  $$
  \aligned
  \int_0^{u}
  \left(\frac{1}{t^4}-\frac{1}{{u}^4}\right)du 
  &=\int_0^{u}
    \left(
    \sum_{m=1}^{\infty}\frac{C\lr{4}_{10m}}{(10m)_4}\frac{{u}^{10m-4}}{(10m-4)!}
    \right)du \\
  &=\sum_{m=1}^{\infty}\frac{C\lr{4}_{10m}}{(10m)_4}\frac{{u}^{10m-3}}{(10m-3)!}.
  \endaligned
  \tag{10.4.2}
  $$
On the other hand, after differentiating (10.1.1) with respect to  $u$, 
dividing by  $t^4$, we have
  $$
  \frac{1}{t^4}
  =\frac{1}{t^4}\frac{dt}{du}
    +\sum_{m=1}^{\infty}(-1)^m\binom{-\frac{1}{2}}{m}
                           t^{10m-4}\frac{dt}{du}.
  \tag{10.4.3}
  $$
Namely,
  $$
  \int_0^{u}
   \left(
    \frac{1}{t^4}-\frac{1}{{u}^4}
   \right)du
  =\left(-\frac13\frac{1}{t^3}
   +  \sum_{m=1}^{\infty}(-1)^m\binom{-\frac{1}{2}}{m}
                           \frac{t^{10m-3}}{10m-3}\right)
   + \frac13\frac{1}{{u}^3}.
  $$
By equating this with (10.4.2) and by using (10.3.1), we see that
  $$
  \align
  &\sum_{m=1}^{\infty}\frac{C\lr{4}_{10m}}{(10m)_4}\frac{{u}^{10m-3}}{(10m-3)!} \\
  &\hskip 30pt =\left(
   -\frac13\sum_{m=0}^{\infty}
                      \frac{C\lr{3}_{10m}}{(10m)_3}\frac{{u}^{10m-3}}{(10m-3)!}
  +\sum_{m=1}^{\infty}(-1)^m\binom{-\frac{1}{2}}{m}
                           \frac{t^{10m-3}}{10m-3}
  \right)
  +\frac13\frac{1}{{u}^3}.
  \endalign
  $$
Since  $C\lr{3}_0=1$, we can remove the terms of negative degree, and we have
  $$
   \sum_{m=1}^{\infty}\frac{ C\lr{3}_{10m}}{(10m)_3}\frac{{u}^{10m-3}}{(10m-3)!}
  +\sum_{m=1}^{\infty}\frac{3C\lr{4}_{10m}}{(10m)_4}\frac{{u}^{10m-3}}{(10m-3)!}
  =3\sum_{m=1}^{\infty}(-1)^m\binom{-\frac{1}{2}}{m}
                           \frac{t^{10m-3}}{10m-3}.
  $$
The right hand side belongs to  $\bold Z[\frac 12]\llg u\rrg$  
because of 1.3.1, 8.1.2 (3), and (8.2.4).  
Adding this with (10.1.3), we see that
  $$
  \frac{C\lr{3}_{10m}}{(10m)_3}+\frac{3C\lr{4}_{10m}}{(10m)_4} 
      \in (10\hskip -1pt\cdot\hskip -3pt 1\hskip -2pt -\hskip -1pt 4)!\bold Z[\tfrac12]
      \subset 3!\bold Z[\tfrac12].
  \tag{10.4.4}
  $$

\newpage \input 10-5.o
\newpage \fp
\inbox{\bf 11}\ \ 
{\bf The Clarke type theorem on  $\bk{y}(\bk{u})$  over  $\bold Z[\frac15]$} \ \ 
\vskip 5pt
\fp
{\bf 11.1. The Clarke type theorem on  $\bk{y}^{1/5}(\bk{u})$.} \ \ \
Being parallel to the Section 10, we consider 
the formal inverse series of  $u\mapsto {y(\bk{u})}^{1/5}$.  
Let
  $$
  y(\bk{u})=-1/{s^5}.
  $$
Then 
  $$
  \align
  u
  &=\int_{\infty}^x \frac{xdx}{2y} 
   =\int_{\infty}^y \frac{x}{2y}\frac{dx}{dy}dy 
   =\int_{\infty}^y \frac{x}{5x^4}dy 
   =\int_{\infty}^y \frac{1}{5x^3}dy \\
  &=\int_{\infty}^y \frac{1}{5(y^2+1)^{3/5}}dy 
   =\int_0^s
    \frac{\frac{5ds}{s^6}}
         {5(\frac1{s^{10}}+1)^{3/5}}  \\
  &=\int_0^s
      \frac{1}{(1+s^{10})^{3/5}}ds \\
  &=\int_0^s
    \bigg(1-\frac{1}{1!}\frac35                  s^{10}
            +\frac{1}{2!}\frac35\frac85           s^{20}
            -\frac{1}{3!}\frac35\frac85\frac{13}5 s^{30}
            +\cdots\bigg)ds \\
  &=\int_0^s
     \bigg(
      1+\sum_{m=1}^{\infty}
        \binom{-\frac{3}{5}}{m} s^{10m}
     \bigg)ds \\
  &=s+\sum_{m=1}^{\infty}
        \binom{-\frac{3}{5}}{m} \frac{s^{10m+1}}{10m+1}.
  \endalign
  $$
For convenience of quotation, we rewrite this:
  $$
  u
  =s+\sum_{m=1}^{\infty}
        \binom{-\frac{3}{5}}{m} \frac{s^{10m+1}}{10m+1}.
  \tag{11.1.1}
  $$
We denote each the Carlitz coefficient of this by
  $$
  f_{10m}=\binom{-\frac{3}{5}}{m}.
  \tag{11.1.2}
  $$
By 1.3.1 and (1.3.2), we see that
  $$
  f_{10m}\in\bold Z[\tfrac15].
  \tag{11.1.3}
  $$
The other coefficients  $f_n$  are  $0$.  
Here we are regarding each the coefficient as a specialization 
of  $f_j$  of the Subsection 2.1.  
Then the divided universal Bernoulli number  $\hat{B}_{10m}/(10m)$  
is specialized to  $D\lr{1}_{10m}/(10m)$.   
We claim that
  $$
  \frac{D\lr{1}_{10m}}{10m}\ \in\ 4!\bold Z[\tfrac15]. 
  \tag{11.1.4}
  $$
To prove this, we apply 3.3.1 for  $p=2$  and  $p=3$.  
In this case, as  $f_{2-1}=0$  and  $f_{3-1}=0$, 
we may consider only the partitions  $U$  such that  $U_{2-1}=U_{3-1}=0$.  
If  $m=1$, then  $D\lr{1}_{10}/10=(3/5)(9!/11)=2^73^57/11$  is divisible by  $4!$.  
If  $m\geqq 2$, then because of
  $$
  w(U)=10m\geqq 20, \ \ 
  d(U)\geqq 1
  $$
we have
  $$
  \ord_2(\tau_{\sssize U})
  \geqq\lfloor \tfrac{20+1-2}{4}\rfloor
  =4C\ \ \ 
  \ord_3(\tau_{\sssize U})
  \geqq\lfloor \tfrac{20+1-2}{6}\rfloor
  =3. 
  $$
Hence, thanks to the expression in 2.2.5, 
all the coefficients  $D\lr{1}_{10m}/(10m)$  are divisible by  $4!$.  
This and (11.1.3) yield (11.1.4).  

If  $p=10m+1$  is a prime number, 
then the coefficient  $f_{p-1}$  coincides with  $(2,2)$-entry of 
the Hasse-Witt matrix with respect to the natural basis of the differential
forms of the first kind on  $C$  of reduction modulo  $p$:
  $$
  f_{p-1}=-\binom{-\frac{3}{5}}{\frac{p-1}{10}}\equiv A_p \mod p.
  \tag{11.1.5}
  $$
See also 18.1.

Under the consideration above, 
by applying Clarke's theorem 2.3.1 to the Laurent development of 
the function  $u\mapsto -y(\bk{u})^{1/5}$  at  $u=0$, namely to
  $$
  \frac1s=\sum_{m=0}^{\infty}D\lr{1}_{10m}\frac{{u}^{10m-1}}{(10m)!},
  \tag{11.1.6}
  $$
we have 
\vskip 3pt
\lbox{65}{360}{
  $$
  \frac{1}{4!}\frac{D\lr{1}_{10m}}{10m}
   \in\hskip 5pt
    -\hskip -10pt
      \sum_{\Sb p\equiv 1\hskip -3pt \mod{5} \\ 10m=a(p-1) \endSb}
      \hskip -5pt
  \frac{{(4!a)|_p}^{-1}\hskip -5pt\mod{p^{1+\ord_pa}}}{p^{1+\ord_pa}}{A_p}^a
  +\bold Z[\tfrac{1}{5}].
  \tag{11.1.7} 
  $$
}
\vskip 1pt
\fp

\newpage \fp
{\bf 11.2. Relation between coefficients of  $\bk{y}(\bk{u})^{1/5}$  and of  $\bk{y}(\bk{u})^{2/5}$.} \ \ \ 
We prove here a relation between the Hurwitz coefficients of   $y^{2/5}(\bk{u})$  
and of  $y(\bk{u})^{1/5}$  for  $y^2=x^5-1$.  
  $$ 
  s=\frac{-1}{y(\bk{u})^{1/5}}\ \ (=u+\cdots), 
  $$
we let 
  $$
  \frac{1}{s^2}
  =y(\bk{u})^{2/5}  
  =\sum_{m=0}^{\infty}\frac{D\lr{2}_{10m}}{(10m)_2}\frac{{u}^{10m-2}}{(10m-2)!}
  \ \ \ \ (D\lr{2}_0=1).
  \tag{11.2.1}
  $$
Then
  $$
  \aligned
  \int_0^{u}
  \left(\frac{1}{s^2}-\frac{1}{{u}^2}\right)du 
  &=\int_0^{u}
    \left(
    \sum_{m=1}^{\infty}\frac{D\lr{2}_{10m}}{(10m)_2}\frac{{u}^{10m-2}}{(10m-2)!}
    \right)du \\
  &=\sum_{m=1}^{\infty}\frac{D\lr{2}_{10m}}{(10m)_2}\frac{{u}^{10m-1}}{(10m-1)!}.
  \endaligned
  \tag{11.2.2}
  $$
On the other hand, after differentiating (11.1.1) with respect to  $u$, 
dividing by  $s^2$, we have
  $$
  \frac{1}{s^2}
  =\frac{1}{s^2}\frac{ds}{du}
    +\sum_{m=1}^{\infty}\binom{-\frac{3}{5}}{m}
                           s^{10m-2}\frac{ds}{du}.
  \tag{11.2.3}
  $$
Namely, 
  $$
  \int_0^{u}
   \left(
    \frac{1}{s^2}-\frac{1}{{u}^2}
   \right)du
  =\left(-\frac{1}{s}
   +  \sum_{m=1}^{\infty}\binom{-\frac{3}{5}}{m}
                           \frac{s^{10m-1}}{10m-1}\right)
   + \frac{1}{u}. 
  $$
By equating this with (11.2.2), and by using (11.1.6), we see that
  $$
  \align
  &\sum_{m=1}^{\infty}\frac{D\lr{2}_{10m}}{(10m)_2}\frac{{u}^{10m-1}}{(10m-1)!} \\
  &\hskip 30pt =\left(
  -\sum_{m=0}^{\infty}\frac{D\lr{1}_{10m}}{10m}\frac{{u}^{10m-1}}{(10m-1)!}
  +\sum_{m=1}^{\infty}\binom{-\frac{3}{5}}{m}
                           \frac{s^{10m-1}}{10m-1}
  \right)
  +\frac{1}{u}.
  \endalign
  $$
Since  $D\lr{1}_0=1$, we can remove the terms of negative degree, and we have
  $$
   \sum_{m=1}^{\infty}\frac{D\lr{1}_{10m}}{10m}\frac{{u}^{10m-1}}{(10m-1)!}
  +\sum_{m=1}^{\infty}\frac{D\lr{2}_{10m}}{(10m)_2}\frac{{u}^{10m-1}}{(10m-1)!}
  =\sum_{m=1}^{\infty}\binom{-\frac{3}{5}}{m}
                           \frac{s^{10m-1}}{10m-1}. 
  $$
The right hand side belongs to  $\bold Z[\frac{1}{5}]\llg u\rrg$  
because of 8.1.2 (3) and (8.3.6).  
Adding this with (11.1.3), we conclude that 
  $$
   \frac{D\lr{1}_{10m}}{10m} 
  +\frac{D\lr{2}_{10m}}{(10m)_2} 
      \in (10\hskip -1pt\cdot\hskip -3pt 1\hskip -2pt -\hskip -1pt 2)!\bold Z[\tfrac12]
      \subset 4!\bold Z[\tfrac15].  
  \tag{11.2.4}
  $$

\newpage \fp
{\bf 11.3. Relation between coefficients of  $\bk{y}(\bk{u})^{2/5}$  and of  $\bk{y}(\bk{u})^{3/5}$.} \ \ \ 
We prove here a relation between the Hurwitz coefficients of  $y^{3/5}(\bk{u})$  
and of  $y(\bk{u})^{2/5}$  for  $y^2=x^5-1$.  
For 
  $$ 
  s=\frac{-1}{y(\bk{u})^{1/5}}\ \ (=u+\cdots),
  $$
we let 
  $$
  \frac{1}{s^3}
  =-y(\bk{u})^{3/5}  
  =\sum_{m=0}^{\infty}\frac{D\lr{3}_{10m}}{(10m)_3}\frac{{u}^{10m-3}}{(10m-3)!}
  \ \ \ \ (D\lr{3}_0=1).
  \tag{11.3.1}
  $$
Then 
  $$
  \aligned
  \int_0^{u}
  \left(\frac{1}{s^3}-\frac{1}{{u}^3}\right)du 
  &=\int_0^{u}
    \left(
    \sum_{m=1}^{\infty}\frac{D\lr{3}_{10m}}{(10m)_3}\frac{{u}^{10m-3}}{(10m-2)!}
    \right)du \\
  &=\sum_{m=1}^{\infty}\frac{D\lr{3}_{10m}}{(10m)_3}\frac{{u}^{10m-2}}{(10m-2)!}.
  \endaligned
  \tag{11.3.2}
  $$
On the other hand, after differentiating (11.1.1) with respect to  $u$, 
dividing by  $s^3$, we have
  $$
  \frac{1}{s^3}
  =\frac{1}{s^3}\frac{ds}{du}
    +\sum_{m=1}^{\infty}\binom{-\frac{3}{5}}{m}
                           s^{10m-3}\frac{ds}{du}
  \tag{11.3.3}
  $$
Namely,
  $$
  \int_0^{u}
   \left(
    \frac{1}{s^3}-\frac{1}{{u}^3}
   \right)du
  =\left(-\frac{1}{2}\frac{1}{s^2}
   +  \sum_{m=1}^{\infty}\binom{-\frac{3}{5}}{m}
                           \frac{s^{10m-2}}{10m-2}\right)
   + \frac{1}{2}\frac{1}{{u}^2}. 
  $$
By equating this with (11.3.2), and by using (11.2.1), we see that
  $$
  \align
  &\sum_{m=1}^{\infty}\frac{D\lr{3}_{10m}}{(10m)_3}\frac{{u}^{10m-2}}{(10m-2)!} \\
  &\hskip 30pt =\left(
   -\frac12\sum_{m=0}^{\infty}
                      \frac{D\lr{2}_{10m}}{(10m)_2}\frac{{u}^{10m-2}}{(10m-2)!}
  +\sum_{m=1}^{\infty}\binom{-\frac{3}{5}}{m}
                           \frac{s^{10m-2}}{10m-2}
  \right)
  +\frac12\frac{1}{{u}^2}.
  \endalign
  $$
Since  $D\lr{2}_0=1$, we can remove the terms of negative degree, and, by multiplying  $2$,  
we have
  $$
   \sum_{m=1}^{\infty}\frac{ D\lr{2}_{10m}}{(10m)_2}\frac{{u}^{10m-2}}{(10m-2)!}
  +\sum_{m=1}^{\infty}\frac{2D\lr{3}_{10m}}{(10m)_3}\frac{{u}^{10m-2}}{(10m-2)!}
  =2\sum_{m=1}^{\infty}\binom{-\frac{3}{5}}{m}
                           \frac{s^{10m-2}}{10m-2}.
  $$
The right hand side belongs to  $\bold Z[\frac{1}{5}]\llg u\rrg$  
because of 8.1.2 (3) and (8.3.6).  
This and (11.1.3) show that 
  $$
  \frac{D\lr{2}_{10m}}{(10m)_2}+\frac{2D\lr{3}_{10m}}{(10m)_3} 
      \in (10\hskip -1pt\cdot\hskip -3pt 1\hskip -2pt -\hskip -1pt 3)!\bold Z[\tfrac12]
      \subset 4!\bold Z[\tfrac15].
  \tag{11.3.4}
  $$

\newpage \fp
{\bf 11.4. Relation between coefficients of  $\bk{y}(\bk{u})^{3/5}$  and of  $\bk{y}(\bk{u})^{4/5}$.} \ \ \ 
We prove here a relation between the Hurwitz coefficients of  $y^{4/5}(\bk{u})$  
and of  $y(\bk{u})^{3/5}$  for  $y^2=x^5-1$.  
For 
  $$ 
  s=\frac{-1}{y(\bk{u})^{1/5}}\ \ (=u+\cdots),
  $$
we let 
  $$
  \aligned
  \frac{1}{s^4}
  =y(\bk{u})^{4/5}  
  =\sum_{m=0}^{\infty}\frac{D\lr{4}_{10m}}{(10m)_4}\frac{{u}^{10m-4}}{(10m-4)!}
  \ \ \ \ (D\lr{4}_0=1).
  \endaligned
  \tag{11.4.1}
  $$
Then
  $$
  \aligned
  \int_0^{u}
  \left(\frac{1}{s^4}-\frac{1}{{u}^4}\right)du 
  &=\int_0^{u}
    \left(
    \sum_{m=1}^{\infty}\frac{D\lr{4}_{10m}}{(10m)_4}\frac{{u}^{10m-4}}{(10m-4)!}
    \right)du \\
  &=\sum_{m=1}^{\infty}\frac{D\lr{4}_{10m}}{(10m)_4}\frac{{u}^{10m-3}}{(10m-3)!}
  \endaligned
  \tag{11.4.2}
  $$
On the other hand, after differentiating (11.1.1) with respect to  $u$, 
dividing by  $s^4$, we have
  $$
  \frac{1}{s^4}
  =\frac{1}{s^4}\frac{ds}{du}
    +\sum_{m=1}^{\infty}\binom{-\frac{3}{5}}{m}
                           s^{10m-4}\frac{ds}{du}.
  \tag{11.4.3}
  $$
Namely, 
  $$
  \int_0^{u}
   \left(
    \frac{1}{s^4}-\frac{1}{{u}^4}
   \right)du
  =-\frac{1}{3}\frac{1}{s^3}
   + \frac{1}{3}\frac{1}{{u}^3}
   +  \sum_{m=1}^{\infty}\binom{-\frac{3}{5}}{m}
                           \frac{s^{10m-3}}{10m-3}. 
  $$
By equating this with (11.4.2), and by using (11.3.1), we see that
  $$
  \aligned
  \sum_{m=0}^{\infty}\frac{D\lr{4}_{10m}}{(10m)_4}&\frac{{u}^{10m-3}}{(10m-3)!} \\
  =&\left(\frac{-1}3 \sum_{m=1}^{\infty}
                     \frac{D\lr{3}_{10m}}{(10m)_3}\frac{{u}^{10m-3}}{(10m-3)!}
  +\sum_{m=1}^{\infty}\binom{-\frac{3}{5}}{m}
                           \frac{s^{10m-3}}{10m-3}
   \right)
  +\frac13\frac{1}{{u}^3}.
  \endaligned
  $$
Since  $D\lr{3}_0=1$, we can remove the terms of negative degree, and, by multiplying  $3$,
we have
  $$
   \sum_{m=1}^{\infty}\frac{ D\lr{3}_{10m}}{(10m)_3}\frac{{u}^{10m-3}}{(10m-3)!}
  +\sum_{m=1}^{\infty}\frac{3D\lr{4}_{10m}}{(10m)_4}\frac{{u}^{10m-3}}{(10m-3)!}
  =3\sum_{m=1}^{\infty}\binom{-\frac{3}{5}}{m}
                           \frac{s^{10m-3}}{10m-3}. 
  $$
The right hand side of this belongs to  $\bold Z[\frac{1}{5}]\llg u\rrg$  
because of 8.1.2 (3) and (8.3.6).  
This and (11.1.3) show that  
  $$
  \frac{D\lr{3}_{10m}}{(10m)_3}+\frac{3D\lr{4}_{10m}}{(10m)_4} 
      \in (10\hskip -1pt\cdot\hskip -3pt 1\hskip -2pt -\hskip -1pt 4)!\bold Z[\tfrac12]
      \subset 4!\bold Z[\tfrac15].  
  \tag{11.4.4}
  $$

\newpage \fp
{\bf 11.5. Relation between coefficients of  $\bk{y}(\bk{u})^{4/5}$  and of  $\bk{y}(\bk{u})$.} \ \ \ 
We prove here a relation between the Hurwitz coefficients of  $y(\bk{u})$  
and of  $y(\bk{u})^{4/5}$  for  $y^2=x^5-1$.  
For
  $$ 
  s=\frac{-1}{y(\bk{u})^{1/5}}\ \ (=u+\cdots)
  $$
we let 
  $$
  \frac{1}{s^5}
  =-y(\bk{u})  
  =\sum_{m=0}^{\infty}\frac{D\lr{5}_{10m}}{(10m)_5}\frac{{u}^{10m-5}}{(10m-5)!}
  \ \ \ \ (D\lr{5}_0=1).
  \tag{11.5.1}
  $$
Then
  $$
  \aligned
  \int_0^{u}
  \left(\frac{1}{s^5}-\frac{1}{{u}^5}\right)du 
  &=\int_0^{u}
    \left(
    \sum_{m=1}^{\infty}\frac{D\lr{5}_{10m}}{(10m)_5}\frac{{u}^{10m-5}}{(10m-5)!}
    \right)du \\
  &=\sum_{m=1}^{\infty}\frac{D\lr{5}_{10m}}{(10m)_5}\frac{{u}^{10m-4}}{(10m-4)!}.
  \endaligned
  \tag{11.5.2}
  $$
On the other hand, after differentiating (11.1.1) with respect to  $u$, dividing by  $s^5$,
we obtain
  $$
  \frac{1}{s^5}
  =\frac{1}{s^5}\frac{ds}{du}
    +\sum_{m=1}^{\infty}\binom{-\frac{3}{5}}{m}
                           s^{10m-5}\frac{ds}{du}.
  \tag{11.5.3}
  $$
Namely, 
  $$
  \int_0^{u}
   \left(
    \frac{1}{s^5}-\frac{1}{{u}^5}
   \right)du
  =-\frac{1}{4}\frac{1}{s^4}
   + \frac{1}{4}\frac{1}{{u}^4}
   +  \sum_{m=1}^{\infty}\binom{-\frac{3}{5}}{m}
                           \frac{s^{10m-4}}{10m-4}.
  $$
By equating this with (11.5.2), and by using (11.4.1), we see that
  $$
  \aligned
  \sum_{m=1}^{\infty}\frac{D\lr{5}_{10m}}{(10m)_5}&\frac{{u}^{10m-4}}{(10m-4)!} \\
  =&\left(\frac{-1}4\sum_{m=0}^{\infty}
                     \frac{D\lr{4}_{10m}}{(10m)_4} \frac{{u}^{10m-4}}{(10m-4)!}
  +\sum_{m=1}^{\infty}\binom{-\frac{3}{5}}{m}
                           \frac{s^{10m-4}}{10m-4}
  \right)
  +\frac14\frac{1}{{u}^4}.
  \endaligned
  $$
Since  $D\lr{4}_0=1$, we can remove the terms of negative degree, and, by multiplying  $4$, 
we have
  $$
   \sum_{m=1}^{\infty}\frac{ D\lr{4}_{10m}}{(10m)_4}\frac{{u}^{10m-4}}{(10m-4)!}
  +\sum_{m=1}^{\infty}\frac{4D\lr{5}_{10m}}{(10m)_5}\frac{{u}^{10m-4}}{(10m-4)!}
  =4\sum_{m=1}^{\infty}\binom{-\frac{3}{5}}{m}
                           \frac{s^{10m-4}}{10m-4}.
  $$
The right hand side of this belongs to  $\bold Z[\frac{1}{5}]\llg u\rrg$  
because of 8.1.2 (3) and (8.3.6).  
This and (11.1.3) shows that
  $$
  \frac{D\lr{4}_{10m}}{(10m)_4}+\frac{4D\lr{5}_{10m}}{(10m)_5} 
      \in (10\hskip -1pt\cdot\hskip -3pt 1\hskip -2pt -\hskip -1pt 5)!\bold Z[\tfrac12]
      \subset 4!\bold Z[\tfrac15]
  \tag{11.5.4}
  $$

\newpage \input 11-6.o
\newpage \fp
\inbox{\bf 12}\ \ 
{\bf  The Clarke type theorem on  $\bk{x}(\bk{u})$  and  $\bk{y}(\bk{u})$  over  $\bold Z$.} \ \ 
\vskip 5pt
\fp
{\bf 12.1. Relation between coefficients of  $\bk{y}(\bk{u})$  and of  $\bk{x}^2(\bk{u})$.}\ \ \
Let 
  $$
  x^2(\bk{u})=\frac1{t^4}
  =\sum_{m=0}^{\infty}
   \frac{C\lr{4}_{10m}}{(10m)_4}\frac{{u}^{10m-4}}{(10m-4)!}.
  \tag{12.1.1}
  $$
If  $D=d/du=(2y/x)d/dx$, then
  $$
  D(x^2)=2xDx=2x\tfrac{2y}{x}=4y.
  \tag{12.1.2}
  $$
Hence
  $$
  \frac{C\lr{4}_{10m}}{(10m)_4}=4\frac{D_{10m}}{10m}. 
  \tag{12.1.3}
  $$
Therefore, (10.5.3) and (11.6.3) show that 
the denominator of  $C\lr{4}_{10m}/(10m)_4$  does not contain any power of 2, 
and that  $C\lr{4}_{10m}/(10m)_4$  has the property that
\vskip 5pt
\fp
\lbox{65}{360}{
  $$
  \frac{C\lr{4}_{10m}}{(10m)_4}\in\hskip 5pt
  -\hskip -10pt\sum_{\Sb p\equiv 1\hskip -3pt\mod{5} \\ 10m=a(p-1) \endSb}
  \hskip -5pt
   \frac{{(3!a)|_p}^{-1}\hskip -5pt\mod{p^{1+\ord_pa}}}{p^{1+\ord_pa}}{A_p}^a
  +\bold Z.
  \tag{12.1.4}
  $$
}
\vskip 3pt
\fp
Furthermore, since  $\frac{1}{4}\frac{C\lr{4}_{10m}}{(10m)_4}=\frac{D_{10m}}{10m}$  by (12.1.3), 
we see that
  $$
  \frac{D_{10m}}{10m}\in
  \bold Z\Big[\frac15,\  \frac1p\ ; \ 
  p\equiv 1\ \text{mod}\ 5, \ p-1|10m\Big]
  \subset\bold Z_{(2)},
  $$
and that the denominator of  $D_{10m}/(10m)$  does not contain any power of  $2$.  
Thus the numerator of  $C\lr{4}_{10m}/(10m)_4$  is divisible by  $4$.  
This consideration and (11.6.3) give rise to
  $$
  \frac{D_{10m}}{10m}
  \in\hskip 5pt
  -\hskip -10pt\sum_{\Sb p\equiv 1\hskip -3pt\mod{5}\\ 10m=a(p-1) \endSb}
   \hskip -10pt
   \frac{{(4!a)|_p}^{-1}\hskip -5pt\mod{p^{1+\ord_pa}}}{p^{1+\ord_pa}}{A_p}^a
   +\bold Z.
  \tag{12.1.5}
  $$
This is the second formula in 6.1.1.

\newpage \input 12-2.o
\newpage \fp
\inbox{\bf 13}\ \ 
{\bf Proof of the Kummer type congruence relation (Yasuda)}\ \ 
\vskip 5pt
\fp
Now we start to prove 7.1.1 (and 7.2.1) following Yasuda \cite{{\bf Ya2}}.
\vskip 3pt
\fp
{\bf 13.1. 13.1. Honda's theorem and formal groups.}\ \ \ 
We show that there exists formal group over  $\bold Z$  
such that its formal logarithm is (10.1.1) (or (11.1.1)).  
\vskip 5pt
\fp
\lbox{151}{360}{{\bf Proposition 13.1.1.}\ \it  
Suppose  $p$  be a fixed prime number. 
Let  $t$  be an indeterminate 
and let  $u(t)\in \bold Q_p[[t]]$.  
If there exists  $\beta\in\bold Z_p$  such that
  $$
  pu(t)-\beta u(t^p)\in p\bold Z_p[[t]]
  \tag{13.1.2}
  $$
then the formal group law  $F$  given by 
  $$
  F(t_1, t_2):=u^{-1}(u(t_1)+u(t_2))
  \tag{13.1.3}
  $$
is defined over  $\bold Z_p$ {\rm (}i.e. $\in\bold Z_p[[t_1, t_2]]${\rm )}. 
Here  $u^{-1}(t)\in\bold Q_p[[t]]$  means the unique power series 
such that  $u^{-1}(u(t))=t$.  
In particular, if  $\alpha\in\bold Z_p$, then we have
  $$
  F(\alpha t)=u^{-1}(\alpha u(t))\in \alpha t+t^2\bold Z_p[[t]].
  \tag{13.1.4}
  $$
}
\vskip 1pt
\fp
This is obtained by applying Honda's theorem \cite{{\bf Ho}, p.223, Theorem 2} as  $n=1$,  
$q=p$, $P=1$, $u=p-\beta T$, $f=u(t)$.   
We can apply this theorem to  $t\mapsto u=u_g$  of (10.1.1) 
and  $s\mapsto u$  of  (11.1.1).  
To do so, we introduce the  $p$-adic  $\varGamma$-function
  $$
  \varGamma_p:\bold Z_p \rightarrow {\bold Z_p}^{\times}.
  \tag{13.1.5}
  $$
This function satisfies, for any positive integer  $n$, that
  $$
  \varGamma_p(n)
  =(-1)^n\prod_{\Sb 1\leqq j< n \\ p\not\hskip 2pt |\hskip 1pt j\endSb}j.
  \tag{13.1.6}
  $$
The most important properties are that
  $$
  \varGamma_p(z+1)
  =\cases
  -z\varGamma_p(z) &  (z\not\in p\bold Z_p)\\
  -\varGamma_p(z) &  (z    \in p\bold Z_p)
  \endcases
  \tag{13.1.7}
  $$
and that, for any positive integer  $\nu$, 
  $$
  z\equiv w \mod{p^{\nu}\bold Z_p}\ \ \text{implies}\ \ \ 
  \varGamma_p(z)\equiv \varGamma_p(w) \mod{p^{\nu}\bold Z_p}.
  \tag{13.1.8}
  $$
For the details, see \cite{{\bf Mo}} or \cite{{\bf R}}.  

Denoting by  $u(t)$  the power series development of
  $$
  u(=u_g)=\int_{\infty}^{(x,y)}\frac{x^{g-1}dx}{2y}
  $$
with respect to  $t=-x^{-\frac{1}{2}}$  we claim that
\vskip 3pt
\fp
\lbox{166}{360}
{{\bf Lemma 13.1.9.}\ \it 
For the curve  $y^2=x^{2g+1}-1$ {\rm (}resp. $y^2=x^{2g+1}-x$ {\rm )} 
and a prime  $p\equiv 1\mod{2g+1}$ {\rm (}resp. $p\equiv 1\mod{4g}$ {\rm )}, 
  $$
  \align
  \beta_p&=-(-1)^{(p-1)/(4g+2)}
  \frac{\varGamma_p(\frac{1}{2})}
       {\varGamma_p(\frac{4g+1}{4g+2})                                 
        \varGamma_p(\frac{2g+2}{4g+2})}\ \ (\ \in {\bold Z_p}^{\times}\ ). \\
  \bigg(\text{resp.} \ \ 
  \beta_p&=\ -(-1)^{(p-1)/(4g)}\ 
  \frac{\varGamma_p(\frac{1}{2})}
       {\varGamma_p(\frac{4g-1}{4g})                                 
        \varGamma_p(\frac{2g+1}{4g})}\hskip 11pt  (\ \in {\bold Z_p}^{\times}\ ).
  \bigg)
  \endalign
  $$
Then the series satisfies
  $$
  pu(t)-\beta_p u(t^p)\in p\bold Z_p[[t]].
  \tag{13.1.10}
  $$
The similar formula is satisfied by the power series development  $s=-y^{-\frac{1}{2g+1}}$  
with respect to  $u$. 
}  
\fp
{\it Proof.} \ 
For simplicity, we prove the statement only for the curve  $y^2=x^5-1$ ($g=2$).  
Let  $f_{10n}/(10n+1)=[t^{10n+1}]u(t)$.  
Then, as is stated in (10.1.2), $f_{10n}=(-1)^n\binom{-\frac{1}{2}}{n}$.  
It suffice for us to prove that, if  $p(10m+1)=10n+1$, then
  $$
  p\frac{f_{10n}}{10n+1}-\beta_p\frac{f_{10m}}{10m+1}\in p \bold Z_p.
  $$
We notice that  $\big\lceil\tfrac{n}{p}\big\rceil=\big\lceil m+\tfrac{1}{10}-\tfrac{1}{10p}\big\rceil=m$.  
This and 
  $$
  p(2\big\lceil\tfrac{n}{p}\big\rceil-1)
  \leqq  2n-1 
  \leqq  p(2\big\lfloor\tfrac{n}{p}\big\rfloor-1)
  \leqq  p\big(2(\big\lceil\tfrac{n}{p}\big\rceil+1)-1\big)
  $$
show that the largest odd integer divisible by  $p$  
not exceed  $2n-1$  is  $p(2m-1)$.  
Hence
  $$
  \allowdisplaybreaks
  \align
  p&\frac{f_{10n}}{10n+1}-\beta_p\frac{f_{10m}}{10m+1} \\
  &=p (-1)^n\binom{-\frac{1}{2}}{n}\frac{1}{10n+1}  
     - \beta_p\cdot (-1)^m\binom{-\frac{1}{2}}{m}\frac{1}{10m+1} \\
  &=\frac{1}{10m+1}
     \Bigg\{(-1)^{pm+\frac{p-1}{10}}
       \frac{\prod_{j=1}^n\big(\hskip -2pt -\tfrac{2j-1}{2}\big)}{n!}
      -\beta_p\cdot (-1)^m
       \frac{\prod_{j=1}^m\big(\hskip -2pt -\frac{2j-1}{2}\big)}{m!}  
     \Bigg\}  \\
  &=\frac{(-1)^{m+\frac{p-1}{10}}}{10m+1}
     \Bigg\{
      \frac{\prod_{j=1,\hskip 1pt p\hskip -1pt\not\hskip 1.8pt|\hskip .5pt 2j-1}^n
                   \big(\hskip -3pt -\frac{2j-1}{2}\big)
             \hskip 3pt
            \prod_{k=1}^m           
                   \big(\hskip -3pt -\frac{p(2k-1)}{2}\big)
           }
           {\big(\prod_{j=1,\hskip 1pt p\hskip -1pt\not\hskip 1.7pt|j}^n j\big)
            \big(\prod_{k=1}^m pk\big)
           } \\
  & \hskip 185pt
     -(-1)^{\frac{p-1}{10}}\beta_p 
      \frac{\prod_{j=1}^m\big(\hskip -3pt -\frac{2j-1}{2}\big)}{m!}  
     \Bigg\}.
  \endalign
  $$
By using (13.1.7) repeatedly, we have
  $$
  \align
  &=\frac{(-1)^{m+\frac{p-1}{10}}}{10m+1}
     \Bigg\{
      \frac{
       \frac{(-1)^n\varGamma_p(-\frac{1}{2}+1)}
            {      \varGamma_p(-\frac{2n-1}{2})}
               \hskip 4pt
             p^m\prod_{k=1}^m     \big(\hskip -3pt -\frac{2k-1}{2}\big)
             }{(-1)^{n+1}\varGamma_p(n+1)\hskip 3pt p^m m!} \\
  \nomathbreak
  & \hskip 195pt  -(-1)^{\frac{p-1}{10}}\beta_p
      \frac{\prod_{j=1}^m\big(\hskip -3pt -\frac{2j-1}{2}\big)}{m!} 
     \Bigg\} \\
  \allowdisplaybreak
  &=\frac{(-1)^{m+\frac{p-1}{10}}}{10m+1}
     \Bigg\{
      -\frac{\varGamma_p(\frac{1}{2})}
            {\varGamma_p(n+1)\varGamma_p(-n+\frac{1}{2})}
       \frac{\prod_{k=1}^m   \big(\hskip -3pt -\frac{2k-1}{2}\big)}
            {m!} \\
  & \hskip 180pt  
     +\frac{\varGamma_p(\frac 12)}
                          {\varGamma_p(\frac{9}{10})\varGamma_p(\frac{6}{10})}
      \frac{\prod_{j=1}^m\big(\hskip -3pt-\frac{2j-1}{2}\big)}{m!} 
     \Bigg\}  \\
  &=\frac{(-1)^{m+\frac{p-1}{10}}}{10m+1}\binom{-\frac{1}{2}}{m}
     \Bigg\{
     -\frac{\varGamma_p(\frac{1}{2})}
           {\varGamma_p(-\frac{1}{10}p(10m+1)+\frac{9}{10})
            \varGamma_p(-\frac{1}{10}p(10m+1)+\frac{6}{10})} \\
  &  \hskip 230pt
     +\frac{\varGamma_p(\frac 12)}
                          {\varGamma_p(\frac{9}{10})\varGamma_p(\frac{6}{10})}
     \Bigg\}.
    \endalign
  $$
Assume  $p^{\nu}||(10m+1)$.  
Then (13.1.7) and (13.1.8) show
  $$
     \frac{\varGamma_p(\frac{1}{2})}
           {\varGamma_p(-\frac{1}{10}p(10m+1)+\frac{9}{10})
            \varGamma_p(-\frac{1}{10}p(10m+1)+\frac{6}{10})} 
     \equiv 
     \frac{\varGamma_p(\frac 12)}
           {\varGamma_p(\frac{9}{10})\varGamma_p(\frac{6}{10})} \mod{p^{\nu+1}\bold Z_p}.
  $$
It is shown by 1.3.3 that
  $$
  \binom{-\frac{1}{2}}{m}\in \bold Z_p.
  $$
Therefore
  $$
  p\frac{f_{10n}}{10n+1}-\beta_p\frac{f_{10m}}{10m+1}
  \in p\bold Z_p.
  $$
and the proof has completed. 
\qed
\vskip 5pt
\fp
{\bf Remark 13.1.11.}\ 
(1) (Yasuda) \   
Regarding the map  $\iota:C\rightarrow J$  to be defined over  $\bold Z_p$, 
we denote its formal completion at  $\infty\mapsto$ "the origin of $J$"  
by  $\hat{\iota}:\hat{C} \rightarrow \hat{J}$.  
Because of the action (5.2.2) (resp. (5.3.2)), 
$\hat{J}$  is decomposed into a product of $1$-dimensional formal groups.  
The composite map  $\hat{\pi}\circ\hat{\iota}:\hat{C}\rightarrow G$  of  $\hat{\iota}$  
with the projection  $\pi:\hat{J}\rightarrow G$, 
where  $G$  is a certain factor of the product, would be an isomorphism of formal groups. 
\fp
(2) The height of the formal group  $F$  over  $\bold Z_p$  associated to  $u(t)$  above is  $1$, 
namely, 
there exists  $b\in\bold Z_p^{\times}$  such that  $u^{-1}(pu(t))\equiv bt^p$  mod $p\bold Z_p[[t]]$.  
This fact follows from
  $$
  pu(t)=p\big(t+\cdots+f_{p-1}\jrac{t^p}{p}+\cdots\big)
  $$
and the fact  $f_{p-1}\in\bold Z_p^{\times}$  shown by (10.1.5).  
The same fact is seen for the formal group associated 
to the series development of  $s\mapsto u$.

\newpage \fp
{\bf 13.2. Hochschilt's formula and Honda's theorem.} \ \ \ 
We recall the following Lemma called Hochschilt's formula in order to show 13.2.2 below.
\vskip 3pt
\fp
\lbox{71}{360}{{\bf Proposition 13.2.1.}\ \ \it
Let  $p$  be a prime, $R$  be a commutative ring of characteristic  $p$,   
and  $D$  be a derivation from  $R$  to itself.  
Then, for  $b\in R$, we have
  $$
  (b D)^p=b^pD^p+\big((bD)^{p-1}(b)\big)\hskip -2pt\cdot\hskip -2pt D.
  $$
}
\vskip 3pt
\fp
Since this is described in \cite{{\bf Ma}, p.197, Theorem 25.5} with a detailed proof, 
we omit it.  

While the following general equality is proved by using the formula above, 
rather weak form for the case of  $g=1$  is described in \cite{{\bf G}} by a different way.  
\vskip 5pt
\fp
\lbox{90}{360}{{\bf Proposition 13.2.2.}\ \ \it
For the curve  $y^2=x^{2g+1}-1$  {\rm (}resp.  $y^2=x^{2g+1}-x$ {\rm )}, 
let  $t=t(u)$  and  $s=s(u)$  be the power series of \hskip 2pt {\rm 8.2.1} 
and of \hskip 2pt {\rm 8.3.1}, respectively.  
Let  $p\equiv 1$ {\rm mod} $(2g + 1)$ {\rm (}resp. {\rm mod} $4g${\rm )}  be a prime. 
If  $\varphi\in \bold Z_p[[t]]$  or  $\varphi\in \bold Z_p[[s]]$  then 
  $$
  \Big(\Big(\jrac{d}{du}\Big)^p-A_p\jrac{d}{du}\Big)\varphi\in p\bold Z_p[[t]]
  \ \ \ \text{\rm or}\ \ \in p\bold Z_p[[s]]. 
  $$
}
\vskip 3pt
\fp
{\it Proof. } \ 
The derivation  $D=\jrac{d}{dt}$  on the ring  $\bold Z_p[[t]]$  induces a derivation 
from  $\bold F_p[[t]]$  to itself.  
In the sequel of this proof, the symbol  \lq\lq$=$"  means the equality 
in  $\bold F_p[[t]]$.  
First of all, we pay attention to the fact
  $$
  \aligned
  &\jrac{du}{dt}=1+\sum_{n=1}^{\infty}(-1)^n\binom{-\frac12}{n}t^{10n}
                \in 1+t\hskip 1pt\bold Z[\tfrac{1}{2}][[t]], \\
  &\jrac{dt}{du}\in 1+t\hskip 1pt\bold Z[\tfrac{1}{2}][[t]]
  \endaligned
  \tag{13.2.3}
  $$
given by (10.1.1) and (10.1.3).  
By using 13.2.1 for this  $b:=\frac{dt}{du}$  and  $D$  above, we see that
  $$
  \align
  0
  &=\Big(\jrac{d}{du}\Big)^{\hskip -2pt p}u \\
  &=(bD)^pu \\ 
  &=\big(b^pD^p+((bD)^{p-1}(b))\cdot D\big)u \\
  &=\Big(\jrac{dt}{du}\Big)^{\hskip -2pt p}\jrac{d\sz{p}u}{dt\sz{p}} 
   +\Big\{\Big(\jrac{dt}{du}\jrac{d}{dt}\Big)^{\hskip -2pt p-1}\jrac{dt}{du}\Big\}\jrac{du}{dt} \\
  &=\Big(\jrac{dt}{du}\Big)^{\hskip -2pt p}\jrac{d\sz{p}u}{dt\sz{p}}
   +\Big\{\Big(\jrac{d}{du}\Big)^{\hskip -2pt p-1}\jrac{dt}{du}\Big\}\jrac{du}{dt} \\
  &=\Big(\jrac{dt}{du}\Big)^{\hskip -2pt p}\jrac{d\sz{p}u}{dt\sz{p}}
   +\jrac{d\sz{p}t}{du\sz{p}}\jrac{du}{dt}. 
  \endalign
  $$
Hence 
  $$
  \jrac{d\sz{p}t}{du\sz{p}}=-\Big(\jrac{du}{dt}\Big)^{\hskip -2pt -p-1}\jrac{d\sz{p}u}{dt\sz{p}}.
  \tag{13.2.4}
  $$
By using 13.2.1 again for  $\varphi$, $D$, and  $b$, we have
  $$
  \allowdisplaybreaks
  \align
  \Big(\jrac{d}{du}\Big)^{\hskip -2pt p}\varphi
  &=(bD)^p\varphi \\
  &=\big(b^pD^p+\big((bD)^{p-1}(b)\big)\cdot D\big)\varphi \\
  &=\Big(\jrac{dt}{du}\Big)^{\hskip -2pt p}\Big(\jrac{d}{dt}\Big)^{\hskip -2pt p}\varphi 
    +\Big(\Big(\jrac{dt}{du}\jrac{d}{dt}\Big)^{\hskip -2pt p-1}\jrac{dt}{du}\Big)\jrac{d}{dt}\varphi \\
  &=\Big(\jrac{du}{dt}\Big)^{\hskip -2pt -p}\jrac{d\sz{p}}{dt\sz{p}}\varphi
    +\jrac{d\sz{p}t}{du\sz{p}}\jrac{du}{dt}\jrac{d}{du}\varphi \\
  &=\jrac{d\sz{p}t}{du\sz{p}}\jrac{du}{dt}\jrac{d}{dt}\varphi.  
  \endalign
  $$
Here we have used (13.2.3). By substituting (13.2.4) into the last formula, we obtain
  $$
  \Big(\jrac{d}{du}\Big)^{\hskip -2pt p}\varphi
  =-\Big(\jrac{du}{dt}\Big)^{\hskip -3pt -p}\jrac{d\sz{p}u}{dt\sz{p}}\cdot\jrac{d}{du}\varphi.
  \tag{13.2.5}
  $$
Since  $(\frac{d}{dt})^p\bold Z_p[[t]]\in p\bold Z_p[t]]$, 
by applying  $(\frac{d}{dt})^p$  to (13.1.10), we have
  $$
  \align
  \jrac{d\sz{p}u}{dt\sz{p}}
  &=\Big(\jrac{d}{dt}\Big)^{\hskip -2pt p-1}\big(\beta_pt^{p-1}u'(t^p)\big) \\
  &=(p-1)!\ u'(t^p) \\
  &=-\beta_pu'(t)^p.
  \endalign
  $$
The properties (13.1.6) and (13.1.8) imply
  $$
  \align
  A_p&=(-1)^{(p-1)/(4g+2)}\binom{\frac{p-1}{2}}{\frac{p-1}{4g+2}} \\
     &=(-1)^{(p-1)/(4g+2)}(-1)
       \frac{\varGamma_p(\frac{p-1}{2}+1)}
            {\varGamma_p(\frac{p-1}{4g+2}+1)
             \varGamma_p(\frac{p-1}{2}-\frac{p-1}{4g+2}+1)} \\
     &=-(-1)^{(p-1)/(4g+2)}
       \frac{\varGamma_p(\frac{p+1}{2})}
            {\varGamma_p(\frac{p+4g+1}{4g+2})
             \varGamma_p(\frac{2gp+2g+2}{4g+2})} \\
     &\equiv -(-1)^{(p-1)/(4g+2)}
       \frac{\varGamma_p(\frac{1}{2})}
            {\varGamma_p(\frac{4g+1}{4g+2})
             \varGamma_p(\frac{2g+2}{4g+2})} \mod{p} \\
     &=\beta_p.
  \endalign
  $$
So that 
  $$
  \jrac{d\sz{p}u}{dt\sz{p}}
  =-A_p\Big(\jrac{du}{dt}\Big)^{\hskip -2pt p}.
  $$
To Substitute this into (13.2.5) shows
  $$
  \Big(\jrac{d}{du}\Big)^{\hskip -2pt p}\varphi
  =A_p\jrac{d}{du}\varphi
  $$ 
as desired. 
We can prove the statement for s = s(u) by the same argument.
\vskip 5pt
\fp
{\bf Remark 13.2.6.} \ 
When the proof of Yasuda was not yet given, we used in order
to prove 13.2.2 a certain weaker formula than 13.2.1 combined with a method of
Carlitz in (\cite{{\bf Ca1}}). 
For convenience to the reader, the weaker formula is mentioned in 18.3.

\newpage \fp
{\bf 13.3. Proof of the Kummer type congruence relation.}\ \ \ 
The tools for the proof of 7.1.1 (and 7.2.1) have been completely prepared. 
To avoid confusion, if we regard
  $$
  x(\bk{u})=\frac{1}{{u}^2}
       +\sum_{n=2}^{\infty}\frac{C_{10n}}{10n}
                           \frac{{u}^{10n-2}}{(10n-2)!}
  $$
to be an element in  $\bold Q_p[[u]]$, we denote this by
  $$
  x\langle u\rangle (\hskip -2pt :=x(\bk{u}))
  $$
Let   $\xi\in \bold Z_p$  be a primitive  $(p-1)$-st root of  $1$.  
Then, as
  $$
  x\langle u\rangle-\xi^2x\langle \xi u\rangle
  =\sum_{\Sb n=1 \\ p-1\hskip -1pt \not\hskip 2pt |10n \endSb}^{\infty}
     \frac{(1-\xi^{10n})C_{10n}}{10n}
     \frac{{u}^{10n-2}}{(10n-2)!}
  \tag{13.3.1}
  $$
it suffice to prove for any positive integer  $a$  and  $D=d/du$  that
  $$
  (D^p-A_pD)^a(x\langle u\rangle-\xi^2x\langle \xi u\rangle)
  \in p^a\bold Z_p\llg u\rrg. 
  \tag{13.3.2}
  $$
Indeed, if  $10n\geqq a+2$  and  $(p-1)\hskip -3pt\not|\hskip 2pt 10n$, 
then the coefficient of  $u^{10n-a-2}/(10n-a-2)!$  is
  $$
  (1-\xi^{10n})\sum_{r=0}^a\binom{r}{a}(-A_p)^{a-r}\frac{C_{10n+r(p-1)}}{10n+r(p-1)}
  $$
and  $1-\xi^{10n}\not\in p\bold Z_p$.  
Thanks to 13.2.2, it suffice for proving (13.3.2) to prove
  $$
  x\langle u(t)\rangle-\xi^2x\langle \xi u(t)\rangle
  \in \bold Z_p[[t]].
  \tag{13.3.3}
  $$
If we set
  $$
  \aligned
  F_{\xi}(t)&={u}^{-1}(\xi u(t))\ \ \
  (\text{${u}^{-1}$  is the formal inverse series of  $t \mapsto u$}) \\
  &=t(\xi u(t)),
  \endaligned
  \tag{13.3.4}
  $$
then 15.3.4 and 15.3.9 yield that
  $$
  x\langle u(t)\rangle-\xi^2x\langle \xi u(t)\rangle
  =\frac{1}{t^2}-\frac{\xi^2}{t(\xi u(t))^2}
  =\frac{1}{t^2}-\frac{\xi^2}{(\xi t+\cdots)^2}
  \in \bold Z_p[[t]]
  $$
because  $F_{\xi}(t)\in \xi t+ t^2\bold Z_p[[t]]$. 
Hence, 7.1.1 have been proved.  
\qed

\newpage \fp
\inbox{\bf 14}\ \ 
{\bf Other Kummer type congruence relations.}\ \ 
\vskip 5pt
\fp
In this Section we do not restrict to the curve  $y^2=x^5-1$, 
and describe hyperelliptic curves defined by either equation of
  $$
  y^2=x^{2g+1}-1, \hskip 30pt 
  y^2=x^{2g+1}-x.
  $$
\vskip 8pt
\fp
{\bf 14.1. On generalized Bernoulli-Hurwitz numbers of higher order.} \ 
Let  $t=-1/x(\bk{u})^{1/2}$  and  $s=-1/y(\bk{u})^{1/5}$  as usual.  
We prove the Kummer type congruence relations 
for the Hurwitz coefficients of  $t^{-\nu}$  ($1\leqq\nu\leqq 4g+2$) 
generalizing the congruence for   $t^{-2}=x(\bk{u})$  proved in the Section 13.  
\vskip 6pt
\fp
{\bf The case of  $\bk{y}(\bk{u})^2=\bk{x}(\bk{u})^{2g+1}-{\bold 1}$.} \ \ 
For the curve  $y^2=x^{2g+1}-1$, 
we recall the numbers  $C\lr{\nu}_{(4g+2)n}$ (defined in the Section 10) 
and  $D\lr{\nu}_{(4g+2)n}$ (defined in the Section 11). 
Namely, for $\nu=1$, $2$, $\cdots$, we let
  $$
  \aligned
  \frac{1}{t^{\nu}}
  &=\frac{1}{{u}^{\nu}}
  +\sum_{n=1}^{\infty}\frac{C\lr{\nu}_{(4g+2)n}}{((4g+2)n)_{\nu}}\frac{{u}^{(4g+2)n-\nu}}{((4g+2)n-\nu)!}, \\
  \frac{1}{s^{\nu}}
  &=\frac{1}{{u}^{\nu}}
  +\sum_{n=1}^{\infty}\frac{D\lr{\nu}_{(4g+2)n}}{((4g+2)n)_{\nu}}\frac{{u}^{(4g+2)n-\nu}}{((4g+2)n-\nu)!}.
  \endaligned
  \tag{14.1.1}
  $$
Of course, if  $(4g+2)\hskip -2pt\not|\hskip 2pt n$  then we assume  $C\lr{\nu}_n=D\lr{\nu}_n=0$.  

Let  $p\equiv 1$ mod $(4g+2)$  be a prime number, 
and  $\xi\in\bold Z_p$  be a primitive  $(p-1)$-st root of $1$.  
By using
  $$
  \frac{1}{t(u)^{\nu}}-\xi^{\nu}\frac{1}{t(\xi u)^{\nu}}\ 
  =\hskip -10pt
   \sum_{\Sb n=1 \\ p-1\hskip -1pt \not\hskip 2pt |(4g+2)n \endSb}^{\infty}
   \hskip -10pt
     \frac{(1-\xi^{(4g+2)n})C\lr{\nu}_{(4g+2)n}}{((4g+2)n)_{\nu}}
     \frac{{u}^{(4g+2)n-\nu}}{((4g+2)n-\nu)!}
  \tag{14.1.2}
  $$
instead of (13.3.1), we have the following.  
\vskip 5pt
\fp
\lbox{125}{360}{{\bf Theorem 14.1.3.}\ \ \it
Let  $p\equiv 1$ {\rm mod} $(2g+1)$  be a prime number, $a$  be a positive integer.  
Let  $\nu$  be an integer such that  $1\leqq \nu\leqq 4g+2$.  
Assume  $(4g+2)n\geqq a+\nu$  and  $p-1\hskip -3pt\not|\hskip 2pt (4g+2)n$.  
Then we have
  $$
  \align
  \sum_{r=0}^a\binom{r}{a}(-A_p)^{a-r}
  \frac{C\lr{\nu}_{(4g+2)n+r(p-1)}}{((4g+2)n+r(p-1))_{\nu}}
  &\equiv 0 \mod{p^a}, \\
  \sum_{r=0}^a\binom{r}{a}(-A_p)^{a-r}
  \frac{D\lr{\nu}_{(4g+2)n+r(p-1)}}{((4g+2)n+r(p-1))_{\nu}}
  &\equiv 0 \mod{p^a}.
  \endalign
  $$
}
\vskip 20pt
\fp

\newpage

\fp
{\bf The case of  $\bk{y}(\bk{u})^2=\bk{x}(\bk{u})^{2g+1}-\bk{x}(\bk{u})$ .} \ \ 
For the curve  $y^2=x^{2g+1}-x$  and  $\nu=1$, $2$, $\cdots$, 
we define  $C\lr{\nu}_{4gn}$  and  $D\lr{\nu}_{4gn}$  by
  $$
  \aligned
  \frac{1}{t^{\nu}}
  &=\frac{1}{{u}^{\nu}}
  +\sum_{n=1}^{\infty}\frac{C\lr{\nu}_{4gn}}{(4gn)_{\nu}}\frac{{u}^{4gn-\nu}}{(4gn-\nu)!}, \\
  \frac{1}{s^{\nu}}
  &=\frac{1}{{u}^{\nu}}
  +\sum_{n=1}^{\infty}\frac{D\lr{\nu}_{4gn}}{(4gn)_{\nu}}\frac{{u}^{4gn-\nu}}{(4gn-\nu)!}.
  \endaligned
  \tag{14.1.4}
  $$
We also assume that, if  $(4g+2)\hskip -2pt\not|\hskip 2pt n$  then  $C\lr{\nu}_n=D\lr{\nu}_n=0$.  

Let  $p\equiv 1$ mod $4g$  be a prime, 
and  $\xi\in\bold Z_p$  be a primitive  $(p-1)$-st root of  $1$. 
By using
  $$
  \frac{1}{t(u)^{\nu}}-\xi^{\nu}\frac{1}{t(\xi u)^{\nu}}\ 
  =\hskip -10pt
   \sum_{\Sb n=1 \\ p-1\hskip -1pt \not\hskip 2pt |4gn \endSb}^{\infty}
   \hskip -10pt
     \frac{(1-\xi^{4gn})C\lr{\nu}_{4gn}}{(4gn)_{\nu}}
     \frac{{u}^{4gn-\nu}}{(4gn-\nu)!}
  \tag{14.1.5}
  $$
instead of (13.3.1), we have
\vskip 5pt
\fp
\lbox{125}{360}{{\bf Theorem 14.1.6.}\ \ \it
Let  $p\equiv 1$  {\rm mod}  $4g$  be a prime number, 
$a$  be a positive integer.  
Let  $\nu$  be an integer such that  $1\leqq \nu\leqq 4g+2$.  
Assume  $4gn=a+\nu$  and  $p-1\hskip -3pt\not|\hskip 2pt 4gn$.  
Then we have
  $$
  \align
  \sum_{r=0}^a\binom{r}{a}(-A_p)^{a-r}
  \frac{C\lr{\nu}_{4gn+r(p-1)}}{(4gn+r(p-1))_{\nu}}
  &\equiv 0 \mod{p^a}, \\
  \sum_{r=0}^a\binom{r}{a}(-A_p)^{a-r}
  \frac{D\lr{\nu}_{4gn+r(p-1)}}{(4gn+r(p-1))_{\nu}}
  &\equiv 0 \mod{p^a}.
  \endalign
  $$
}

\newpage \fp
{\bf 14.2. On the Hurwitz coefficients of  $\bk{t}(\bk{u}_g)$  and of  $\bk{s}(\bk{u}_g)$.}\ \ 
We consider now the Hurwitz coefficients  $c\lr{h}_m$,  $d\lr{h}_m$  
of powers of  $t=t(u)$  in 8.2.1 and  $s=s(t)$  in 8.3.1, namely, those of
  $$
  \frac{t(u)^h}{h!}=\sum_{m=h}^{\infty}c\lr{h}_m\frac{{u}^m}{m!}, \hskip 30pt
  \frac{s(u)^h}{h!}=\sum_{m=h}^{\infty}d\lr{h}_m\frac{{u}^m}{m!}
  \hskip 10pt\in \bold Z_p\llg u\rrg.
  \tag{14.2.1}
  $$
Here, obviously  $c\lr{h}_h=d\lr{h}_h=1$.  
The results are as follows.
\vskip 3pt
\fp
\lbox{97}{360}{{\bf Lemma 14.2.2.}\ \ \it 
If  $\nu>0$  and  $m\geqq a+1$, then
  $$
  \aligned
  &\sum_{r=0}^a\binom{a}{r}(-A_p)^{a-r}c\lr{h}_{m+r(p-1)}\equiv 0 \mod{p^a},  \\ 
  &\sum_{r=0}^a\binom{a}{r}(-A_p)^{a-r}d\lr{h}_{m+r(p-1)}\equiv 0 \mod{p^a}.
  \endaligned
  $$
}
\fp
{\it Proof.} \ 
By 8.1.2 (2), (3), and 8.2.1, we have
  $$
  \bold Z_p\llg u\rrg=\bold Z_p\llg t\rrg.  
  \tag{14.2.3}
  $$
For  $D=d/du$  and any integer  $h\geqq 1$, we show that
  $$
  (D^p-A_pD)\Big(\frac{t^h}{h!}\Big)\in b + p\bold Z_p\llg u\rrg,
  \hskip 20pt  (b\in \bold Z_p)
  \tag{14.2.4}
  $$
by induction.  
If  $h=1$, this is checked by 13.2.2 and (14.2.3).   
Since
  $$
  \align
  D(D^p-A_pD)\Big(&\jrac{t^{h+1}}{(h+1)!}\Big)
   =(D^p-A_pD)\Big(\jrac{t^h}{h!}\jrac{dt}{du}\Big) \\
  &=\Big(D^p\jrac{t^h}{h!}\Big)\jrac{dt}{du}
   + \sum_{j=1}^{p-1}\binom{p}{j}\Big(D^{p-j}\jrac{t^h}{h!}\Big)\Big(D^j\jrac{dt}{du}\Big)
     + \jrac{t^h}{h!}\Big(D^p\jrac{dt}{du}\Big) \\
  &\hskip 30pt 
    - A_pD\Big(\jrac{t^h}{h!}\Big)\jrac{dt}{du} - A_p\jrac{t^h}{h!}\Big(D\jrac{dt}{du}\Big) \\
  & \in \Big\{(D^p-A_pD)\jrac{u^h}{h!}\Big\}\jrac{dt}{du}
      + \jrac{t^h}{h!}\Big\{(D^p-A_pD)\jrac{dt}{du}\Big\} 
      + p\bold Z_p\llg u\rrg
  \endalign
  $$
by using the hypothesis of induction, (13.2.4), and 13.2.2, 
we see that this belongs to  $p\bold Z_p\llg u\rrg$.  
Hence, we have
  $$
  (D^p-A_pD)\Big(\jrac{t^{h+1}}{(h+1)!}\Big)\in b+p\bold Z_p\llg u\rrg
  \hskip 30pt (b\in \bold Z_p).
  $$
Thus, we have concluded (14.2.4).  
Using (14.2.4) repeatedly, we have, for  $a>0$, that
  $$
  (D^p-A_pD)^a\Big(\frac{t^h}{h!}\Big)\in b + p^a\bold Z_p\llg u\rrg 
  \hskip 20pt  (b\in \bold Z_p). 
  \tag{14.2.5}
  $$
Because of
  $$
  (D^p-A_pD)^a\Big(\frac{t^h}{h!}\Big)
  =\sum_{m=h}^{\infty}\bigg\{\sum_{r=0}^a\binom{a}{r}(-A_p)^{a-r}{c\lr{h}_{m+r(p-1)}}\bigg\}
   D^a\jrac{u^m}{m!},
  $$
we obtain the first formula.  
The second one is proved similarly.  
\qed

\newpage \fp
{\bf 14.3. The Vandiver-Carlitz type congruence relation.}\ \ \ 
In this subsection, we prove the congruence relation of Vandiver-Carlitz type.  
Vandiver gave in \cite{{\bf V}} such the congruence for Bernoulli numbers, 
and Carlitz investigated in \cite{{\bf Ca4}} for Hurwitz numbers.  
\vskip 5pt
\fp
\lbox{110}{360}{{\bf Proposition 14.3.1.}\ \ \it
Let  $a>0$  and  $n=a+1$  are integers.  
For the numbers  $\{C\lr{\nu}_n\}$  and  $\{D\lr{\nu}_n\}$ 
defined in {\rm (14.1.1)} or {\rm (14.1.4)}, we have
  $$
  \align
  \sum_{r=0}^a\binom{a}{r}(-A_p)^{a-r}
  C\lr{\nu}_{n+r(p-1)}
  &\equiv 0 \mod{p^{a-\nu}}, \\
  \sum_{r=0}^a\binom{a}{r}(-A_p)^{a-r}
  D\lr{\nu}_{n+r(p-1)}
  &\equiv 0 \mod{p^{a-\nu}}.
  \endalign
  $$
}
\fp
{\bf Remark 14.3.2.}\ \
These congruence relations need not the condition  $(p-1)\hskip -3pt\not|\hskip 2pt n$    
in contrary to 7.1.1, 7.2.1, 14.1.3, 14.1.6. 
However, since 14.3.1 can be proved without (13.1.10) and 13.1.1, 
it would be not so deep results.
\vskip 3pt
\fp
{\it Proof. } \ 
As in 10.1, we let  $[t^{n+1}/(n+1)]u(t)=f_n$.   
Here we use still the notation  $u=u_g$, $D=d/du$.  
Since
  $$
  \align
  \sum_{m=0}^{\infty}C\lr{1}_{m}\frac{u^m}{m!}
  &=\frac{u}{t(u)} 
   =\frac{1}{t}
   \sum_{m=\nu}^{\infty}f_h\frac{t^{h+1}}{h+1} \\
  &=\sum_{h=0}^{\infty}
        \frac{f_h h!}{h+1} \sum_{m=h}^{\infty}c\lr{h}_m\frac{u^m}{m!} 
   =\sum_{m=0}^{\infty}
      \bigg(\sum_{h=0}^m \frac{f_h h!}{h+1}c\lr{h}_m\bigg)
    \frac{u^m}{m!}
  \endalign
  $$
we have
  $$
  C\lr{1}_m=\sum_{h=0}^m \frac{f_h h!}{h+1}c\lr{h}_m.  
  \tag{14.3.3}
  $$
Therefore, 
  $$
  \align 
  \sum_{r=0}^a\binom{a}{r}(-A_p)^{a-r}C\lr{1}_{m+r(p-1)}
  &=\sum_{r=0}^a\binom{a}{r}(-A_p)^{a-r}
    \sum_{h=0}^{m+r(p-1)}\frac{f_h h!}{h+1}c\lr{h}_{m+r(p-1)}  \\
  &=\sum_{h=0}^{m+r(p-1)}
    \frac{f_h h!}{h+1}
    \sum_{m=0}^a
     \binom{a}{r}(-A_p)^{a-r}c\lr{h}_{m+r(p-1)}.
  \endalign
  $$
If  $m\geqq a+1$, then the inner sum is divisible by  $p^a$  because of (14.2.2). 
If  $p^w||(h+1)$ ($w\geqq 1$), 
then  $p^w-1\leqq h$, and  $\ord_p(p^w-1)!\leqq \ord_ph!$.  
By (1.1.2), we see that
  $$
  \ord_p(p^w-1)!
  = \jrac{(p^w-1)-(p-1)w}{p-1}
  = p^{w-1}+p^{w-2}+\cdots+p+1-w. 
  $$
In our situation, we may assume  $p\leqq 3$.  
Hence, the above is  $\geqq w$  if  $w\geqq 2$.  
The worst case is when  $w=1$, so that
  $$
  w\leqq\ord_p h!+1, \hskip 30pt (w\geqq 1).
  \tag{14.3.4}
  $$
Therefore
  $$
  \sum_{r=0}^a\binom{a}{r}(-A_p)^{a-r}C\lr{1}_{m+r(p-1)}
  \ \in \ p^{a-1}\bold Z_{(p)}\hskip 30pt 
  (m\geqq a+1). 
  \tag{14.3.5}
  $$
This is just the case of  $\nu=1$  of the first formula in 14.3.1. 
The general case is proved similarly as follows.  
Since  $u=\sum_{h=0}^{\infty}f_h\frac{t^{h+1}}{h+1}$, we see
  $$
   u^{\nu}
  =\sum_{h=0}^{\infty}\sum_{\Sb h_1+h_2+\cdots\\ +h_{\nu}=h\endSb}
   \frac{f_{h_1}}{h_1+1}
   \frac{f_{h_2}}{h_2+1}
   \cdots
   \frac{f_{h_{\nu}}}{h_{\nu}+1}
  t^{h+\nu}.
  $$
Hence
  $$
  \align
  \Big(\frac{u}{t}\Big)^{\nu}\hskip -5pt
  &=\hskip -3pt
   \sum_{h=0}^{\infty}
     \sum_{\Sb h_1+h_2+\cdots\\ +h_{\nu}=h\endSb}
     \hskip -3pt
     \frac{f_{h_1}}{h_1+1}
     \frac{f_{h_2}}{h_2+1}
     \cdots
     \frac{f_{h_{\nu}}}{h_{\nu}+1}
   t^h  \\
  &=\hskip -3pt
   \sum_{h=0}^{\infty}
     \sum_{\Sb h_1+h_2+\cdots\\ +h_{\nu}=h\endSb}
     \hskip -3pt
     \frac{f_{h_1}}{h_1+1}
     \frac{f_{h_2}}{h_2+1}
     \cdots
     \frac{f_{h_{\nu}}}{h_{\nu}+1}
    \hskip 2pt h!
   \sum_{m=h}^{\infty}c\lr{h}_m\frac{u^m}{m!}\\
  &=\sum_{m=0}^{\infty}
   \Bigg(\sum_{h=0}^m
    \hskip 5pt h!
    \hskip -10pt
     \sum_{\Sb h_1+h_2+\cdots\\ +h_{\nu}=h\endSb}
      \hskip -3pt
      \frac{f_{h_1}}{h_1+1}
      \frac{f_{h_2}}{h_2+1}
      \cdots
      \frac{f_{h_{\nu}}}{h_{\nu}+1}
      c\lr{h}_m
    \Bigg)
   \frac{u^m}{m!},
  \endalign
  $$
and that
  $$
  \frac{1}{{t}^{\nu}}
  =\sum_{m=0}^{\infty}
   \frac{1}{(m)_{\nu}}
   \Bigg(
    \sum_{h=0}^m
    \hskip 5pt h!
    \hskip -10pt
     \sum_{\Sb h_1+h_2+\cdots\\ +h_{\nu}=h\endSb}
     \hskip -3pt
      \frac{f_{h_1}}{h_1+1}
      \frac{f_{h_2}}{h_2+1}
      \cdots
      \frac{f_{h_{\nu}}}{h_{\nu}+1}
      c\lr{h}_m
    \Bigg)
   \frac{u^{m-\nu}}{(m-\nu)!}.
  $$
This development yields that
  $$
  C\lr{\nu}_m
  =\sum_{h=0}^m
    \hskip 5pt h!
    \hskip -10pt
     \sum_{\Sb h_1+h_2+\cdots\\ +h_{\nu}=h\endSb}
      \hskip -3pt
      \frac{f_{h_1}}{h_1+1}
      \frac{f_{h_2}}{h_2+1}
      \cdots
      \frac{f_{h_{\nu}}}{h_{\nu}+1}
      c\lr{h}_m. 
  \tag{14.3.6}
  $$
and that 
  $$
  \align
  \hskip 30pt
  \sum_{r=0}^a\binom{a}{r}(-A_p)^{a-r}C\lr{\nu}_{n+r(p-1)} 
  &=\sum_{h=0}^m 
   \Bigg(\hskip -6pt
    \hskip 5pt h!
    \hskip -10pt
     \sum_{\hskip 5pt \Sb h_1+h_2+\cdots\\ +h_{\nu}=h\endSb}
      \hskip -4pt
      \frac{f_{h_1}}{h_1+1}
      \frac{f_{h_2}}{h_2+1}
      \cdots
      \frac{f_{h_{\nu}}}{h_{\nu}+1}
   \Bigg) \\
  &\hskip 55pt
   \cdot\bigg\{
   \sum_{r=0}^a\binom{a}{r}(-A_p)^{a-r}c\lr{h}_{n+r(p-1)}
   \bigg\}.
  \tag{14.3.7}
  \endalign
  $$
If we assume  $p^{w_j}||(h_j+1)$, then, as in (14.3.4), we conclude
  $$
  w_1+\cdots+w_{\nu}
  \leqq \ord_p({h_1}!{h_2}!\cdots{h_{\nu}}!)+\nu
  \leqq \ord_p(h!)+\nu
  $$
because  $w_j\leqq \ord_p{h_j}!+1$.    
This shows that the right hand side of (14.3.7) is divisible by  $p^{a-\nu}$.  
Thus, we have proved the first congruence in (14.3.1).  
The second one is proved similarly.
\qed

\newpage \fp\inbox{\bf 15}\ \ {\bf Numerical examples for classical numbers}
\vskip 3pt
\fp
{\bf 15.1. Bernoulli numbers.} \ 
First several values of  $C_{2n}=(-1)^{n-1}2^{2n}B_{2n}$ ($B_{2n}$  is the  $2n$-th Bernoulli number) 
of the curve  $y^2=x-1$  ($g=0$)  are as follows:
$$
\allowdisplaybreaks
\align
 2^2B_2\ &=\frac{1}{3}\cdot 2, \ \hskip 69pt
-2^4B_4  = \frac{1}{3\cdot 5}\cdot 2^3, \\
 2^6B_6\ &=\frac{1}{3\cdot 7}\cdot 2^5, \ \hskip 52pt 
-2^8B_8  =\frac{1}{3\cdot 5}\cdot 2^7, \\
 2^{10}B_{10}&=\frac{1}{3\cdot 11}\cdot 2^9\cdot 5, \ \hskip 27pt
-2^{12}B_{12}= \frac{1}{3\cdot 5\cdot 7\cdot 13}\cdot 2^{11}\cdot 691, \\
 2^{14}B_{14}&=\frac{1}{3}\cdot 2^{15}\cdot 7, \ \hskip 40pt
-2^{16}B_{16}= \frac{1}{3\cdot 5\cdot 17}\cdot 2^{15}\cdot 3617, \\
 2^{18}B_{18}&=\frac{1}{3\cdot 5\cdot 19}\cdot 2^{17}\cdot 43867, \\
-2^{20}B_{20}&=\frac{1}{3\cdot 5\cdot 11}\cdot 2^{19}\cdot 283\cdot 617, \\
 2^{22}B_{22}&=\frac{1}{3\cdot 23}\cdot 2^{21}\cdot 11\cdot 131 \cdot 593, \\
-2^{24}B_{24}&=\frac{1}{3\cdot 5\cdot 7\cdot 13}\cdot 2^{23}\cdot 103 \cdot 2294797, \\
 2^{26}B_{26}&=\frac{1}{3}\cdot 2^{25}\cdot 13 \cdot 657931, \\
-2^{28}B_{28}&=\frac{1}{3\cdot 5\cdot 29}\cdot 2^{27}\cdot 7\cdot 9349 \cdot 362903, \\
 2^{30}B_{30}&=\frac{1}{3\cdot 7\cdot 11\cdot 31}\cdot 2^{29}\cdot 5\cdot 1721 \cdot 1001259881, \\
-2^{32}B_{32}&=\frac{1}{3\cdot 5\cdot 17}\cdot 2^{31}\cdot 37\cdot 683\cdot 305065927, \\
 2^{34}B_{34}&=\frac{1}{3}\cdot 2^{33}\cdot 2^{33}\cdot 17\cdot 151628697551, \\
-2^{36}B_{36}&=\frac{1}{3\cdot 5\cdot 7\cdot 13\cdot 19\cdot 37}
             \cdot 2^{35}\cdot 26315271553053477373, \\
 2^{38}B_{38}&=\frac{1}{3}\cdot 2^{37}\cdot 19\cdot 154210205991661, \\
-2^{40}B_{40}&=\frac{1}{3\cdot 5\cdot 11\cdot 41}
             \cdot 2^{39}\cdot 137616929\cdot 1897170067619, \\ 
 2^{42}B_{42}&=\frac{1}{3\cdot 7\cdot 43}
             \cdot 2^{41}\cdot 1520097643918070802691, \\
-2^{44}B_{44}&=\frac{1}{3\cdot 5\cdot 23}
             \cdot 2^{43}\cdot 11\cdot 59\cdot 8089\cdot 2947939\cdot 1798482437, \\
 2^{46}B_{46}&=\frac{1}{3\cdot 47}
             \cdot 2^{45}\cdot 23\cdot 383799511\cdot 67568238839737, \\
-2^{48}B_{48}&=\frac{1}{3\cdot 5\cdot 7\cdot 13\cdot 17}
             \cdot 2^{47}\cdot 653\cdot 56039\cdot 153289748932447906241, \\
 2^{50}B_{50}&=\frac{1}{3\cdot 11}
             \cdot 2^{49}\cdot 5^2\cdot 417202699\cdot 47464429777438199.
\endalign 
$$

\newpage \input c3.o
\newpage \fp
{\bf 15.3. Hurwitz numbers for the curve  $\bk{y}^2=\bk{x}^3-\bk{x}$.} \ 
First several values of  $C_{4n}=2^{4n}E_{4n}$  
($E_{4n}$  is the $4n$-th original Hurwitz number) for the curve  $y^2=x^3-x$  
($g=1$) are as follows:
$$
\allowdisplaybreaks
\align
2^4E_4&=\frac{1}{5}\cdot 2^3\cdot 3, \\
2^8E_8&=\frac{1}{5}\cdot 2^7\cdot 3, \\
2^{12}E_{12}&=\frac{1}{5\cdot 13}\cdot 2^{11}\cdot 3^4\cdot 7, \\
2^{16}E_{16}&=\frac{1}{5\cdot 17}\cdot 2^{15}\cdot 3^4\cdot 7^2\cdot 11, \\
2^{20}E_{20}&=\frac{1}{5}\cdot 2^{19}\cdot 3^6\cdot 7^2\cdot 11, \\
2^{24}E_{24}&=\frac{1}{5\cdot 13}\cdot 2^{23}\cdot 3^7\cdot 7^3\cdot 11^2\cdot 19, \\
2^{28}E_{28}&=\frac{1}{5\cdot 29}\cdot 2^{27}\cdot 3^9\cdot 7^4\cdot 11^2\cdot 19\cdot 23, \\
2^{32}E_{32}&=\frac{1}{5\cdot 17}\cdot 2^{31}\cdot 3^{10}\cdot 7^4\cdot 11^2\cdot 19\cdot 23\cdot 223, \\
2^{36}E_{36}&=\frac{1}{5\cdot 13\cdot 37}
     \cdot 2^{35}\cdot 3^{14}\cdot 7^5\cdot 11^3\cdot 19\cdot 23\cdot 31\cdot 61, \\
2^{40}E_{40}&=\frac{1}{5\cdot 41}
     \cdot 2^{39}\cdot 3^{13}\cdot 7^5\cdot 11^3\cdot 19^2\cdot 23\cdot 31\cdot 2381, \\
2^{44}E_{44}&=\frac{1}{5}
     \cdot 2^{43}\cdot 3^{15}\cdot 7^6\cdot 11^4\cdot 19^4\cdot 23\cdot 31, \\
2^{48}E_{48}&=\frac{1}{5\cdot 13\cdot 17}
     \cdot 2^{47}\cdot 3^{16}\cdot 7^5\cdot 11^4\cdot 19^2\cdot 23^2\cdot 31\cdot 43\cdot 1162253, \\
2^{52}E_{52}&=\frac{1}{5\cdot 53}
     \cdot 2^{51}\cdot 3^{18}\cdot 7^7\cdot 11^4\cdot 13\cdot 19^2\cdot 23^2\cdot 31\cdot 43
     \cdot 47\cdot 8887, \\
2^{56}E_{56}&=\frac{1}{5\cdot 29}
     \cdot 2^{55}\cdot 3^{19}\cdot 7^8\cdot 11^5\cdot 19^2\cdot 23^2\cdot 31\cdot 43
     \cdot 47\cdot 61\cdot 52289, \\
2^{60}E_{60}&=\frac{1}{5\cdot 13\cdot 61}
     \cdot 2^{59}\cdot 3^{22}\cdot 7^8\cdot 11^5\cdot 19^3\cdot 23^2 \cdot 31
     \cdot 43\cdot 47\cdot 2630966033, \\
2^{64}E_{64}&=\frac{1}{5\cdot 17}
     \cdot 2^{63}\cdot 3^{22}\cdot 7^9\cdot 11^5\cdot 19^3\cdot 23^2\cdot 31^2
     \cdot 43\cdot 47\cdot 59\cdot 109\cdot 814903,\\
2^{68}E_{68}&=\frac{1}{5}
     \cdot 2^{67}\cdot 3^{24}\cdot 7^9\cdot 11^6\cdot 17\cdot 19\cdot 23^2
     \cdot 31^2\cdot 43\cdot 47\cdot 59\cdot 80232721, \\
2^{72}E_{72}&=\frac{1}{5\cdot 13\cdot 37\cdot 73}
     \cdot 2^{71}\cdot 3^{25}\cdot 7^{10}\cdot 11^6\cdot 19^3\cdot 23^3\cdot 31^2
     \cdot 43\cdot 47\cdot 59\cdot 67\\ & \cdot 48316510111193, \\
2^{76}E_{76}&=\frac{1}{5}
     \cdot 2^{75}\cdot 3^{27}\cdot 7^{10}\cdot 11^6\cdot 19^4\cdot 23^3\cdot 31^2
     \cdot 43\cdot 47\cdot 59\cdot 67\cdot 71\cdot 3469\cdot 1330177.
\endalign
$$

\newpage \fp\inbox{\bf 16}\ \ {\bf Numerical examples for new numbers}
\vskip 3pt
\fp
{\bf 16.1.  $\bk{x}(\bk{u})$  of the curve  $\bk{y}^2=\bk{x}^5-\bold 1$.} \ \ \ 
First several values of  $C_{10n}$  for the curve  $y^2=x^5-1$  ($g=2$) are as follows:
$$
\allowdisplaybreaks
\align
C_{10}&=\frac{1}{11}\cdot 2^8 \cdot 3^2 \cdot 5^2 \cdot 7, \\
C_{20}&=\frac{-1}{11}\cdot
 2^{18}\cdot
 3^9   \cdot
 5^4   \cdot
 7     \cdot
 13    \cdot
 17,   \\
C_{30}&=\frac{1}{11 \cdot 31} \cdot
 2^{28} \cdot
 3^{14} \cdot
 5^7    \cdot
 7^3    \cdot
 13^2   \cdot
 17     \cdot
 19     \cdot
 23^2,  \\
C_{40}&=\frac{-1}{11 \cdot 41}  \cdot
 2^{38} \cdot
 3^{17} \cdot
 5^9    \cdot
 7^3    \cdot
 13^2   \cdot
 17^2   \cdot
 19^2   \cdot
 23     \cdot
 29     \cdot
 37     \cdot
 31991,  \\
C_{50}&=\frac{1}{11}  \cdot
 2^{47} \cdot
 3^{23} \cdot
 5^{12} \cdot
 7^6    \cdot
 13^3   \cdot
 17     \cdot
 19^2   \cdot
 23^2   \cdot
 29     \cdot
 37     \cdot
 43     \cdot
 47     \cdot
 4999,  \\
C_{60}&=\frac{-1}{11 \cdot 31 \cdot 61}  \cdot
 2^{59} \cdot
 3^{28} \cdot
 5^{15} \cdot
 7^6    \cdot
 13^4   \cdot
 17^2   \cdot
 19^3   \cdot
 23^2   \cdot
 29^2   \cdot
 37     \cdot
 43 \\ &\cdot
 47     \cdot
 53     \cdot
 351453077,   \\
C_{70}&=\frac{1}{11 \cdot 71}  \cdot
 2^{72} \cdot
 3^{31} \cdot
 5^{16} \cdot
 7^9    \cdot
 13^5   \cdot
 17^3   \cdot
 19^3   \cdot
 23^2   \cdot
 29^2   \cdot
 37     \cdot
 43 \\ &\cdot
 47     \cdot
 53     \cdot
 59     \cdot
 67     \cdot
 6740734411, \\
C_{80}&=\frac{-1}{11 \cdot 41} \cdot
 2^{78} \cdot
 3^{34} \cdot
 5^{19} \cdot
 7^8    \cdot
 13^6   \cdot
 17^3   \cdot
 19^4   \cdot
 23^3   \cdot
 29^2   \cdot
 37^2   \cdot
 43 \\ &\cdot
 47     \cdot
 53     \cdot
 59     \cdot
 67     \cdot
 73     \cdot
 109    \cdot 
 460903 \cdot
 121384433,  \\
C_{90}&=\frac{1}{11 \cdot 31} \cdot
 2^{87} \cdot
 3^{42} \cdot
 5^{21} \cdot
 7^{10} \cdot
 13^6   \cdot
 17^4   \cdot
 19^4   \cdot
 23^3   \cdot
 29^3   \cdot
 37^2   \cdot
 43^2 \\ &\cdot
 47     \cdot
 53     \cdot
 59     \cdot
 67^2   \cdot
 73     \cdot
 79     \cdot
 83     \cdot
 131    \cdot
 881    \cdot
 2799606697,   \\
C_{100}&=\frac{-1}{11 \cdot 101} \cdot
 2^{97} \cdot
 3^{47} \cdot
 5^{24} \cdot
 7^{11} \cdot
 13^7   \cdot
 17^3   \cdot
 19^5   \cdot
 23^4   \cdot
 29^3   \cdot
 37^2   \cdot
 43^2 \\ &\cdot
 47^2   \cdot
 53     \cdot
 59     \cdot
 67     \cdot
 73     \cdot
 79     \cdot
 83     \cdot
 89     \cdot
 97     \cdot
 10343  \cdot
 1938718187373563,   \\
C_{110}&=\frac{1}{11}  \cdot
 2^{107} \cdot
 3^{51} \cdot
 5^{27} \cdot
 7^{13} \cdot
 13^8   \cdot
 17^4   \cdot
 19^5   \cdot
 23^4   \cdot
 29^3   \cdot
 37     \cdot
 43^2 \\ &\cdot
 47^2   \cdot
 53^2   \cdot
 59     \cdot
 67     \cdot
 73     \cdot
 79     \cdot
 83     \cdot
 89     \cdot
 97     \cdot
 103    \cdot
 107 \\ &\cdot
 3019729 \cdot
 865724129494813, \\
C_{120}&=\frac{-1}{11 \cdot 31 \cdot 41 \cdot 61}  \cdot
 2^{119} \cdot
 3^{56} \cdot
 5^{29} \cdot
 7^{13} \cdot
 13^9   \cdot
 17^5   \cdot
 19^6   \cdot
 23^5   \cdot
 29^4   \cdot
 37^2   \cdot
 43^2 \\ &\cdot
 47^2   \cdot
 53^2   \cdot
 59^2   \cdot
 67     \cdot
 73     \cdot
 79     \cdot
 83     \cdot
 89     \cdot
 97     \cdot
 103    \cdot
 107    \cdot
 109    \cdot
 113 \\ &\cdot
 863833294249 \cdot
 7389430581319, \\
C_{130}&=\frac{1}{11 \cdot 131}  \cdot
 2^{128}\cdot
 3^{61} \cdot
 5^{32} \cdot
 7^{15} \cdot
 13^{11}\cdot
 17^5   \cdot
 19^6   \cdot
 23^5   \cdot
 29^4   \cdot
 37^2   \cdot
 43^2 \\ &\cdot
 47^2   \cdot
 53^2   \cdot
 59^2   \cdot
 67     \cdot
 73     \cdot
 79     \cdot
 83     \cdot
 89     \cdot
 97     \cdot
 103    \cdot
 107    \cdot
 109    \cdot
 113    \cdot
 127 \\& \cdot
 5303   \cdot
 97785319  \cdot
 175363749323953511, \\
C_{140}&=\frac{-1}{11 \cdot 71} \cdot
 2^{139}\cdot
 3^{65} \cdot
 5^{34} \cdot
 7^{15} \cdot
 13^{10}\cdot
 17^6   \cdot
 19^7   \cdot
 23^6   \cdot
 29^4   \cdot
 37^2   \cdot
 43^3 \\ &\cdot
 47     \cdot
 53^2   \cdot
 59^2   \cdot
 67^2   \cdot
 73     \cdot
 79     \cdot
 83     \cdot
 89     \cdot
 97     \cdot
 103    \cdot
 107    \cdot
 109    \cdot
 113    \cdot
 127 \\ &\cdot
 137    \cdot
 3191   \cdot
 79927801 \cdot
 2927519326077590415331021, \\
C_{150}&=\frac{1}{11 \cdot 31 \cdot 151} \cdot
 2^{150}\cdot
 3^{70} \cdot
 5^{37} \cdot
 7^{17} \cdot
 13^{12}\cdot
 17^5   \cdot
 19^7   \cdot
 23^6   \cdot
 29^5   \cdot
 37^3   \cdot
 43^3 \\ &\cdot
 47^2   \cdot
 53^2   \cdot
 59^2   \cdot
 67^2   \cdot
 73^2   \cdot
 79     \cdot
 83     \cdot
 89     \cdot
 97     \cdot
 103    \cdot
 107    \cdot
 109    \cdot
 113    \cdot
 127 \\ &\cdot
 137    \cdot
 139    \cdot
 50951  \cdot
 450127 \cdot
 1464426640811 \cdot
 58871719018640089, \\
C_{160}&=\frac{-1}{11 \cdot 41} \cdot
 2^{158} \cdot
 3^{72}  \cdot
 5^{40}  \cdot
 7^{17}  \cdot
 13^{12} \cdot
 17^6    \cdot
 19^8    \cdot
 23^6    \cdot
 29^5    \cdot
 37^3    \cdot
 43^3 \\ &\cdot
 47^2    \cdot
 53^2    \cdot
 59^2    \cdot
 67^2    \cdot
 73^2    \cdot
 79^2    \cdot
 83      \cdot
 89      \cdot
 97      \cdot
 103     \cdot
 107     \cdot
 109     \cdot
 113     \cdot
 127     \cdot
 137 \\ &\cdot
 139     \cdot
 149     \cdot
 157     \cdot
 5473709 \cdot
 22543502622365730931551293201565706511, \\
C_{170}&=\frac{1}{11} \cdot
 2^{167} \cdot
 3^{78}  \cdot
 5^{42}  \cdot
 7^{19}  \cdot
 13^{12} \cdot
 17^8    \cdot
 19^8    \cdot
 23^7    \cdot
 29^5    \cdot
 37^3    \cdot
 43^3 \\ &\cdot
 47^2    \cdot
 53^3    \cdot
 59^2    \cdot
 67^2    \cdot
 73^2    \cdot
 79^2    \cdot
 83^2    \cdot
 89      \cdot
 97      \cdot
 103     \cdot
 107     \cdot
 109     \cdot
 113     \cdot
 127 \\ &\cdot
 137     \cdot
 139     \cdot
 149     \cdot
 157     \cdot
 163     \cdot
 167 \\ &\cdot
 587     \cdot
 22573   \cdot
 18793   \cdot
 246289  \cdot
 311203545376580358674935387, \\
C_{180}&=\frac{-1}{11 \cdot 31 \cdot 61 \cdot 181} \cdot
 2^{177} \cdot
 3^{87}  \cdot
 5^{45}  \cdot
 7^{19}  \cdot
 13^{15} \cdot
 17^7    \cdot
 19^9    \cdot
 23^7    \cdot
 29^6    \cdot
 37^3    \cdot
 43^4 \\ &\cdot
 47^2    \cdot
 53^3    \cdot
 59^3    \cdot
 67^2    \cdot
 73^2    \cdot
 79^2    \cdot
 83^2    \cdot
 89^2    \cdot
 97      \cdot
 103     \cdot
 107     \cdot
 109     \cdot
 113     \cdot
 127 \\ &\cdot
 137     \cdot
 139     \cdot
 149     \cdot
 157     \cdot
 163     \cdot
 167     \cdot
 173 \\ &\cdot
 239     \cdot
 1471    \cdot
 1579    \cdot
 7030999221688667065861742323016843138707.
\endalign
$$

\newpage 

\fp
Each  $C_{10n}$  is writen by a certain integer  $G_{10n}$ 
(von Staudt-Clausen type theorem) as follows:
$$
\alignat {10}
C_{ 10}&=\frac{6     }{11}&&                &&                &&                &&                &&               &&               &&                &&                   &&+G_{10}, \\
C_{ 20}&=\frac{6^  2 }{11}&&                &&                &&                &&                &&               &&               &&                &&                   &&+G_{20}, \\
C_{ 30}&=\frac{6^  3 }{11}&&+\frac{10  }{31}&&                &&                &&                &&               &&               &&                &&                   &&+G_{30}, \\
C_{ 40}&=\frac{6^  4 }{11}&&                &&+\frac{ 7  }{41}&&                &&                &&               &&               &&                &&                   &&+G_{40}, \\
C_{ 50}&=\frac{6^{ 5}}{11}&&                &&                &&                &&                &&               &&               &&                &&                   &&+G_{50}, \\
C_{ 60}&=\frac{6^{ 6}}{11}&&+\frac{10^2}{31}&&                &&+\frac{ 1  }{61}&&                &&               &&               &&                &&                   &&+G_{60}, \\
C_{ 70}&=\frac{6^{ 7}}{11}&&                &&                &&                &&+\frac{32  }{71}&&               &&               &&                &&                   &&+G_{70}, \\
C_{ 80}&=\frac{6^{ 8}}{11}&&                &&+\frac{ 7^2}{41}&&                &&                &&               &&               &&                &&                   &&+G_{80}, \\
C_{ 90}&=\frac{6^{ 9}}{11}&&+\frac{10^3}{31}&&                &&                &&                &&               &&               &&                &&                   &&+G_{90}, \\
C_{100}&=\frac{6^{10}}{11}&&                &&                &&                &&                &&+\frac{46}{101}&&               &&                &&                   &&+G_{100}, \\
C_{110}&=\frac{6^{11}}{11}&&                &&                &&                &&                &&               &&               &&                &&                   &&+G_{110}, \\
C_{120}&=\frac{6^{12}}{11}&&+\frac{10^4}{31}&&+\frac{ 7^3}{41}&&+\frac{ 1  }{61}&&                &&               &&               &&                &&                   &&+G_{120}, \\
C_{130}&=\frac{6^{13}}{11}&&                &&                &&                &&                &&               &&+\frac{64}{131}&&                &&                   &&+G_{130}, \\
C_{140}&=\frac{6^{14}}{11}&&                &&                &&                &&+\frac{32^2}{71}&&               &&               &&                &&                   &&+G_{140}, \\
C_{150}&=\frac{6^{15}}{11}&&+\frac{10^5}{31}&&                &&                &&                &&               &&               && +\frac{52}{151}&&                   &&+G_{150}, \\
C_{160}&=\frac{6^{16}}{11}&&                &&+\frac{ 7^4}{41}&&                &&                &&               &&               &&                &&                   &&+G_{160}, \\
C_{170}&=\frac{6^{17}}{11}&&                &&                &&                &&                &&               &&               &&                &&                   &&+G_{170}, \\
C_{180}&=\frac{6^{18}}{11}&&+\frac{10^6}{31}&&                &&+\frac{ 1  }{61}&&                &&               &&               &&                &&+\frac{ 37}{181}\  &&+G_{180}.
\endalignat 
$$

\newpage \input c5y.o
\newpage \fp
{\bf 16.3. $\bk{x}(\bk{u})$  of the curve  $\bk{y}^2=\bk{x}^5-\bk{x}$.} \ \ \ 
First several values of  $C_{8n}$  for the curve  $y^2=x^5-x$  ($g=2$) are as follows:
$$
\allowdisplaybreaks
\align
C_{8}&=
  2^7 \cdot
  5, \\
C_{16}&=\frac{-1}{17}\cdot
  2^{15}  \cdot
  3^2     \cdot
  5^2     \cdot
  7^2     \cdot
 11       \cdot
 13,\\
C_{24}&=
  2^{22}  \cdot
  3^4     \cdot
  5^3     \cdot
  7^3     \cdot
 11^3     \cdot
 13       \cdot
 19, \\
C_{32}&=\frac{-1}{17}\cdot
  2^{31}  \cdot
  3^6     \cdot
  5^6     \cdot
  7^4     \cdot
 11       \cdot
 13^2     \cdot
 19       \cdot
 23       \cdot
 29       \cdot
1741,\\
C_{40}&=\frac{1}{41}\cdot
  2^{40} \cdot
  3^8   \cdot
  5^8   \cdot
  7^5   \cdot
 11^2   \cdot
 13^2   \cdot
 19^2   \cdot
 23     \cdot
 29     \cdot
 31     \cdot
 37     \cdot
 5693,\\
C_{48}&=\frac{-1}{17}\cdot
  2^{46}\cdot
  3^{10}\cdot
  5^8   \cdot,
  7^5   \cdot
 11^3   \cdot
 13^3   \cdot
 19^2   \cdot
 23^2   \cdot
 29     \cdot
 31     \cdot
 37     \cdot
 43     \cdot
 41957857,\\
C_{56}&=
  2^{54}\cdot
  3^{12}\cdot
  5^{11}\cdot
  7^8   \cdot
 11^5   \cdot
 13^4   \cdot
 19     \cdot
 23^2   \cdot
 29     \cdot
 31     \cdot
 37     \cdot
 43     \cdot
 47     \cdot
 53     \cdot
 715991,\\
C_{64}&=\frac{-1}{17}\cdot
  2^{63}\cdot
  3^{14}\cdot
  5^{12}\cdot
  7^9   \cdot
 11^3   \cdot
 13^4   \cdot
 19^2   \cdot
 23^2   \cdot
 29^2   \cdot
 31^2   \cdot
 37     \cdot
 43     \cdot
 47 \\ &\cdot
 53     \cdot
 59     \cdot
 61     \cdot
 89     \cdot
 32591401,\\
C_{72}&=\frac{1}{73}\cdot
  2^{72}\cdot
  3^{16}\cdot
  5^{13}\cdot
  7^{10}\cdot
 11^4   \cdot
 13^5   \cdot
 19^2   \cdot
 23^3   \cdot
 29^2   \cdot
 31^2   \cdot
 37     \cdot
 43     \cdot
 47 \\ &\cdot
 53     \cdot
 59     \cdot
 61     \cdot
 67     \cdot
 19346595547931, \\
C_{80}&=\frac{-1}{17\cdot 41}\cdot
  2^{80}\cdot
  3^{18}\cdot
  5^{17}\cdot
  7^{11}\cdot
 11^5   \cdot
 13^6   \cdot
 19^3   \cdot
 23^3   \cdot
 29^2   \cdot
 31^2   \cdot
 37^2   \cdot
 43     \cdot
 47 \\ &\cdot
 53     \cdot
 59     \cdot
 61^2   \cdot
 67     \cdot
 71     \cdot
 5826608412403.
\endalign
$$

\newpage 

\fp
Each  $C_{8n}$  is written by a certain integer  $G_{8n}$ 
(von Staudt-Clausen type theorem) as follows:
$$
\alignat {11}
C_{  8}&                   && &&               && &&              &&                 && &&               &&                 &&                   && \text{\lq\lq an integer"}, \\
C_{ 16}&=\frac{11     }{17}&& &&               && &&              &&                 && &&               &&                 &&                   &&+G_{16}, \\
C_{ 24}&                   && &&               && &&              &&                 && &&               &&                 &&                   && \text{\lq\lq an integer"}, \\
C_{ 32}&=\frac{11^  2 }{17}&& &&               && &&              &&                 && &&               &&                 &&                   &&+G_{32}, \\
C_{ 40}&=                  && &&\frac{35  }{41}&& &&              &&                 && &&               &&                 &&                   &&+G_{40}, \\
C_{ 48}&=\frac{11^  3 }{17}&& &&               && &&              &&                 && &&               &&                 &&                   &&+G_{48}, \\
C_{ 56}&                   && &&               && &&              &&                 && &&               &&                 &&                   && \text{\lq\lq an integer"}, \\
C_{ 64}&=\frac{11^  4 }{17}&& &&               && &&              &&                 && &&               &&                 &&                   &&+G_{64}, \\
C_{ 72}&=                  && &&               && &&\frac{2  }{73}&&                 && &&               &&                 &&                   &&+G_{72}, \\
C_{ 80}&=\frac{11^  6 }{17}&&+&&\frac{35^2}{41}&& &&              &&                 && &&               &&                 &&                   &&+G_{80}, \\
C_{ 88}&=                  && &&               && &&              &&\ \frac{18  }{89}&& &&               &&                 &&                   &&+G_{88}, \\
C_{ 96}&=\frac{11^  7 }{17}&& &&               && &&              &&                 &&+&&\frac{10  }{97}&&                 &&                   &&+G_{96}, \\
C_{104}&                   && &&               && &&              &&                 && &&               &&                 &&                   && \text{\lq\lq an integer"}, \\
C_{112}&=\frac{11^  8 }{17}&& &&               && &&              &&                 && &&               &&+\frac{18  }{113}&&                   &&+G_{112}, \\
C_{120}&=                  && &&\frac{35^3}{41}&& &&              &&                 && &&               &&                 &&                   &&+G_{120}, \\
C_{128}&=\frac{11^  9 }{17}&& &&               && &&              &&                 && &&               &&                 &&                   &&+G_{128}, \\
C_{136}&=                  && &&               && &&              &&                 && &&               &&                 &&\ \frac{131}{137}\ &&+G_{136}, \\
C_{144}&=\frac{11^{10}}{17}&& &&               &&+&&\frac{2^2}{73}&&                 && &&               &&                 &&                   &&+G_{144}.
\endalignat 
$$

\newpage \fp
{\bf 16.4.  $\bk{x}(\bk{u})$  of the curve  $\bk{y}^2=\bk{x}^7-\bold 1$.} \ \ \ 
First several values of  $C_{14n}$  for the curve  $y^2=x^7-1$  ($g=3$) are as follows:
$$
\allowdisplaybreaks
\align
C_{14}&=
  2^{11} \cdot
  3^4    \cdot
  5      \cdot
  7^2    \cdot
 11,    \\
C_{28}&=\frac{-1}{29}\cdot
  2^{25} \cdot
  3^{11} \cdot
  5^4    \cdot
  7^4    \cdot
  11^2   \cdot
  13^2   \cdot
  17     \cdot
  19     \cdot
  23,   \\
C_{42}&=\frac{1}{43}\cdot
  2^{39} \cdot
  3^{16} \cdot
  5^6    \cdot
  7^6    \cdot
 11^3    \cdot
 13^3    \cdot
 17^2    \cdot
 19^2    \cdot
 23      \cdot
 31^2    \cdot
 37      \cdot
 61,\\
C_{56}&=\frac{-1}{29}\cdot
  2^{53} \cdot
  3^{22} \cdot
  5^9    \cdot
  7^9    \cdot
 11^4    \cdot
 13^4    \cdot
 17^3    \cdot
 19      \cdot
 23      \cdot
 31      \cdot
 37 \\  &\cdot
 41      \cdot
 47      \cdot
 53      \cdot
179      \cdot
2273,\\
C_{70}&=\frac{1}{71}\cdot
  2^{67} \cdot
  3^{28} \cdot
  5^{11} \cdot
  7^{11} \cdot
 11^6    \cdot
 13^5    \cdot
 17^4    \cdot
 19^2    \cdot
 23^2    \cdot
 31^2    \cdot
 37  \\ &\cdot
 41      \cdot
 47      \cdot
 53      \cdot
 59      \cdot
 61      \cdot
 67      \cdot
 1093    \cdot
 31513,\\
C_{84}&=\frac{-1}{29 \cdot 43}\cdot
 2^{81}  \cdot,
 3^{35}  \cdot
 5^{13}  \cdot
 7^{13}  \cdot
 11^7    \cdot
 13^6    \cdot
 17^3    \cdot
 19^3    \cdot
 23^3    \cdot
 31^2    \cdot
 37^2 \\ &\cdot
 41^2    \cdot
 47      \cdot
 53      \cdot
 59      \cdot
 61      \cdot
 67      \cdot
 73      \cdot
 79      \cdot
 14953011635323,  \\
C_{98}&=\frac{1}{29}\cdot
2^{95}  \cdot
3^{40}  \cdot
5^{16}  \cdot
7^{16}  \cdot
11^8  \cdot
13^7  \cdot
17^4  \cdot
19^4  \cdot
23^4  \cdot
31^3  \cdot
37^2 \\ &\cdot
41^2  \cdot
47^2  \cdot
53  \cdot
59  \cdot
61  \cdot
67  \cdot
73  \cdot
79  \cdot
83  \cdot
89  \cdot
5413  \cdot
15937  \cdot
27361,  \\
C_{112}&=\frac{-1}{29\cdot 113}\cdot 
2^{109}  \cdot
3^{47}  \cdot
5^{18}  \cdot
7^{19}  \cdot
11^{11}  \cdot
13^8  \cdot
17^5  \cdot
19^3  \cdot
23^4  \cdot
31^3  \cdot
37^2 \\ &\cdot
41^2  \cdot
47^2  \cdot
53^2  \cdot
59  \cdot
61  \cdot
67  \cdot
73  \cdot
79  \cdot
83  \cdot
89  \cdot
97  \cdot
101  \cdot
103  \cdot
107 \\ &\cdot
109  \cdot
7752349172767981, \\
C_{126}&=\frac{1}{43\cdot 127}\cdot
2^{123}  \cdot
3^{52}  \cdot
5^{22}  \cdot
7^{20}  \cdot
11^{13}  \cdot
13^9  \cdot
17^6  \cdot
19^4  \cdot
23^5  \cdot
31^4  \cdot
37^3 \\ &\cdot
41^3  \cdot
47^2  \cdot
53^2  \cdot
59^2  \cdot
61^2  \cdot
67  \cdot
73  \cdot
79  \cdot
83  \cdot
89  \cdot
97  \cdot
101  \cdot
103  \cdot
107 \\ &\cdot
109  \cdot
47713  \cdot
139246440988973,  \\
C_{140}&=\frac{-1}{29\cdot 71}\cdot
2^{137}  \cdot
3^{57}  \cdot
5^{24}  \cdot
7^{22}  \cdot
11^{13}  \cdot
13^{10}  \cdot
17^7  \cdot
19^5  \cdot
23^6  \cdot
31^4  \cdot
37^3 \\ &\cdot
41^3  \cdot
47  \cdot
53  \cdot
59  \cdot
61  \cdot
67  \cdot
73  \cdot
79  \cdot
83  \cdot
89  \cdot
97  \cdot
101  \cdot
103  \cdot
107 \\ &\cdot
109  \cdot
131  \cdot
137  \cdot
215936153785645363436712289,  \\
C_{154}&=
2^{151}  \cdot
3^{64}  \cdot
5^{26}  \cdot
7^{25}  \cdot
11^{16}  \cdot
13^{11}  \cdot
17^8  \cdot
19^6  \cdot
23^6  \cdot
31^3  \cdot
37^4 \\ &\cdot
41^3  \cdot
47^2  \cdot
53^2  \cdot
59^2  \cdot
61^2  \cdot
67^2  \cdot
73^2  \cdot
79  \cdot
83  \cdot
89  \cdot
97  \cdot
101  \cdot
103  \cdot
107 \\ &\cdot
109  \cdot
131  \cdot
137  \cdot
139  \cdot
149  \cdot
151  \cdot
141376783296257984516233421,  \\
C_{168}&=\frac{-1}{29\cdot 43}\cdot
2^{165}  \cdot
3^{70}  \cdot
5^{29}  \cdot
7^{27}  \cdot
11^{16}  \cdot
13^{11}  \cdot
17^7  \cdot
19^5  \cdot
23^7  \cdot
31^4  \cdot
37^4 \\ &\cdot
41^4  \cdot
47^2  \cdot
53^3  \cdot
59^2  \cdot
61^2  \cdot
67^2  \cdot
73^2  \cdot
79^2  \cdot
83^2  \cdot
89  \cdot
97  \cdot
101  \cdot
103  \cdot
107 \\ &\cdot
109  \cdot
131  \cdot
137  \cdot
139  \cdot
149  \cdot
151  \cdot
157^2  \cdot
163  \cdot
811 \\ &\cdot
417793  \cdot
1745978749  \cdot
41834306314956317,  \\
C_{182}&=
2^{181}  \cdot
3^{76}  \cdot
5^{31}  \cdot
7^{29}  \cdot
11^{17}  \cdot
13^{14}  \cdot
17^8  \cdot
19^6  \cdot
23^7  \cdot
31^4  \cdot
37^4 \\ &\cdot
41^4  \cdot
47^2  \cdot
53^3  \cdot
59^3  \cdot
61  \cdot
67^2  \cdot
73^2  \cdot
79^2  \cdot
83^2  \cdot
89^2  \cdot
97  \cdot
101  \cdot
103  \cdot
107 \\ &\cdot
109  \cdot
131  \cdot
137  \cdot
139  \cdot
149  \cdot
151  \cdot
157  \cdot
163  \cdot
167  \cdot
173  \cdot
179 \\ &\cdot
2538509  \cdot
7602187  \cdot
167935399  \cdot
79956683481346979,  \\
C_{196}&=\frac{-1}{29\cdot 197}\cdot
2^{201}  \cdot
3^{82}  \cdot
5^{34}  \cdot
7^{32}  \cdot
11^{19}  \cdot
13^{15}  \cdot
17^9  \cdot
19^7  \cdot
23^8  \cdot
31^6  \cdot
37^5 \\ &\cdot
41^4  \cdot
47^3  \cdot
53^3  \cdot
59^3  \cdot
61^2  \cdot
67^2  \cdot
73^2  \cdot
79^2  \cdot
83^2  \cdot
89^2  \cdot
97^2  \cdot
101  \cdot
103  \cdot
107 \\ &\cdot
109  \cdot
131  \cdot
137  \cdot
139  \cdot
149  \cdot
151  \cdot
157  \cdot
163  \cdot
167  \cdot
173  \cdot
179  \cdot
181  \cdot
191  \cdot
193 \\ &\cdot
727297  \cdot
64778760224034269867033351318033,  \\
C_{210}&=\frac{1}{43\cdot 71\cdot 211}\cdot
2^{207}  \cdot
3^{88}  \cdot
5^{36}  \cdot
7^{35}  \cdot
11^{19}  \cdot
13^{16}  \cdot
17^{10}  \cdot
19^8  \cdot
23^9  \cdot
31^5  \cdot
37^5 \\ &\cdot
41^5  \cdot
47^3  \cdot
53^3  \cdot
59^3  \cdot
61^2  \cdot
67^3  \cdot
73^2  \cdot
79^2  \cdot
83^2  \cdot
89^2  \cdot
97^2  \cdot
101^2  \cdot
103^2  \cdot
107 \\ &\cdot
109  \cdot
131  \cdot
137  \cdot
139  \cdot
149  \cdot
151  \cdot
157  \cdot
163  \cdot
167  \cdot
173  \cdot
179  \cdot
181  \cdot
191  \cdot
193 \\ &\cdot
199  \cdot
537631123386707971  \cdot
187949832897792288197625137771,  \\
C_{224}&=\frac{-1}{29\cdot 113}\cdot
2^{211}  \cdot
3^{93}  \cdot
5^{38}  \cdot
7^{37}  \cdot
11^{21}  \cdot
13^{17}  \cdot
17^{11}  \cdot
19^7  \cdot
23^9  \cdot
31^6  \cdot
37^6 \\ &\cdot
41^5  \cdot
47^3  \cdot
53^4  \cdot
59^3  \cdot
61^3  \cdot
67^3  \cdot
73^3  \cdot
79^2  \cdot
83^2  \cdot
89^2  \cdot
97^2  \cdot
101^2  \cdot
103^2  \cdot
107^2 \\ &\cdot
109^2  \cdot
131  \cdot
137  \cdot
139  \cdot
149  \cdot
151  \cdot
157  \cdot
163  \cdot
167  \cdot
173  \cdot
179  \cdot
181  \cdot
191  \cdot
193 \\ &\cdot
199  \cdot
2143  \cdot
980275197267685204162754635787499610695199551689,  \\
C_{238}&=\frac{1}{239}\cdot
2^{235}  \cdot
3^{102}  \cdot
5^{42}  \cdot
7^{38}  \cdot
11^{22}  \cdot
13^{18}  \cdot
17^{13}  \cdot
19^8  \cdot
23^{10}  \cdot
31^6  \cdot
37^6 \\ &\cdot
41^5  \cdot
47^4  \cdot
53^4  \cdot
59^4  \cdot
61^2  \cdot
67^4  \cdot
73^3  \cdot
79^2  \cdot
83^2  \cdot
89^2  \cdot
97^2  \cdot
101^2  \cdot
103^2  \cdot
107^2 \\ &\cdot
109^2  \cdot
131  \cdot
137  \cdot
139  \cdot
149  \cdot
151  \cdot
157  \cdot
163  \cdot
167  \cdot
173  \cdot
179  \cdot
181  \cdot
191  \cdot
193 \\ &\cdot
199  \cdot
223  \cdot
227  \cdot
229  \cdot
233 \\ &\cdot
1481  \cdot
817289119749725617  \cdot
129311862685631303361301385041,  \\
C_{252}&=\frac{-1}{29\cdot 43\cdot 127}\cdot
2^{213}  \cdot
3^{107}  \cdot
5^{44}  \cdot
7^{41}  \cdot
11^{24}  \cdot
13^{19}  \cdot
17^{11}  \cdot
19^9  \cdot
23^{10}  \cdot
31^7  \cdot
37^6 \\ &\cdot
41^6  \cdot
47^4  \cdot
53^4  \cdot
59^4  \cdot
61^3  \cdot
67^3  \cdot
73^3  \cdot
79^3  \cdot
83^3  \cdot
89^2  \cdot
97^2  \cdot
101^2  \cdot
103^2  \cdot
107^2 \\ &\cdot
109^2  \cdot
131  \cdot
137  \cdot
139  \cdot
149  \cdot
151  \cdot
157  \cdot
163  \cdot
167  \cdot
173  \cdot
179  \cdot
181  \cdot
191  \cdot
193 \\ &\cdot
199  \cdot
223  \cdot
227  \cdot
229  \cdot
233  \cdot
241  \cdot
38356192325687 \\ &\cdot
983117314012755559943  \cdot
209316897257954284206743809693.
\endalign
$$

\newpage

\fp
Each  $C_{14n}$  is written by an integer  $G_{14n}$  (von Staudt-Clausen type theorem) as follows:
$$
\alignat {11}  
C_{ 14}&                    && &&                && &&               &&                 &&                  &&                &&                &&                   && \text{\lq\lq an integer"}, \\
C_{ 28}&= \frac{ 4    }{29} && &&                && &&               &&                 &&                  &&                &&                &&                   &&+G_{28}, \\
C_{ 42}&=                   && &&\frac{10  }{43} && &&               &&                 &&                  &&                &&                &&                   &&+G_{42}, \\
C_{ 56}&= \frac{ 4^ 2 }{29} && &&                && &&               &&                 &&                  &&                &&                &&                   &&+G_{56}, \\
C_{ 70}&=                   && &&                && &&\frac{20  }{71}&&                 &&                  &&                &&                &&                   &&+G_{70}, \\
C_{ 84}&= \frac{ 4^ 3 }{29} &&+&&\frac{10^2}{43} && &&               &&                 &&                  &&                &&                &&                   &&+G_{84}, \\
C_{ 98}&                    && &&                && &&               &&                 &&                  &&                &&                &&                   && \text{\lq\lq an integer"}, \\
C_{112}&= \frac{ 4^ 4 }{29} && &&                && &&               &&+\frac{ 3  }{113}&&                  &&                &&                &&                   &&+G_{112}, \\
C_{126}&=                   && &&\frac{10^3}{43} && &&               &&                 &&+\frac{41  }{127} &&                &&                &&                   &&+G_{126}, \\
C_{140}&= \frac{ 4^ 5 }{29} && &&                &&+&&\frac{20^2}{71}&&                 &&                  &&                &&                &&                   &&+G_{140}, \\
C_{154}&                    && &&                && &&               &&                 &&                  &&                &&                &&                   && \text{\lq\lq an integer"}, \\
C_{168}&= \frac{ 4^ 6 }{29} &&+&&\frac{10^4}{43} && &&               &&                 &&                  &&                &&                &&                   &&+G_{168}, \\
C_{182}&                    && &&                && &&               &&                 &&                  &&                &&                &&                   && \text{\lq\lq an integer"}, \\
C_{196}&= \frac{ 4^ 7 }{29} && &&                && &&               &&                 &&                  &&+\frac{56 }{197}&&                &&                   &&+G_{196}, \\
C_{210}&=                   && &&\frac{10^5}{43} &&+&&\frac{20^3}{71}&&                 &&                  &&                &&+\frac{180}{211}&&                   &&+G_{210}, \\
C_{224}&= \frac{ 4^ 8 }{29} && &&                && &&               &&+\frac{ 3^2}{113}&&                  &&                &&                &&                   &&+G_{224}, \\
C_{238}&=                   && &&                && &&               &&                 &&                  &&                &&                &&\ \frac{118}{239}\ &&+G_{238}, \\
C_{252}&= \frac{ 4^ 9 }{29} &&+&&\frac{10^6}{43} && &&               &&                 &&+\frac{41^2}{127} &&                &&                &&                   &&+G_{252}. 
\endalignat 
$$

\newpage \fp
{\bf 16.5.  $\bk{x}(\bk{u})$  of the curve  $\bk{y}^2=\bk{x}^7-\bk{x}$.} \ \ \ 
First several values of  $C_{12n}$  for the curve  $y^2=x^7-x$ ($g=3$) are as follows:
$$
\allowdisplaybreaks
\align
C_{12}&=\frac{1}{13}\cdot
2^{10} \cdot
3^5    \cdot
5^2   \cdot
7,  \\
C_{24}&=\frac{-1}{13}\cdot
2^{22}   \cdot
3^{10}   \cdot
5^2   \cdot
7^3   \cdot
11^2   \cdot
17^2   \cdot
19,  \\
C_{36}&=\frac{1}{13\cdot 37}\cdot
2^{36}   \cdot
3^{17}   \cdot
5^6   \cdot
7^4   \cdot
11^3   \cdot
17^2   \cdot
19   \cdot
23   \cdot
29   \cdot
31   \cdot
41,  \\
C_{48}&=\frac{-1}{13\cdot 61}\cdot
2^{46}   \cdot
3^{22}   \cdot
5^6   \cdot
7^5   \cdot
11^4   \cdot
17^2   \cdot
19^2   \cdot
23^2   \cdot
29   \cdot
31   \cdot
41 \\ &\cdot
43   \cdot
226843,  \\
C_{60}&=\frac{1}{13\cdot 61}\cdot
2^{56}   \cdot
3^{30}   \cdot
5^{11}   \cdot
7^8   \cdot
11^5   \cdot
17^3   \cdot
19^5   \cdot
23^2   \cdot
29^2   \cdot
31   \cdot
41 \\ &\cdot
43   \cdot
47   \cdot
53   \cdot
1489,  \\
C_{72}&=\frac{-1}{13\cdot 37}  \cdot
2^{73}   \cdot
3^{34}   \cdot
5^{10}   \cdot
7^{10}   \cdot
11^6   \cdot
17^4   \cdot
19^3   \cdot
23^3   \cdot
29^2   \cdot
31^2   \cdot
41 \\ &\cdot
43   \cdot
47   \cdot
53   \cdot
59   \cdot
67   \cdot
9833   \cdot
38618729,  \\
C_{84}&=\frac{1}{13}\cdot
2^{84}   \cdot
3^{42}   \cdot
5^{14}   \cdot
7^{12}   \cdot
11^7   \cdot
17^3   \cdot
19^4   \cdot
23^3   \cdot
29^2   \cdot
31^2   \cdot
41^2 \\ &\cdot
43   \cdot
47   \cdot
53   \cdot
59   \cdot
67   \cdot
71   \cdot
79   \cdot
257   \cdot
311   \cdot
10445899,  \\
C_{96}&=\frac{-1}{13\cdot 97}\cdot
2^{94}   \cdot
3^{47}   \cdot
5^{14}   \cdot
7^{12}   \cdot
11^8   \cdot
17^4   \cdot
19^4   \cdot
23^4   \cdot
29^3   \cdot
31^3   \cdot
41^2 \\ &\cdot
43^2   \cdot
47^2   \cdot
53   \cdot
59   \cdot
67   \cdot
71   \cdot
79   \cdot
83   \cdot
89   \cdot
157   \cdot
47563   \cdot
46457   \cdot
653693,  \\
C_{108}&=\frac{1}{13\cdot 37\cdot 109}\cdot
2^{105}   \cdot
3^{53}   \cdot
5^{18}   \cdot
7^{15}   \cdot
11^9   \cdot
17^6   \cdot
19^5   \cdot
23^4   \cdot
29^3   \cdot
31^3   \cdot
41^2 \\ &\cdot
43^2   \cdot
47^2   \cdot
53^2   \cdot
59   \cdot
67   \cdot
71   \cdot
79   \cdot
83   \cdot
89   \cdot
101   \cdot
103 \\ &\cdot
647   \cdot
255307815673498109,  \\
C_{120}&=\frac{-1}{13\cdot 61}\cdot
2^{116}   \cdot
3^{59}   \cdot
5^{18}   \cdot
7^{17}   \cdot
11^9   \cdot
17^6   \cdot
19^6   \cdot
23^5   \cdot
29^4   \cdot
31^3   \cdot
41^2 \\ &\cdot
43^2   \cdot
47^2   \cdot
53^2   \cdot
59^2   \cdot
67   \cdot
71   \cdot
79   \cdot
83   \cdot
89   \cdot
101   \cdot
103   \cdot
107   \cdot
113 \\ &\cdot
331   \cdot
4759   \cdot
8486377   \cdot
205745894701,  \\
C_{132}&=\frac{1}{13}\cdot
2^{132}   \cdot
3^{63}   \cdot
5^{22}   \cdot
7^{18}   \cdot
11^{12}   \cdot
17^6   \cdot
19^6   \cdot
23^5   \cdot
29^4   \cdot
31^4   \cdot
41^3 \\ &\cdot
43^3   \cdot
47^2   \cdot
53^2   \cdot
59^2   \cdot
67   \cdot
71   \cdot
79   \cdot
83   \cdot
89   \cdot
101   \cdot
103   \cdot
107   \cdot
113 \\ &\cdot
127   \cdot
519205082831736259350305213,  \\
C_{144}&=\frac{-1}{13\cdot 37\cdot 73}\cdot
2^{145}   \cdot
3^{70}   \cdot
5^{22}   \cdot
7^{19}   \cdot
11^{13}   \cdot
17^7   \cdot
19^7   \cdot
23^6   \cdot
29^3   \cdot
31^4   \cdot
41^3 \\ &\cdot
43^3   \cdot
47^3   \cdot
53^2   \cdot
59^2   \cdot
67^2   \cdot
71^2   \cdot
79   \cdot
83   \cdot
89   \cdot
101   \cdot
103   \cdot
107   \cdot
113 \\ &\cdot
127   \cdot
131   \cdot
137   \cdot
139   \cdot
487   \cdot
116447   \cdot
70754581159   \cdot
13276324698781,  \\
C_{156}&=\frac{1}{13\cdot 157}\cdot
2^{152}   \cdot
3^{75}   \cdot
5^{26}   \cdot
7^{22}   \cdot
11^{14}   \cdot
17^8   \cdot
19^8   \cdot
23^6   \cdot
29^4   \cdot
31^4   \cdot
41^3 \\ &\cdot
43^3   \cdot
47^3   \cdot
53^2   \cdot
59^2   \cdot
67^2   \cdot
71^2   \cdot
79   \cdot
83   \cdot
89   \cdot
101   \cdot
103   \cdot
107   \cdot
113 \\ &\cdot
127   \cdot
131   \cdot
137   \cdot
139   \cdot
149   \cdot
151 \\ &\cdot
449   \cdot
9839   \cdot
2761009   \cdot
22771433   \cdot
14897675203748611,  \\
C_{168}&=\frac{-1}{13}\cdot
2^{167}  \cdot
3^{87}   \cdot
5^{26}   \cdot,
7^{25}   \cdot
11^{15}   \cdot
17^7   \cdot
19^8   \cdot
23^7   \cdot
29^4   \cdot
31^5   \cdot
41^4 \\ &\cdot
43^3   \cdot
47^3   \cdot
53^3   \cdot
59^2   \cdot
67^2   \cdot
71^2   \cdot
79^2   \cdot
83^2   \cdot
89   \cdot
101   \cdot
103   \cdot
107   \cdot
113 \\ &\cdot
127   \cdot
131   \cdot
137   \cdot
139   \cdot
149   \cdot
151   \cdot
163 \\ &\cdot
2153   \cdot
2179   \cdot
2909   \cdot
5303   \cdot
28573   \cdot
215723   \cdot
486014655083,   \\
C_{180}&=\frac{1}{13\cdot 37\cdot 61\cdot 181}\cdot 
2^{181}  \cdot
3^{90}   \cdot
5^{31}   \cdot
7^{25}   \cdot
11^{16}   \cdot
17^8    \cdot
19^9    \cdot
23^7   \cdot
29^5   \cdot
31^5   \cdot
41^4 \\ &\cdot
43^4   \cdot
47^4   \cdot
53^3   \cdot
59^3   \cdot
67^2   \cdot
71^2   \cdot
79^2   \cdot
83^2   \cdot
89^2   \cdot
101   \cdot
103   \cdot
107   \cdot
113 \\ &\cdot
127   \cdot
131   \cdot
137   \cdot
139   \cdot
149   \cdot
151   \cdot
163   \cdot
167   \cdot
173 \\ &\cdot
1409   \cdot
4356533   \cdot
41429650536998423202242258179333,  \\
C_{192}&=\frac{-1}{13\cdot 97\cdot 193}\cdot
2^{190}  \cdot
3^{97}  \cdot
5^{30}  \cdot
7^{26}  \cdot
11^{17}  \cdot
17^{10}  \cdot
19^{10}  \cdot
23^8  \cdot
29^5  \cdot
31^6  \cdot
41^4 \\ &\cdot
43^4  \cdot
47^4  \cdot
53^3  \cdot
59^3  \cdot
67^2  \cdot
71^2  \cdot
79^2  \cdot
83^2  \cdot
89^2  \cdot
101  \cdot
103  \cdot
107  \cdot
113 \\ &\cdot
127  \cdot
131  \cdot
137  \cdot
139  \cdot
149  \cdot
151  \cdot
163  \cdot
167  \cdot
173  \cdot
179 \\ &\cdot
1637  \cdot
6553  \cdot
91528069  \cdot
138182422934743639  \cdot
2002149133519714847,  \\
C_{204}&=\frac{1}{13}\cdot
2^{201}  \cdot
3^{100}  \cdot
5^{34}  \cdot
7^{29}  \cdot
11^{18}  \cdot
17^{11}  \cdot
19^{10}  \cdot
23^8  \cdot
29^6  \cdot
31^6  \cdot
41^3 \\ &\cdot
43^4  \cdot
47^4  \cdot
53^3  \cdot
59^3  \cdot
67^3  \cdot
71^2  \cdot
79^2  \cdot
83^2  \cdot
89^2  \cdot
101^2  \cdot
103  \cdot
107  \cdot
113 \\ &\cdot
127  \cdot
131  \cdot
137  \cdot
139  \cdot
149  \cdot
151  \cdot
163  \cdot
167  \cdot
173  \cdot
179  \cdot
191  \cdot
197  \cdot
199 \\ &\cdot
3215816668951  \cdot
690281420686969  \cdot
355828513600997578367,  \\
C_{216}&=\frac{-1}{13\cdot 37\cdot 73\cdot 109}\cdot
2^{213}  \cdot
3^{106}  \cdot
5^{34}  \cdot
7^{31}  \cdot
11^{19}  \cdot
17^{10}  \cdot
19^{11}  \cdot
23^9  \cdot
29^6  \cdot
31^6  \cdot
41^4 \\ &\cdot
43^4  \cdot
47^4  \cdot
53^4  \cdot
59^3  \cdot
67^3  \cdot
71^3  \cdot
79^2  \cdot
83^2  \cdot
89^2  \cdot
101^2  \cdot
103^2  \cdot
107^2  \cdot
113 \\ &\cdot
127  \cdot
131  \cdot
137  \cdot
139  \cdot
149  \cdot
151  \cdot
163  \cdot
167  \cdot
173  \cdot
179  \cdot
191  \cdot
197  \cdot
199  \cdot
211 \\ &\cdot
668273  \cdot
236655047  \cdot
8245796813443668293 \\ &\cdot
2048021738357246699741070257. 
\endalign
$$

\newpage

\fp
Each  $C_{12n}$  is written by an integer  $G_{12n}$ (von Staudt-Clausen type theorem) as follows:
$$
\alignat {9}  
C_{ 12}&=\frac{7     }{13} &&                   &&                &&                 &&                   &&                &&                &&                 &&+G_{12}, \\
C_{ 24}&=\frac{7^  2 }{13} &&                   &&                &&                 &&                   &&                &&                &&                 &&+G_{24}, \\
C_{ 36}&=\frac{7^  3 }{13} &&                   &&                &&                 &&                   &&                &&                &&                 &&+G_{36}, \\
C_{ 48}&=\frac{7^  4 }{13} &&                   &&                &&                 &&                   &&                &&                &&                 &&+G_{48}, \\
C_{ 60}&=\frac{7^  5 }{13} &&+\frac{51    }{61} &&                &&                 &&                   &&                &&                &&                 &&+G_{60}, \\
C_{ 72}&=\frac{7^  6 }{13} &&                   &&+\frac{6  }{73} &&                 &&                   &&                &&                &&                 &&+G_{72}, \\
C_{ 84}&=\frac{7^  7 }{13} &&                   &&                &&                 &&                   &&                &&                &&                 &&+G_{84}, \\
C_{ 96}&=\frac{7^  8 }{13} &&                   &&                &&+\frac{79  }{97} &&                   &&                &&                &&                 &&+G_{96}, \\
C_{108}&=\frac{7^  9 }{13} &&                   &&                &&                 &&+\frac{103  }{109} &&                &&                &&                 &&+G_{108}, \\
C_{120}&=\frac{7^{10}}{13} &&+\frac{51^ 2 }{61} &&                &&                 &&                   &&                &&                &&                 &&+G_{120}, \\
C_{132}&=\frac{7^{11}}{13} &&                   &&                &&                 &&                   &&                &&                &&                 &&+G_{132}, \\
C_{144}&=\frac{7^{12}}{13} &&                   &&+\frac{6^2}{73} &&                 &&                   &&                &&                &&                 &&+G_{144}, \\
C_{156}&=\frac{7^{13}}{13} &&                   &&                &&                 &&                   &&+\frac{22}{157} &&                &&                 &&+G_{156}, \\
C_{168}&=\frac{7^{14}}{13} &&                   &&                &&                 &&                   &&                &&                &&                 &&+G_{168}, \\
C_{180}&=\frac{7^{15}}{13} &&+\frac{51^ 3 }{61} &&                &&+\frac{79^2}{97} &&                   &&                &&+\frac{18}{181} &&                 &&+G_{180}, \\
C_{192}&=\frac{7^{16}}{13} &&                   &&                &&                 &&                   &&                &&                &&+\frac{179}{193}\ &&+G_{196}, \\
C_{204}&=\frac{7^{17}}{13} &&                   &&                &&                 &&                   &&                &&                &&                 &&+G_{204}, \\
C_{216}&=\frac{7^{18}}{13} &&                   &&+\frac{6^3}{73} &&                 &&+\frac{103^2}{109} &&                &&                &&                 &&+G_{216}. 
\endalignat 
$$

\newpage \fp
{\bf 16.6.  $\bk{x}(\bk{u})$  of the curve  $\bk{y}^2=\bk{x}^9-\bold 1$.} \ \ \ 
First several values of  $C_{18n}$  for the curve  $y^2=x^9-1$ ($g=4$) are as follows:
$$
\align
C_{18}
&=\frac1{19}\cdot
2^{16} \cdot
3^8 \cdot
5^3  \cdot
7^2  \cdot
11   \cdot
13, \\
C_{36}
&=\frac{-1}{19\cdot 37}\cdot
2^{34}  \cdot
3^{17}  \cdot
5^7     \cdot
7^4     \cdot
11^3    \cdot
13^2    \cdot
17^2    \cdot
23      \cdot
29      \cdot
31      \cdot
53,    \\
C_{54}
&=\frac1{19}\cdot
2^{51}  \cdot
3^{26}  \cdot
5^{12}  \cdot
7^8     \cdot
11^3    \cdot
13^4    \cdot
17^3    \cdot
23^3    \cdot
29^2    \cdot
31      \cdot
41      \cdot
43      \cdot
47,    \\
C_{72}
&=\frac{-1}{19\cdot 37\cdot 73}\cdot
2^{70}  \cdot
3^{34}  \cdot
5^{16}  \cdot
7^{11}  \cdot
11^5    \cdot
13^5    \cdot
17^4    \cdot
23^3    \cdot
29^2    \cdot
31^2    \cdot
41 \\  &\cdot
43      \cdot
47      \cdot
53      \cdot
59      \cdot
61      \cdot
67      \cdot
277     \cdot
35107, \\ 
C_{90}
&=\frac1{19}\cdot
2^{89}  \cdot
3^{44}  \cdot
5^{21}  \cdot
7^{13}  \cdot
11^7    \cdot
13^5    \cdot
17^5    \cdot
23^3    \cdot
29^3    \cdot
31^2    \cdot
41^2    \cdot
43^2 \\ &\cdot
47      \cdot
53      \cdot
59      \cdot
61      \cdot
67      \cdot
71      \cdot
79      \cdot
83      \cdot
113     \cdot
199     \cdot
2069,  \\ 
C_{108}
&=\frac{-1}{19\cdot 37\cdot 109}\cdot
2^{105}  \cdot
3^{53}  \cdot
5^{29}  \cdot
7^{17}  \cdot
11^7    \cdot
13^7    \cdot
17^6    \cdot
23^4    \cdot
29^3    \cdot
31^3    \cdot
41^2 \\ & \cdot
43^2    \cdot
47^2    \cdot
53^2    \cdot
59      \cdot
61^2    \cdot
67      \cdot
71      \cdot
79      \cdot
83      \cdot
89      \cdot
97      \cdot
101     \cdot
103     \cdot
251     \cdot
2681009, \\ 
C_{126}
&=\frac1{19\cdot 127}\cdot
2^{126} \cdot
3^{63}  \cdot
5^{31}  \cdot
7^{21}  \cdot
11^{10} \cdot
13^8    \cdot
17^7    \cdot
23^5    \cdot
29^5    \cdot
31^4    \cdot
41^3 \\ &\cdot
43^2    \cdot
47^2    \cdot
53^2    \cdot
59^2    \cdot
61^2    \cdot
67      \cdot
71      \cdot
79      \cdot
83      \cdot
89      \cdot
97      \cdot
101     \cdot
103     \cdot
107 \\ &\cdot
113     \cdot
22639911725059, \\
C_{144}
&=\frac{-1}{19\cdot 37\cdot 73}\cdot
2^{144} \cdot
3^{71}  \cdot
5^{35}  \cdot
7^{22}  \cdot
11^{14} \cdot
13^{11} \cdot
17^8    \cdot
23^6    \cdot
29^3    \cdot
31^4    \cdot
41^3 \\ &\cdot
43^3    \cdot
47^3    \cdot
53^2    \cdot
59^2    \cdot
61^2    \cdot
67^2    \cdot
71^2    \cdot
79      \cdot
83      \cdot
89      \cdot
97      \cdot
101     \cdot
103     \cdot
107 \\ &\cdot
113     \cdot
131     \cdot
137     \cdot
139     \cdot
257     \cdot
2063    \cdot
14375966528323,  \\
C_{162}
&=\frac1{19\cdot 163}\cdot
2^{160} \cdot
3^{84}  \cdot
5^{39}  \cdot
7^{26}  \cdot
11^{12} \cdot
13^{11} \cdot
17^9    \cdot
23^7    \cdot
29^4    \cdot
31^5    \cdot
41^3 \\ &\cdot
43^3    \cdot
47^3    \cdot
53^3    \cdot
59^2    \cdot
61^2    \cdot
67^2    \cdot
71^2    \cdot
79^2    \cdot
83      \cdot
89      \cdot
97      \cdot
101     \cdot
103     \cdot
107 \\ &\cdot
113     \cdot
131     \cdot
137     \cdot
139     \cdot
149     \cdot
151     \cdot
157     \cdot
683     \cdot
1868051 \cdot
4779877 \cdot
5414895901,\\  
C_{180}
&=\frac{-1}{19\cdot 37\cdot 181}\cdot
2^{180} \cdot
3^{92}  \cdot
5^{45}  \cdot
7^{28}  \cdot
11^{14} \cdot
13^{12} \cdot
17^{10} \cdot
23^7    \cdot
29^5    \cdot
31^5    \cdot
41^4 \\ &\cdot
43^4    \cdot
47^3    \cdot
53^3    \cdot
59^3    \cdot
61^2    \cdot
67^2    \cdot
71^2    \cdot
79^2    \cdot
83^2    \cdot
89^2    \cdot
97      \cdot
101     \cdot
103     \cdot
107 \\ &\cdot
113     \cdot
131     \cdot
137     \cdot
139     \cdot
149     \cdot
151     \cdot
157     \cdot
167     \cdot
173     \cdot
179 \\ &\cdot
80051837 \cdot
3425819816490635592983, \\
C_{198}
&=\frac1{19\cdot 199}\cdot
2^{196} \cdot
3^{99}  \cdot
5^{47}  \cdot
7^{32}  \cdot
11^{18} \cdot
13^{14} \cdot
17^{11} \cdot
23^8    \cdot
29^7    \cdot
31^6    \cdot
41^4 \\ &\cdot
43^4    \cdot
47^4    \cdot
53^3    \cdot
59^3    \cdot
61^3    \cdot
67^2    \cdot
71^2    \cdot
79^2    \cdot
83^2    \cdot
89^2    \cdot
97^2    \cdot
101     \cdot
103     \cdot
107 \\ &\cdot
113     \cdot
131     \cdot
137     \cdot
139     \cdot
149     \cdot
151^2   \cdot
157     \cdot
167     \cdot
173     \cdot
179     \cdot
191     \cdot
193     \cdot
197 \\ &\cdot
3943    \cdot
64682082077 \cdot
9591335434411.
\endalign
$$

\newpage

\fp
Each  $C_{18n}$  is written by an integer  $G_{18n}$  (von Staudt-Clausen type theorem) as follows:
$$
\alignat {10}  
C_{ 18}&=\frac{7     }{19} && &&              && &&               &&                 &&                  &&                  &&                  &&+G_{18}, \\
C_{ 36}&=\frac{7^  2 }{19} &&+&&\frac{5  }{37}&& &&               &&                 &&                  &&                  &&                  &&+G_{36}, \\
C_{ 54}&=\frac{7^  3 }{19} && &&              && &&               &&                 &&                  &&                  &&                  &&+G_{54}, \\
C_{ 72}&=\frac{7^  4 }{19} &&+&&\frac{5^2}{37}&&+&&\frac{67  }{73}&&                 &&                  &&                  &&                  &&+G_{76}, \\
C_{ 90}&=\frac{7^  5 }{19} && &&              && &&               &&                 &&                  &&                  &&                  &&+G_{90}, \\
C_{108}&=\frac{7^  6 }{19} &&+&&\frac{5^3}{37}&& &&               &&                 &&                  &&                  &&                  &&+G_{108}, \\
C_{126}&=\frac{7^  7 }{19} && &&              && &&               &&+\frac{95 }{109} &&                  &&                  &&                  &&+G_{126}, \\
C_{144}&=\frac{7^  8 }{19} &&+&&\frac{5^4}{37}&&+&&\frac{67^2}{73}&&                 &&                  &&                  &&                  &&+G_{144}, \\
C_{162}&=\frac{7^  9 }{19} && &&              && &&               &&                 &&+\frac{ 7  }{163} &&                  &&                  &&+G_{162}, \\
C_{180}&=\frac{7^{10}}{19} &&+&&\frac{5^5}{37}&& &&               &&                 &&                  &&+\frac{20   }{181}&&                  &&+G_{180}, \\
C_{198}&=\frac{7^{11}}{19} && &&              && &&               &&                 &&                  &&                  &&+\frac{114  }{199}&&+G_{198}.
\endalignat 
$$

\newpage \fp
\inbox{\bf 17}\ \ {\bf Non-hyperelliptic curves}\fp
\vskip 5pt
\fp
{\bf 17.1. Algebraic curves with completely ramified at infinity.}\ \ \ 
We here describe how the natural variable  $u$  is chosen 
for a general curve of cyclotomic type.

Let  $a$  and  $b$  be coprime pair of positive integers, and let
  $$
  f(x,y)=y^a-x^b-\sum_{(i,j)}\lambda_{ia+jb}x^iy^j
  \tag{17.1.1}
  $$
be separable polynomial, where the pair  $(i,j)$  runs 
through the integers such that
  $$
  0\leqq i<b-1, \ \ 0\leqq j<a-1, \ \ ia+jb<ab
  \tag{17.1.2}
  $$
We discuss with the algebraic curve defined by
  $$
  C:f(x,y)=0.
  \tag{17.1.3}
  $$
Here  $C$  is regarded naturally as a curve with unique point  $\infty$  at infinity. 
The genus of  $C$  is given by  $g=(a-1)(b-1)/2$.   
On this curve  $C$, the set
     $$
     \dfrac{x^{i-1}y^{a-j-1}dx}{f_y(x,y)}
     \tag{17.1.4}
     $$
forms a basis of the differential forms of the first kind, 
where  $f_y(x,y)=\frac{\partial f}{\partial y}(x,y)$.  
By choosing a generator of the fundamental group of  $C$, 
we define the period matrix  $\left[\omega'\ \omega''\right]$.   
We also define the lattice of periods in  $\bold C^g$  by
     $$
     \Lambda
     := \omega' \tp{\left[\matrix \bold Z & \bold Z & \cdots & \bold Z
                  \endmatrix \right]}
     +  \omega''\tp{\left[\matrix \bold Z & \bold Z & \cdots & \bold Z
                  \endmatrix \right]}\ (\subset \bold C^g).
     $$
We denote the Jacobian variety of  $C$  by  $J$, and the symmetric product of  $g$  copies
of  $C$  by  $\text{Sym}^g(C)$.  
Then we have a birational map
  $$
  \align
  \text{Sym}^g(C) &\rightarrow \text{Pic}^{\circ}(C)=J \\
  (P_1, \ldots, P_g)&\mapsto \text{the class of}\ P_1+\cdots +P_g-g\cdot\infty.
  \endalign 
  $$ 
As an analytic manifold, $J$  is identified with  $\bold C^g/\Lambda$.   
We denote by  $\kappa$  the natural map  $\bold C^g \rightarrow \bold C^g/\Lambda=J$.  
The map
  $$
  \iota : Q \mapsto Q-\infty
  $$
gives an embedding of  $C$  into  $J$.  
The pull-back  $\kappa^{-1}\iota(C)$  of the image of  $\iota$  by  $\kappa$  is 
a universal Abelian covering of  $C$.   
The birational map is represented analytically by sending 
each  $(P_1, \ldots, P_g)\in\text{Sym}^g(C)$  to the point  $\bk{u}\mod{\Lambda}\in\bold  C^g/\Lambda$, 
where
  $$
  \bk{u}=(u_1, \ldots, u_g)=
  \left(\int_{\infty}^{P_1}+\cdots+\int_{\infty}^{P_g}\right)
  (\omega_1, \ \ldots,  \ \omega_g).
  $$
For each point  $\bk{u}\in\kappa^{-1}\iota(C)$, 
we denote by
  $$
  x(\bk{u}), \qquad y(\bk{u})
  \tag{17.1.5}
  $$
the value of  $(x, y)$-coordinate of  $C$  
such that  $\kappa(\bk{u})=\iota\left(x(\bk{u}), \ y(\bk{u})\right)$  
We regard a rational expression of  $x(\bk{u})$  and  $y(\bk{u})$  as 
a function on  $\kappa^{-1}\iota(C)$.   
Basically, we have
\vskip 3pt
\fp
\lbox{75}{360}
{{\bf Lemma 17.1.6.} \ \ \it
The Laurent development of  $x(\bk{u})$  and  $y(\bk{u})$  
around  $\bk{u}=(0, \cdots, 0)$  with respect to  $u_g$   is of the form
  $$
  x(\bk{u})=\frac 1{{u_g}^a}+\Dg{u_g}{-a+1}, \quad
  y(\bk{u})=-\frac 1{{u_g}^b}+\Dg{u_g}{-b+1}.
  $$
}
\fp
Moreover, we have a result correspond to 4.2.2, and we see that
\vskip 3pt
\fp
\lbox{41}{360}
{It is natural to take  $u:=u_g$  as a local parameter around  $\bk{u}=(0, \cdots, 0)$  
on   $\iota^{-1}\kappa(C)$.  
}
\vskip 3pt
\fp
We call a curve  $C$  defined by  $f(x, y)=0$, where
  $$
  f(x,y)=y^a-x^b+1, \ \ \text{or} \ \  f(x,y)=y^a-x^b+x,
  \tag{17.1.7}
  $$
to be a curve {\it of cyclotomic type}.  
For such a curve  $C$, the Hurwitz coefficients of 
the power series development of the functions  $u=u_g\mapsto x(\bk{u})$, $u_g\mapsto y(\bk{u})$  
with respect to  $u$  satisfy the Clarke type theorem (von Staudt-Clausen type theorem 
plus the extension of von Staudt second theorem) and the Kummer type theorem.  
We can prove them by the same argument of this paper.

\newpage \fp
{\bf 17.2. $\bk{x}(\bk{u})$  of the curve  $\bk{y}^3=\bk{x}^5-\bold 1$.} \ \ \ 
Let consider the curve  $y^3=x^5-1$ ($g=4$).  
Its local parameter defined in 17.1 is
  $$
  u_4=\int_{\infty}^{(x,y)}\frac{x^2}{3y^2}. 
  $$
The function  $x(\bk{u})$  for this curve satisfies
  $$
  \frac{du_4}{dx}=\frac{x^2dx}{3y^2},
  $$
namely
\vskip 3pt
\fp
\lbox{45}{360}{
  $$
  \ \ \ x(\bk{u})^6x'(\bk{u})^3=27(x(\bk{u})^5-1)^2
   \ \ \ \ \text{($'$\  means  $\frac{d}{du_4}$)}.
  $$
}\fp
\vskip 3pt
\fp

\newpage

\fp
First several values of the Hurwitz coefficients  $C_m$  of 
the power series development of  $u\mapsto x(\bk{u})$  are as follows. 
Here if  $15\hskip -2pt\not|\hskip 2pt m$, then  $C_m=0$.  
$$
\allowdisplaybreaks
\align
C_{15}&=
 2^7  \cdot
 3^6  \cdot
 5^3  \cdot
 7   \cdot
 11,  \\
C_{30}&=\frac{1}{31} \cdot
 2^{19}  \cdot
 3^{14}  \cdot
 5^7  \cdot
 7^3  \cdot
 11^2  \cdot
 13^2  \cdot
 17   \cdot
 19^2  \cdot
 23,  \\
C_{45}&=
 2^{31}  \cdot
 3^{21}  \cdot
 5^{10}  \cdot
 7^6  \cdot
 11^3  \cdot
 13^3  \cdot
 17^2  \cdot
 19^2  \cdot
 29   \cdot
 37   \cdot
 41   \cdot
 787,  \\
C_{60}&=\frac{1}{31\cdot 61} \cdot
 2^{43}  \cdot
 3^{28}  \cdot
 5^{15}  \cdot
 7^9  \cdot
 11^5  \cdot
 13^5  \cdot
 17^3  \cdot
 19^3  \cdot
 23   \cdot
 29   \cdot
 37^2  \cdot
 41   \cdot
 43 \\ &\cdot
 47   \cdot
 53   \cdot
 2671,  \\
C_{75}&=
 2^{55}  \cdot
 3^{37}  \cdot
 5^{18}  \cdot
 7^{11}  \cdot
 11^6  \cdot
 13^5  \cdot
 17^4  \cdot
 19^2  \cdot
 23^2  \cdot
 29^2  \cdot
 37^2  \cdot
 41   \cdot
 43 \\ &\cdot
 47   \cdot
 53   \cdot
 59   \cdot
 67   \cdot
 71   \cdot
 149   \cdot
 5843,  \\
C_{90}&=\frac{1}{31} \cdot
 2^{67}  \cdot
 3^{44}  \cdot
 5^{21}  \cdot
 7^{13}  \cdot
 11^9  \cdot
 13^5  \cdot
 17^5  \cdot
 19^3  \cdot
 23   \cdot
 29^4  \cdot
 37^2  \cdot
 41^2  \cdot
 43^2 \\ &\cdot
 47   \cdot
 53   \cdot
 59   \cdot
 67   \cdot
 71   \cdot
 73   \cdot
 79   \cdot
 83   \cdot
 482561869,  \\
C_{105}&=
 2^{79}  \cdot
 3^{51}  \cdot
 5^{25}  \cdot
 7^{17}  \cdot
 11^9  \cdot
 13^6  \cdot
 17^6  \cdot
 19^5  \cdot
 23^2  \cdot
 29^3  \cdot
 37^2  \cdot
 41^2  \cdot
 43^2 \\ &\cdot
 47^2 \cdot
 59   \cdot
 67   \cdot
 71   \cdot
 73   \cdot
 79   \cdot
 83   \cdot
 89   \cdot
 97   \cdot
 101   \cdot
 4132022673901,  \\
C_{120}&=\frac{1}{31\cdot 61} \cdot
 2^{91}  \cdot
 3^{59}  \cdot
 5^{29}  \cdot
 7^{18}  \cdot
 11^9  \cdot
 13^8  \cdot
 17^6  \cdot
 19^5  \cdot
 23^3  \cdot
 29^4  \cdot
 37^3  \cdot
 41^2  \cdot
 43^2 \\ &\cdot
 47^2 \cdot
 53   \cdot
 59   \cdot
 67   \cdot
 71   \cdot
 73   \cdot
 79   \cdot
 83   \cdot
 89   \cdot
 97   \cdot
 101   \cdot
 103   \cdot
 107   \cdot
 109   \cdot
 113 \\ &\cdot
 113153   \cdot
 31199395512997,  \\
C_{135}&=
 2^{103}  \cdot
 3^{66}  \cdot
 5^{34}  \cdot
 7^{20}  \cdot
 11^{12}  \cdot
 13^9  \cdot
 17^7  \cdot
 19^6  \cdot
 23^2  \cdot
 29^4  \cdot
 37^3  \cdot
 41^3  \cdot
 43^3 \\ &\cdot
 47^2 \cdot
 53   \cdot
 59^2  \cdot
 67   \cdot
 71   \cdot
 73   \cdot
 79   \cdot
 83^2  \cdot
 89   \cdot
 97   \cdot
 101   \cdot
 103   \cdot
 107   \cdot
 109   \cdot
 113   \cdot
 127 \\ &\cdot
 131    \cdot
 1567   \cdot
 970351   \cdot
 280656653,  \\
C_{150}&=\frac{1}{31\cdot 151} \cdot
 2^{115}  \cdot
 3^{73}  \cdot
 5^{37}  \cdot
 7^{24}  \cdot
 11^{13}  \cdot
 13^{11}  \cdot
 17^8  \cdot
 19^5  \cdot
 23^3  \cdot
 29^5  \cdot
 37^3  \cdot
 41^3  \cdot
 43^3 \\ &\cdot
 47^3 \cdot
 53   \cdot
 59^2  \cdot
 67^2  \cdot
 71^2  \cdot
 73^2  \cdot
 79   \cdot
 83   \cdot
 89   \cdot
 97   \cdot
 101   \cdot
 103   \cdot
 107   \cdot
 109   \cdot
 113   \cdot
 127 \\ &\cdot
 131   \cdot
 137   \cdot
 139   \cdot
 1993   \cdot
 13703297   \cdot
 394936930338121.
\endalign
$$

\newpage

\fp
Each  $C_{15n}$  is written by an integer  $G_{15n}$  as follows 
(von Staudt-Clausen type theorem):
$$
\alignat {7}
C_{ 15}&                      &&   &&                       &&   &&                   &&   && \text{\lq\lq an integer"}, \\
C_{ 30}&=\frac{2         }{31}&&   &&                       &&   &&                   &&\ +&& G_{30}, \\
C_{ 45}&                      &&   &&                       &&   &&                   &&   && \text{\lq\lq an integer"}, \\
C_{ 60}&=\frac{2\cdot 4  }{31}&&\ +&&\frac{6           }{61}&&   &&                   &&\ +&& G_{60}, \\
C_{ 75}&                      &&   &&                       &&   &&                   &&   && \text{\lq\lq an integer"}, \\
C_{ 90}&=\frac{2\cdot 4^2}{31}&&   &&                       &&   &&                   &&\ +&& G_{90}, \\
C_{105}&                      &&   &&                       &&   &&                   &&   && \text{\lq\lq an integer"}, \\
C_{120}&=\frac{2\cdot 4^3}{31}&&\ +&&\frac{11\cdot 12  }{61}&&   &&                   &&\ +&& G_{120}, \\
C_{135}&                      &&   &&                       &&   &&                   &&   && \text{\lq\lq an integer"}, \\
C_{150}&=\frac{2\cdot 4^4}{31}&&   &&                       &&\ +&&\frac{136    }{151}&&\ +&& G_{150}.
\endalignat 
$$

\newpage \fp
{\bf 17.3. $\bk{x}(\bk{u})$  of the curve  $\bk{y}^4=\bk{x}^3-\bk{x}$.} \ \ \ 
Let consider the curve  $y^4=x^3-x$ ($g=3$).  
Its local parameter defined in 17.1 is
  $$
  u_3=\int_{\infty}^{(x,y)}\frac{xdx}{4y^3}.
  $$
The function  $x(\bk{u})$  for this curve satisfies
  $$
  \frac{du_3}{dx}=\frac{x}{4y^3},
  $$
namely, 
\vskip 3pt
\fp
\lbox{45}{360}{
  $$
  \ \ \ x(\bk{u})x'(\bk{u})^4=4^4(x^2(\bk{u})-1)^3
   \ \ \ \ \text{($'$\  means  $\frac{d}{du_3}$)}.
  $$
}\fp
\vskip 3pt
\fp

\newpage

\fp
First several values of the Hurwitz coefficients  $C_m$  
of the power series development of  $u\mapsto x(\bk{u})$  are as follows.  
Here if  $8\hskip -2pt\not|\hskip 2pt m$, then  $C_m=0$.  
$$
\allowdisplaybreaks
\align
C_{8}&=
 2^6,  \\
C_{16}&=\frac{-1}{17}\cdot
 2^{16}  \cdot
 3^3  \cdot
 5^2  \cdot
 7    \cdot
 11,   \\
C_{24}&=
 2^{21}  \cdot
 3^6  \cdot
 5^2  \cdot
 7^4  \cdot
 11    \cdot
 13    \cdot
 19,   \\
C_{32}&=\frac{-1}{17}\cdot
 2^{30}  \cdot
 3^9  \cdot
 5^4  \cdot
 7^4  \cdot
 11^3  \cdot
 13^2  \cdot
 19    \cdot
 23    \cdot
 53,   \\
C_{40}&=\frac{1}{41}\cdot
 2^{38}  \cdot
 3^{12}  \cdot
 5^7  \cdot
 7^5  \cdot
 11^3  \cdot
 13^2  \cdot
 19    \cdot
 23    \cdot
 29    \cdot
 31    \cdot
 10079,   \\
C_{48}&=\frac{-1}{17}\cdot
 2^{48}  \cdot
 3^{15}  \cdot
 5^6  \cdot
 7^5  \cdot
 11^4  \cdot
 13^3  \cdot
 19^2  \cdot
 23    \cdot
 29    \cdot
 31    \cdot
 37    \cdot
 43    \cdot
 631    \cdot
 743,   \\
C_{56}&=
 2^{52}  \cdot
 3^{18}  \cdot
 5^8  \cdot
 7^8  \cdot
 11^4  \cdot
 13^4  \cdot
 19^2  \cdot
 23^2  \cdot
 29    \cdot
 31    \cdot
 37    \cdot
 43    \cdot
 47    \cdot
 6278521,   \\
C_{64}&=\frac{-1}{17}\cdot
 2^{62}  \cdot
 3^{21}  \cdot
 5^{10}  \cdot
 7^8  \cdot
 11^5  \cdot
 13^3  \cdot
 19^3  \cdot
 23^2  \cdot
 29^2  \cdot
 31    \cdot
 37    \cdot
 43    \cdot
 47 \\ &\cdot
 53    \cdot
 59    \cdot
 4383871061,   \\
C_{72}&=\frac{1}{73}\cdot
 2^{70}  \cdot
 3^{24}  \cdot
 5^{10}  \cdot
 7^{10}  \cdot
 11^7  \cdot
 13^4  \cdot
 19^3  \cdot
 23^3  \cdot
 29^2  \cdot
 31^2  \cdot
 37    \cdot
 43^2  \cdot
 47 \\ &\cdot
 53    \cdot
 59    \cdot
 61    \cdot
 67    \cdot
 131    \cdot
 262369,   \\
C_{80}&=\frac{-1}{17\cdot 41}\cdot
 2^{80}  \cdot
 3^{27}  \cdot
 5^{13}  \cdot
 7^{11}  \cdot
 11^7  \cdot
 13^4  \cdot
 19^4  \cdot
 23^4  \cdot
 29^2  \cdot
 31^2  \cdot
 37^2  \cdot
 43    \cdot
 47 \\ &\cdot
 53    \cdot
 59    \cdot
 61    \cdot
 67    \cdot
 71    \cdot
 439    \cdot
 547    \cdot
 667333,   \\
C_{88}&=\frac{1}{89}\cdot
 2^{84}  \cdot
 3^{30}  \cdot
 5^{14}  \cdot
 7^{12}  \cdot
 11^8  \cdot
 13^5  \cdot
 19^4  \cdot
 23^3  \cdot
 29^2  \cdot
 31^2  \cdot
 37^2  \cdot
 43    \cdot
 47 \\ &\cdot
 53    \cdot
 59    \cdot
 61    \cdot
 67    \cdot
 71    \cdot
 79    \cdot
 83    \cdot
 1121453    \cdot
 21825851419,   \\
C_{96}&=\frac{-1}{17\cdot 97}\cdot
 2^{96}  \cdot
 3^{33}  \cdot
 5^{14}  \cdot
 7^{12}  \cdot
 11^8  \cdot
 13^6  \cdot
 19^4  \cdot
 23^4  \cdot
 29^3  \cdot
 31^3  \cdot
 37^2  \cdot
 43^2  \cdot
 47 \\ &\cdot
 53    \cdot
 59    \cdot
 61    \cdot
 67    \cdot
 71    \cdot
 79    \cdot
 83    \cdot
 139    \cdot
 20840707    \cdot
 194526553349,   \\
C_{104}&=
 2^{104}  \cdot
 3^{36}  \cdot
 5^{16}  \cdot
 7^{14}  \cdot
 11^9  \cdot
 13^7  \cdot
 19^5  \cdot
 23^4  \cdot
 29^3  \cdot
 31^3  \cdot
 37^2  \cdot
 43^2  \cdot
 47^2 \\ &\cdot
 53    \cdot
 59    \cdot
 61    \cdot
 67    \cdot
 71    \cdot
 79    \cdot
 83    \cdot
 1993    \cdot
 126570133    \cdot
 4007528953,   \\
C_{112}&=\frac{-1}{17\cdot 113}\cdot
 2^{113}  \cdot
 3^{39}  \cdot
 5^{18}  \cdot
 7^{16}  \cdot
 11^9  \cdot
 13^7  \cdot
 19^5  \cdot
 23^4  \cdot
 29^3  \cdot
 31^3  \cdot
 37^2  \cdot
 43^2  \cdot
 47^2 \\ &\cdot
 53^2  \cdot
 59    \cdot
 61    \cdot
 67    \cdot
 71    \cdot
 79    \cdot
 83    \cdot
 101    \cdot
 103    \cdot
 107    \cdot
 167  \\ &\cdot
 70855976395479100717919,   \\
C_{120}&=\frac{1}{41}\cdot
 2^{116}  \cdot
 3^{42}  \cdot
 5^{18}  \cdot
 7^{17}  \cdot
 11^9  \cdot
 13^8  \cdot
 19^6  \cdot
 23^5  \cdot
 29^4  \cdot
 31^3  \cdot
 37^3  \cdot
 43^2  \cdot
 47^2 \\ &\cdot
 53^2  \cdot
 59    \cdot
 61    \cdot
 67    \cdot
 71    \cdot
 79    \cdot
 83    \cdot
 101    \cdot
 103    \cdot
 107    \cdot
 109 \\ &\cdot
 541    \cdot
 117101    \cdot
 1595557    \cdot
 27462441284549.
\endalign
$$

\newpage 

\fp
Each  $C_{8n}$  is written by an integer  $G_{8n}$  as follows 
(von Staudt-Clausen type theorem):
$$
\alignat {11}
C_{  8}&                        &&   &&                       && &&              &&                 && &&               &&                 &&                   && \text{\lq\lq an integer"}, \\
C_{ 16}&=\frac{11          }{17}&&   &&                       && &&              &&                 && &&               &&                 &&                   &&+G_{16}, \\
C_{ 24}&                        &&   &&                       && &&              &&                 && &&               &&                 &&                   && \text{\lq\lq an integer"}, \\
C_{ 32}&=\frac{11\cdot 15  }{17}&&   &&                       && &&              &&                 && &&               &&                 &&                   &&+G_{32}, \\
C_{ 40}&=                       &&   &&\frac{29          }{41}&& &&              &&                 && &&               &&                 &&                   &&+G_{40}, \\
C_{ 48}&=\frac{11\cdot 15^2}{17}&&   &&                       && &&              &&                 && &&               &&                 &&                   &&+G_{48}, \\
C_{ 56}&                        &&   &&                       && &&              &&                 && &&               &&                 &&                   && \text{\lq\lq an integer"}, \\
C_{ 64}&=\frac{11\cdot 15^3}{17}&&   &&                       && &&              &&                 && &&               &&                 &&                   &&+G_{64}, \\
C_{ 72}&=                       &&   &&                       && &&\frac{1  }{73}&&                 && &&               &&                 &&                   &&+G_{72}, \\
C_{ 80}&=\frac{11\cdot 15^4}{17}&&\ +&&\frac{29\cdot 10^2}{41}&& &&              &&                 && &&               &&                 &&                   &&+G_{80}, \\
C_{ 88}&=                       &&   &&                       && &&              &&\ \frac{28  }{89}&& &&               &&                 &&                   &&+G_{88}, \\
C_{ 96}&=\frac{11\cdot 15^5}{17}&&   &&                       && &&              &&                 &&+&&\frac{94  }{97}&&                 &&                   &&+G_{96}, \\
C_{104}&                        &&   &&                       && &&              &&                 && &&               &&                 &&                   && \text{\lq\lq an integer"}, \\
C_{112}&=\frac{11\cdot 15^6}{17}&&   &&                       && &&              &&                 && &&               &&+\frac{73  }{113}&&                   &&+G_{112}, \\
C_{120}&=                       &&   &&\frac{29\cdot 10^3}{41}&& &&              &&                 && &&               &&                 &&                   &&+G_{120}.
\endalignat 
$$

\newpage \fp
\inbox{\bf 18}\ \ 
{\bf Appendices}\ \ 
\vskip 5pt
\fp
{\bf 18.1. Relations on binomial coefficients.} \ \ \ 
The Hasse-Witt matrix with respect to the basis (4.1.1) 
for the curve  $y^2=x^5-1$  mod  $p$  is diagonal (\cite{{\bf Yu}}), 
and its  $(2,2)$-entry is, as in 6.3, given by
  $$
  A_p
  =(-1)^{(p-1)/10}\cdot
   \dbinom{(p-1)/2}{(p-1)/10} \mod{p}.
  $$
We show that it coincide with the value (11.1.4).  
In other words, we claim, for a prime  $p\equiv 1$ mod $5$  
by letting  $p=10m+1$, that
  $$
  \frac{(2m-1)!!}{m!2^m}
  \equiv (-1)^m\binom{5m}{m} \mod{p}.
  $$
Since this is equivalent to
  $$
  \frac{(2m-1)!!}{m!}
  \equiv (-2)^m\binom{5m}{m} \mod{p}
  $$
we prove this. Because
  $$
  2(5m),\ \ 2(5m-1),\cdots, 2(5m-(m-1))  \mod{p}
  $$
are coincide with
  $$
  -1, \ \ -3, \cdots, \ -(2m-1) \mod{p}
  $$
respectively. 
Thus, we see the equation above.   
Similarly, we can prove
  $$
  \frac{(5m-3)!!}{m!(-5)^m}
  \equiv (-1)^m\binom{5m}{m} \mod{p}.
  $$

Although, the facts above are proved by using basic properties of 
$p$-adic $\varGamma$-function, we have described them quite elementary way.

\newpage \fp
{\bf 18.2. Links with certain Eisenstein type series without justification.}\ \ \ 
It is a natural question whether the numbers  $C_n$  and  $D_n$  relate a kind of $L$-function.  
In the below, we explain that such a question is not completely nonsense.  

We recall the proof of Hurwitz formula, that is
  $$
  \sum_{\lambda\in \bold Z+\bold Z\sqrt{-1}, \hskip 1pt \neq 0}
  \frac 1{\lambda^{4n}}
  =\frac{\varpi^{4n}}{(4n)!}2^{4n}E_{4n}
  \tag{18.2.1}
  $$
along \cite{{\bf Hu1}}(see also \cite{{\bf AIK}, pp.193-198}), where
  $$
  \varpi
  =\int_1^{\infty}\frac{dx}{\sqrt{x^3-x}}
  (>0).
  $$
After equating the two kinds of developments
  $$
  \aligned
  \wp(u)&=\frac 1{u^2}
  + \sum_{n=2}^{\infty}\frac{2^nE_n}{n}\frac{u^{n-2}}{(n-2)!},\\
  \wp(u)&=\frac 1{u^2}
  +\hskip -10pt 
  \sum_{\ell\in \bold Z\varpi+\bold Z\varpi\sqrt{-1}, \hskip 1pt \neq 0}
  \hskip -3pt
  \left(\frac 1{(u-\ell)^2}-\frac 1{\ell^2}\right),
  \endaligned
  \tag{18.2.2}
  $$
by removing  $1/u^2$, and by differentiating them  $4n-2$ times, 
we have (18.2.1) by substituting  $u=0$.  
This argument is similar to the proof of Euler's formula,
namely the expression of the special value  $\zeta(2m)$  of the Riemann zeta function 
by a Bernoulli number and  $\pi$, by using two kinds of developments for  $1/\sin^2(u)$.  

For the case of a hyperelliptic curve of arbitrary genus  $g$, $x(\bk{u})$  
has Laurent expansion around each lattice 
point  $\ell=(\ell_1, \ell_2, \cdots, \ell_g)\in\Lambda$ ($\subset \bold C^g$):
  $$
  x(\bk{u})=\frac 1{(u_g-\ell_g)^2}+\cdots.
  \tag{18.2.3}
  $$
Although we do not have any justification, we could expect such the formula that
  $$
  x(\bk{u})=\frac 1{{u_g}^2} 
  + \sumaste_{\ell\in\Lambda,\hskip 1pt \neq 0}
  \left(\frac 1{(u_g-\ell_g)^2}-\frac 1{{\ell_g}^2}\right),
  \tag{18.2.4}
  $$
where  {\large *}  means the sum is not justified.  
If such the formula exists by the similar argument as in  $\wp(u)$, 
we have the following type formula:
\vskip 3pt
\fp
\lbox{65}{360}
{For the curve  $y^2=x^{2g+1}-1$, 
  $$
  \sumaste_{\lambda\in L,\hskip 1pt \neq 0}
  \frac 1{\lambda^{2(2g+1)n}}
  =\frac{\Omega^{2(2g+1)n}}{(2(2g+1)n)!}C_{2(2g+1)n}
  \tag{18.2.5}
  $$
}
\fp
where
  $$
  \dsize\Omega=\int_1^{\infty}\frac{x^{g-1}dx}{y}>0
  $$
and  $L$  is  $1/\Omega$  times the image of the projection of  $\Lambda\subset\bold C^g$  
to the  $g$-th factor  $\bold C$ 
(For example, if  $2g+1$  is a prime, then  $L=\bold Z[e^{2\pi i/(2g+1)}]$).  

Furthermore,  
\vskip 3pt
\fp
\lbox{64}{360}
{For the curve  $y^2=x^{2g+1}-x$, 
  $$
  \sumaste_{\lambda\in L, \hskip 1pt \neq 0}
  \frac 1{\lambda^{4gn}}
  =\frac{\Omega^{4gn}}{(4gn)!}C_{4gn}
  \tag{18.2.6}
  $$
}
\fp
where
  $$
  \Omega=\int_1^{\infty}\frac{x^{g-1}dx}{y}>0
  $$ 
and  $L$  is  $1/\Omega$  times the image of the projection of  $\Lambda\subset\bold C^g$  
to the  $g$-th factor  $\bold C$.  

If such results is obtained, it would be a discovery of new $L$-functions 
which is in the different category of Hecke's $L$-function associated to Gr\"ossen characters. 
The author also expect that this story relates to the thing written in \cite{{\bf I}, p.240}.  
\vskip 50pt

\newpage \fp
{\bf 18.3. Related formula with Hochschilt's formula.}\ \ \ 
The following formula is similar to 13.2.2.  
Although this is weaker than 13.2.2, it is much general.  
\vskip 5pt
\fp
\lbox{88}{360}{{\bf Proposition 18.3.1.}\ \ \it 
Let  $p$  be a prime.  
Let  $K$  be a function field over the finite field  $\hskip 2pt\bold F_p$  of  $p$  elements.  
We denote by  $K^p$  the set of  $p$-powers of the elements of  $K$.  
Let  $\varphi$  be a generator of the extension  $K/K^p$, 
namely, an element such that  $K=K^p(\varphi)$.  
Suppose  $D$  be a derivation from  $K$  to itself. 
Then
  $$
  \align
  D^{p-1}\varphi&=a+b\varphi, \ \ \ (a,\ b\in K^p). 
  \endalign
  $$
}
\vskip 1pt
\fp
{\it Proof. } \  
The field  $K$  is a vector space over  $K^p$  of dimension  $p$, 
and  $1$, $\varphi$, $\varphi^2$, $\cdots$, $\varphi^{p-1}$  form 
a basis of this space.  
So we may write
  $$
  D\varphi=c_0+c_1\varphi+c_2\varphi^2+\cdots+c_{p-1}\varphi^{p-1}, \ \ \ \ (c_j\in K^p).
  $$
Since  $D(\varphi^j)=j\varphi^{j-1}D\varphi$,  we have
  $$
  D
  \left[\matrix
  1 & \varphi & \varphi^2 & \cdots & \varphi^{p-1}
  \endmatrix\right]
  =
  \left[\matrix
  1 & \varphi & \varphi^2 & \cdots & \varphi^{p-1}
  \endmatrix\right]
  T,
  $$
where
  $$
  T=
  \left[
  \matrix
  0     & c_0     & 2c_{p-1}\varphi^p & 3c_{p-2}\varphi^p & 4c_{p-3}\varphi^p  & \cdots & (p-1)c_2\varphi^p     \\
  0     & c_1     & 2c_0              & 3c_{p-1}\varphi^p & 4c_{p-2}\varphi^p  & \cdots & (p-1)c_3\varphi^p     \\
  0     & c_2     & 2c_1              & 3c_0              & 4c_{p-1}\varphi^p  & \cdots & (p-1)c_4\varphi^p     \\
  0     & c_3     & 2c_2              & 3c_1              & 4c_0               & \cdots & (p-1)c_5\varphi^p     \\
  0     & c_4     & 2c_3              & 3c_2              & 4c_1               & \cdots & (p-1)c_6\varphi^p     \\
  \vdots& \vdots  & \vdots            & \vdots            & \vdots             & \ddots & \vdots                \\
  0     & c_{p-3} & 2c_{p-4}          & 3c_{p-5}          & 4c_{p-6}           & \cdots & (p-1)c_{p-1}\varphi^p \\
  0     & c_{p-2} & 2c_{p-3}          & 3c_{p-4}          & 4c_{p-5}           & \cdots & (p-1)c_0              \\
  0     & c_{p-1} & 2c_{p-2}          & 3c_{p-3}          & 4c_{p-4}           & \cdots & (p-1)c_1              \\
  \endmatrix
  \right].
  $$
It suffice for proving 18.3.1 to show the following.    
\vskip 3pt
\fp
\lbox{71}{360}
{{\bf Lemma 18.3.2.} \ \ \it
Under the situation above, the matrix of size  $(p-1)\times (p-1)$  
given by removing the first row and the first column of 
  $$
  T^{p-1}\mod{p}
  $$ 
is a scalar matrix.   
}
\vskip 3pt
\fp
The following proof is due to Hiroshi Suzuki.  
We denote for a square matrix  $M$  by  $M_{i_1,i_2,\cdots,i_j}$  
the matrix consists of the $i_1$, $i_2$, $\cdots$, $i_j$-th rows  
with the  $i_1$, $i_2$, $\cdots$, $i_j$-th columns of  $M$, 
and  $|\hskip 7pt|$  means taking the determinant.  
It is well-known that the characteristic polynomial of  $M$  
is given by 
  $$
  |tI-M|
  =t^p+\sum_{j=1}^p(-1)^{p-j}
  \bigg(\sum_{1\leqq i_1<i_2<\cdots <i_j\leqq p}|M_{i_1,i_2,\cdots,i_j}|\bigg)
  t^{p-j}.
  \tag{18.3.3}
  $$
To save complication, we regard the suffix of a row or column by modulo  $p$.   
Let
  $$
  \align
  S&=
  \left[\matrix
  0 &   &   &        &      \\
    & 1 &   &        &      \\
    &   & 2 &        &      \\
    &   &   & \ddots &      \\
    &   &   &        &  p-1 \\
  \endmatrix\right], \\
  P=&
  \left[
  \matrix
  c_1             & c_0     & c_{p-1}\varphi^p & c_{p-2}\varphi^p & c_{p-3}\varphi^p  & \cdots & c_2\varphi^p     \\
  c_2             & c_1     & c_0              & c_{p-1}\varphi^p & c_{p-2}\varphi^p  & \cdots & c_3\varphi^p     \\
  c_3             & c_2     & c_1              & c_0              & c_{p-1}\varphi^p  & \cdots & c_4\varphi^p     \\
  c_4             & c_3     & c_2              & c_1              & c_0               & \cdots & c_5\varphi^p     \\
  c_5             & c_4     & c_3              & c_2              & c_1               & \cdots & c_6\varphi^p     \\
  \vdots          & \vdots  & \vdots            & \vdots            & \vdots             & \ddots & \vdots                \\
  c_{p-2}         & c_{p-3} & c_{p-4}          & c_{p-5}          & c_{p-6}           & \cdots & c_{p-1}\varphi^p \\
  c_{p-1}         & c_{p-2} & c_{p-3}          & c_{p-4}          & c_{p-5}           & \cdots & c_0              \\
  c_0\varphi^{-p} & c_{p-1} & c_{p-2}          & c_{p-3}          & c_{p-4}           & \cdots & c_1              \\
  \endmatrix
  \right].
  \endalign
  $$
So that  $T=PS$.   
Then we can easily see, for  $1\leqq i_1 < i_2 < \cdots < i_j\leqq p$  that
  $$
  |P_{i_1+1,i_2+1,\cdots, i_j+1}|=|P_{i_1,i_2,\cdots, i_j}|.
  \tag{18.3.4}
  $$
Moreover, we see that
  $$
  |T_{i_1,i_2,\cdots, i_j}|=(i_1-1)(i_2-1)\cdots(i_j-1)|P_{i_1,i_2,\cdots, i_j}|.
  \tag{18.3.5}
  $$
Hence
  $$
  \aligned
  |&tI-T| \\
  &=t^p+
  \sum_{j=1}^{p-1}
  \bigg((-1)^{p-j}\hskip -9pt\sum_{(i_{1},\cdots,i_{j})\ \ \ }
  \hskip -16pt\sum_{\ \ b\in\bold F_p}(b+i_{1})\cdots(b+i_{j})\bigg)
  \hskip -1pt|P_{i_1,\cdots,i_j}|\hskip 3pt t^j+|T|.
  \endaligned
  \tag{18.3.6}
  $$
Here  $|T|=0$.  
The inner sum across over  $(i_1,\cdots,i_j)$  
with  $1\leqq k_1 < \cdots < k_j \leqq p$  modulo 
the equivalence relation such that additions of the same number 
for all the entries of  $(k_1, \cdots, k_j)$  is not change its class.  
Of course, we are assuming that, if an entry exceed  $p$, 
we subtract  $p$  from it and move it to the first entry circularly.  
Since  $p$  is a prime, each lengths are  $1$  or  $p$.   
The length is  $1$  if and only if  $j=p$.  
Thus it suffice for our aim only to prove
  $$
  \sum_{b\in\bold F_p}(b+i_{1})\cdots(b+i_{j})=\cases
  &0 \ \ \ \ \ (\text{if  $1\leqq j < p-1$\ )},\\
  &-1 \ \ \  (\text{if  $j=p-1$\ )}.
  \endcases
  \tag{18.3.7}
  $$
This is shown by considering the fundamental symmetric polynomial 
of  $p$  variables over  $\bold F_p$.  
If  $m=p-1$, then the orbit is unique and we can chose 
as its representative  $(2, 3, \cdots, p)$.   
Therefore we conclude
  $$
  |tI-T|=t^p-|P_{2,3,\cdots,p}|t.
  \tag{18.3.8}
  $$
The Cayley-Hamilton theorem states  $T^p=|P_{2,3,\cdots,p}|T$.  
Since  $|T_{2,3,\cdots,p}|\neq 0$, the proof of (18.3.2) have completed.  
\qed

\newpage \fp
{\bf 18.4. Problems (or exercises?).}\ \ \ 
Problems not yet proved or things to be
investigated which the author have are as follows:
\vskip 5pt
\fp
(1) Looking at the numerical examples, the signature of  $C_n$  or  $D_n$  
are alternative or stable as in the cases of the Bernoulli and Hurwitz numbers.  
Prove this phenomena. 
\vskip 5pt
\fp
(2) For example, about the curve  $y^2=x^5-1$, if  $p\equiv 1\mod 5$  
and  $p-1\not|\ 10m$  then  $p\not|\ C_{10m}$.  
\vskip 5pt
\fp
(3) How are the orders of $2$-part of the numbers  $C_n$  and  $D_n$?   
For the case of the Bernoulli and Hurwitz numbers, they are exactly determined.  
For instance, is it true that, for  $y^3=x^5-1$ ($g=4$), $\ord_2 C_{15n}=12n-3$?

\newpage \input refs.o
\enddocument
\bye